\documentclass{article}
\usepackage{mathtools}
\usepackage{mathrsfs,enumerate}
\usepackage{amsmath}
\usepackage{amssymb}
\usepackage[numbers]{natbib}
\usepackage{hyperref}
\usepackage{dsfont}
\usepackage[draft]{todonotes}
\usepackage[a4paper,top=2cm,bottom=2cm,left=2cm,right=2cm,bindingoffset=5mm]{geometry}

\numberwithin{equation}{section}

\usepackage{color}
\usepackage{dsfont}
\usepackage[latin1]{inputenc}  
\textwidth16cm \numberwithin{equation}{section} \vspace{8cm}
\usepackage{hyperref}

\newtheorem{theorem}{Theorem}[section]
\newtheorem{proposition}[theorem]{Proposition}
\newtheorem{corollary}[theorem]{Corollary}
\newtheorem{lemma}[theorem]{Lemma}
\newtheorem{definition}[theorem]{Definition}
\newtheorem{remark}[theorem]{Remark}

\newcommand{\R}{{\mathbb R}}

\def\proof{{\bf Proof}.\,\,}
\newcommand{\cL}{{\mathcal L}}
\newcommand{\cM}{{\mathcal M}}
\newcommand{\cP}{{\mathcal P}}
\newcommand{\cI}{{\mathcal I}}
\def\cT{{\mathcal T}}
\newcommand{\vep}{\varepsilon}

\newcommand{\vfi}{\varphi}

\def\lg{\langle}
\def\rg{\rangle}
\def\lgx{\lg x\rg}

\def\limk{L^\infty(\lg x\rg^{-k})}
\def\L1k{L^1(\lg x\rg^{k})}

\newcommand{\be}{\begin{equation}}
\newcommand{\ee}{\end{equation}}
\newcommand{\rife}[1]{{(\ref{#1})}}

\newcommand{\dive}{{\rm div}}

\newcommand{\into}{{\int_{\Omega}}}

\newcommand{\intd}{{\int_{\R^d}}}

\newcommand{\la}{\lambda}
\newcommand{\al}{\alpha}

\def\m{\noalign{\medskip}}
\def\de{\delta}

\def\vep{\varepsilon}

\def\vfi{\varphi}
\def\qed{{\unskip\nobreak\hfil\penalty50
          \hskip2em\hbox{}\nobreak\hfil\mbox{\rule{1ex}{1ex} \qquad}
   \parfillskip=0pt
   \finalhyphendemerits=0\par\medskip}}
 \def\la{\lambda}
 \def\ga{\gamma}
\def\rife#1{(\ref{#1})}

\title{Decay rates of convergence for Fokker-Planck equations with confining drift}

\author{
Alessio Porretta\thanks{Department of Mathematics, University of Rome Tor Vergata. 
Via della Ricerca Scientifica 1, 00133 Roma, Italy.    Email: \texttt{porretta@mat.uniroma2.it}. The author is member of GNAMPA research group of Indam. "Partially supported by  the MIUR Excellence Department Project  Math@TOV" awarded to the Department of Mathematics, University of Rome Tor Vergata.} 
 }

\begin{document}

\maketitle
\begin{abstract} We consider Fokker-Planck equations in the whole Euclidean space, driven by Levy processes, under the action of confining drifts, as in the classical Ornstein-Ulhenbeck model. We introduce a new PDE method to get exponential or sub-exponential decay rates, as time goes to infinity, of zero average  solutions, under   some diffusivity condition on the Levy process, which includes the fractional Laplace operator as a model example. Our approach  relies on the long time oscillation estimates of the adjoint problem and applies to (the possible superposition of) both local and nonlocal diffusions, as well as to  strongly or weakly confining drifts. 
Our results extend, with a unifying perspective,  many previous works based on different analytic or probabilistic methods, with several interesting connections.
On one hand, we make a link  between the (nonlinear) PDE methods used for the long time behavior of Hamilton-Jacobi equations and the decay estimates of Fokker-Planck equations;  on another hand,  we give a purely  analytical approach towards some oscillation decay estimates which were obtained so far only with probabilistic coupling methods. 
\end{abstract}

%\vskip1em
% 
\vskip1em

%\tableofcontents

% \vskip1em

\section{Introduction}

This  article is concerned with Fokker-Planck equations in $\R^d$ and with the issue of time-decay rates of zero average solutions, which results as effect of the presence of a confining drift. To fix the ideas, the model problem that we will discuss is the following
\be\label{model}
\begin{cases}
m_t + \left(-\Delta\right)^{\frac\sigma2} m -\dive( b(t,x) m)     =  0 &  \text{ in } ( 0,T ) \times \R^d, \\ 
m( 0,x ) = m_0 (x),\qquad  \,&  \text{ in }   \R^d
  \end{cases}
\ee
where 
$\left(-\Delta\right)^{\frac\sigma2} m$ is the usual Laplace operator for $\sigma=2$ or the fractional Laplacian for $\sigma\in (0,2)$. 
%The purpose of the article is to present a new approach, entirely analytic and relatively simple, which leads to optimal decay rates of zero average solutions and applies to a variety of different regimes, ranging from local to nonlocal diffusions and including strongly or weakly confining drifts.

If $b(t,x)=b(x)$,  the decay of zero average solutions is equivalent, in particular, to the problem of long time convergence of unit mass solutions towards stationary invariant measures. This  is  a fundamental question in the study of Fokker-Planck equations, which attracted a huge interest all over the years. As is well-known, the confining property of the drift plays a crucial role in this convergence, being related to the ergodicity property of the underlying  stochastic process.
The typical result which is  proved in this context,  for $m_0$ with zero average,  is that 
\be\label{mX}
\|m(t)\|_X \leq \varpi(t) \|m_0\|_X \,,\qquad \forall t>0\,,\quad \hbox{with $\varpi(t)\mathop{\to}\limits^{t\to \infty} 0$,}
\ee
where $\varpi(t)$ decays exponentially if the confinement force is sufficiently strong ($b\cdot x \gtrsim |x|^2$ at infinity),   or  it possibly decays  sub-exponentially if $b\cdot x \gtrsim |x|^\gamma$ with $\gamma\in (0,2)$. The space $X$ is typically a $L^p$ weighted space, with $p=1$  being the most natural choice, in a way that
\rife{mX} implies decay of the total mass.
%the weight  may be linked to stationary measures  or to the Lyapunov functions for the process, depending on the approach which is used to establish such an estimate. 

It is impossible here to mention all relevant contributions in the literature, where an estimate such as \rife{mX} has been proved. Let us just  point out two main  approaches that were  developed so far for this kind of analysis. 

Starting from the case of Laplace operator and drift  $b= \nabla V(x)$, a first approach relies on functional analytic methods, exploiting spectral gap results in $L^2$ (weighted) spaces, logarithmic Sobolev (and Poincar\'e) inequalities and, more generally, entropy-energy methods. We refer the reader to many influential papers such as   \cite{Arnold+al},  \cite{Gross}, \cite{Liggett}, \cite{rock}, \cite{tovi}. Further  refinements of spectral gap methods, relying on  ideas in semigroup theory, were  developed more recently, see e.g. \cite{Gualdani+al}, \cite{Kavian+al}. Even if the above results are mostly related to local operators, 
methods based on entropy decay or spectral gaps were also used for nonlocal diffusions, see \cite{BK}, \cite{GI}, \cite{tristani}.

A second approach, which is different in spirit,  stands on ideas and probabilistic methods  developed in the study of recurrence and decay properties of Markov processes.
This is often referred to as the Meyn \& Tweedie approach (see e.g. \cite{MeynTweed}, \cite{tweed}) and this is intrinsically related to Harris' ergodic theorem on Markov chains, as presented in the enlightening version of Hairer and Mattingly \cite{HM}.  The main  tools in  this approach are the existence of a Lyapunov function at infinity and the 
(local) strict positivity of the operator, see e.g.  \cite{Douc+al}, \cite{CaMi}. Extensions to nonlocal diffusions were recently developed in  \cite{Lafleche}. 
We also address the reader to the  influential paper \cite{Bakry+al} for a comparison between the different approaches mentioned above. 
Finally, it is to mention that  similar estimates can be also deduced  from bilateral bounds on the fundamental solution  (see e.g.  \cite{MM}, \cite{O}, \cite{Zh} for divergence free drift terms) or, with a probabilistic approach, from the study of decay rates of the transition probabilities associated to the process, see e.g.  \cite{Eberle}, \cite{Eberle+al}, \cite{mateusz}. We will comment more about that, in what follows.
\vskip1em
%The purpose of this article is to present a new approach, entirely analytic and relatively simple, which leads to optimal decay rates of zero average solutions and applies to a variety of different regimes, ranging from local to nonlocal diffusions and including strongly or weakly confining drifts.

The purpose of this article is to introduce  a new method for proving decay estimates such as \rife{mX}, for solutions of \rife{model}. This is an entirely PDE method, relatively simple, which applies to a variety of different regimes, ranging from local to nonlocal diffusions and including both geometrical and subgeometrical rates of convergence. 
%seems to be very flexible to cover both local and nonlocal diffusions, both geometrical and subgeometrical rates of convergence. 

The motivation of our work comes from the study of the long time behavior of mean-field game systems \cite{EJP}; those are systems of PDEs coupling Fokker-Planck and Hamilton-Jacobi equations, which describe Nash equilibria in large populations of interacting agents (see \cite{LL}, \cite{HCM}, \cite{CP}). In that context, the Fokker-Planck equation is driven by the optimal strategies chosen by the agents, which results in the fact that the drift term in \rife{model} depends on the solution of the coupled Hamilton-Jacobi equation.  With this in mind, the goal  is to obtain decay estimates for Fokker-Planck equations under hopefully cheap conditions  on the drift term $b(t,x)$; in particular, our  results apply to time-dependent fields $b$ and avoid  any direct assumption on $\dive(b)$. 

Before explaining the novelty of our approach, let us give a sample of the kind of results that we prove, when specialized to the model problem \rife{model}. We start with the case of purely local diffusion, where the results are slightly more general.

\begin{theorem}\label{localthm}  Let us  set $\lg x\rg:= \sqrt{1+ |x|^2}$. For $k>0$,  assume that $m_0\in L^1(\R^d, \lgx^k\,dx)$ be such that $\intd m_0\,dx=0$.  Let  $m$ be  a weak solution to the Cauchy problem
\be\label{lappa}
\begin{cases}
m_t -\Delta m -\dive( b(t,x) m)     =  0 &  \text{ in } ( 0,T ) \times \R^d, \\ 
m( 0,x ) = m_0 (x),\qquad  \,&  \text{ in }   \R^d
  \end{cases}
\ee
  where we assume that $b(t,x)$ is a locally bounded function satisfying 
 \be\label{bdiss0}
\exists \,\, \alpha\,,\, R>0\,:\,\quad b(t,x)\cdot x \geq  \alpha \, |x |^{\gamma}   \qquad \forall x  \in \R^d:\, |x|\geq R\,,\, \forall t>0\,,
\ee
and
\be\label{disso0}
\exists \,\, c_0>0\,:\,\quad (b(t,x)-b(t,y))\cdot (x-y) \geq  - c_0 |x-y|    \qquad \forall x,y  \in \R^d\,, \,\forall t>0\,.
\ee
Then we have:
  \begin{itemize}
  
\item[(a)] If $\gamma\geq 2$ in \rife{bdiss0}, then 
$$
\| m(t)\|_{L^1(\lg x\rg^k)} \leq K \, e^{-\omega t} \, \|m_0\|_{L^1({\langle x\rangle^k})}
$$
for some constant $K,\omega$ only depending on $\alpha, c_0,k, d$.  

\item[(b)] If $\gamma \in (0, 2)$ in \rife{bdiss0}, then for any $\bar k>k$  such that $m_0\in L^1(\R^d, \lgx^{\bar k}\,dx)$ we have 
$$
\| m(t)\|_{L^1(\lg x\rg^k)} \leq K \,  (1+t)^{-q}  \, \|m_0\|_{L^1({\langle x\rangle^{\bar k}})}\qquad \hbox{where $q=  \frac{\bar k-k}{2-\gamma}$.}
$$
for some $K$ depending on $\alpha, \gamma, c_0,k, \bar k, d$.  
\end{itemize}
\end{theorem}

Theorem \ref{localthm} recovers, with a unifying approach,  the results which are known in the literature, including the case of slow, or degenerate, confining drift, which was studied previously both with the entropy method  by Toscani and Villani \cite{tovi} and with   more sophisticated tools of semigroup  theory by Kavian, Mischler and Dao \cite{Kavian+al}.  We stress however that we avoid any $C^1$ regularity of the drift $b$ and any requirement on $\dive(b)$ or $Db$; the one-sided condition \rife{disso0} plays somehow a similar role, although it appears to be a milder requirement which, for instance,  allows for any (space-time) bounded perturbation of dissipative vector fields. 
We also mention that, at the expense of a mild restriction on the parameter $k$ in the weight function, the previous result holds even replacing \rife{disso0} with more general conditions, e.g. assuming $(b(t,x)-b(t,y))\cdot (x-y) \geq  - c_0 (|x-y| \vee |x-y|^2)$, see Remark \ref{gendrift}. 
   
\vskip1em
Whilst Theorem \ref{localthm} only gives a sample of the results that we prove,  the main outcome of our approach appears in the extension of this kind of estimates to general Levy processes, including possibly nonlocal diffusions. This is the model result that we prove  for the specific example of the fractional Laplacian, to be compared with \cite[Thm 1.4]{Lafleche}.

\begin{theorem}\label{nonlocalthm}  Let $\sigma\in (0,2)$ and  $k\in (0,\sigma)$. Assume that $m_0\in L^1(\R^d, \lgx^k\,dx)$ be such that $\intd m_0\,dx=0$.  Let  $m$  be  a weak solution to the Cauchy problem \rife{model},   where we assume that $b$ is a locally Lipschitz vector field satisfying  \rife{bdiss0} and 
either of the following conditions:
\vskip0.4em
(i)  $\sigma\in (1,2)$ and  \rife{disso0} holds true.
\vskip0.4em
(ii)  $\sigma\in (0,1]$  and   $(b(t,x)-b(t,y))\cdot (x-y) \geq  - c_0 |x-y| (|x-y|\wedge 1)^{1-\sigma+\de} $, for  some $\de>0$, $c_0>0$.
\vskip0.4em
Then we have:
  \begin{itemize}
  
\item[(a)] If $\gamma\geq 2$ in \rife{bdiss0}, then 
$$
\| m(t)\|_{L^1(\lg x\rg^k)} \leq K \, e^{-\omega t} \, \|m_0\|_{L^1({\langle x\rangle^k})}
$$
for some constant $K,\omega$ only depending on $\alpha, \sigma,  \de, c_0, k, d$.  

\item[(b)] If $\gamma \in ( ((2-\sigma)\vee 1), 2)$ in \rife{bdiss0}, then for any $k<\bar k<	\sigma$  such that $m_0\in L^1(\R^d, \lgx^{\bar k}\,dx)$ we have 
$$
\| m(t)\|_{L^1(\lg x\rg^k)} \leq K \,  (1+t)^{-q}  \, \|m_0\|_{L^1({\langle x\rangle^{\bar k}})}\qquad \hbox{where $q=  \frac{\bar k-k}{2-\gamma}$.}
$$
for some $K$ depending on $\alpha, \gamma,  \sigma, \de, c_0, k, \bar k, d$.  
\end{itemize}
\end{theorem}

Let us stress that the local Lipschitz  character of the drift, which is assumed  in the above statement, was only used to make sure that $m(t)\in L^1$.  
More general statements, which extend Theorem \ref{nonlocalthm}, will be given  in Section \ref{FoPl}  for  Levy operators: we refer the reader to Theorem \ref{decaym}, in Section 4. In that case, the decay is  formulated in the total variation norm for measure-valued solutions. Correspondingly, we will adopt a suitable formulation of solution which can be applied to a far more general setting; as an example, one can consider the case of fractional Laplacian with additional second order degenerate diffusions.

\vskip1em
Compared with the existing literature, our strategy is certainly close to (and even inspired from) the probabilistic approach mentioned above, although we develop a  purely PDE method which has an interest in its own. There are two main cornerstones in the approach that we suggest:

(i) the decay estimates for the Fokker-Planck equation \rife{model} are entirely deduced by duality from oscillation decay estimates  obtained on the solutions of the adjoint problem
 \be\label{model-adjoint}
\begin{cases}
u_t + \left(-\Delta\right)^{\frac\sigma2} u+  b(t,x) \cdot Du     =  0 &  \text{ in } ( 0,T ) \times \R^d, \\ 
u( 0,x ) = u_0 (x),\qquad  \,&  \text{ in }   \R^d\,.
  \end{cases}
\ee
(ii) for solutions of \rife{model-adjoint}, we  prove decay estimates  which are dual to estimates of the type \rife{mX}. Namely, we show that
\be\label{osci0}
\begin{split}
 & \qquad [u(t)]_{\langle x\rangle^k}    \leq \varpi(t) \, [u_0]_{\langle x\rangle^k}\qquad \hbox{with $\varpi(t)\mathop{\to}\limits^{t\to \infty} 0$,}
\\
\m
& \hbox{where  $[u]_{\langle x\rangle^k}=  \sup\limits_{x,y\in \R^{d}}\,\, \frac{|u(x)-u(y)|}{\langle x\rangle^k+ \langle y\rangle^k}$.}
\end{split}
\ee
Estimates such as \rife{osci0} are an evidence of the ergodicity of the underlying process. In the PDE perspective, global (in time) bounds for the oscillation of solutions are usually exploited in the study of long time behavior of Hamilton-Jacobi equations and the convergence to the ergodic problem (see e.g. \cite{BS},  \cite{BCCI2}, \cite{BLT}, \cite{CLN}). Those results typically stand on either classical  gradient estimates (through the celebrated Bernstein method) or doubling variables techniques developed in the theory of viscosity solutions, specifically by Ishii and Lions \cite{IL}. This kind of methods was later extended to nonlocal operators in \cite{BCCI}, \cite{Jak+al}. However, quantitative estimates like \rife{osci0} are new in this framework, compared to those earlier works. Otherwise,   estimates similar to \rife{osci0} can be found in the probabilistic literature in the form of  contraction estimates for transition probabilities in Wasserstein's metrics,  see e.g. \cite{Eberle}, \cite{Eberle+al}, \cite{mateusz}, \cite{wang}. In that context,  the  proofs rely on probabilistic tools based on coupling methods (\cite{ChenLi}, \cite{Lin-Rog}). 

{   In our approach,  the  key idea  is to obtain \rife{osci0} by proving that 
\be\label{maxpri}
u(t,x)-u(t,y) \leq \varpi(t) \, [u_0]_{\langle x\rangle^k} \big[ \langle x\rangle^k+ \langle y\rangle^k   + \psi(|x-y|)\big]
\ee
for some bounded function $\psi(r)$ which is locally H\"older near $r=0$.
Notice that this estimate combines the effects of the Lyapunov function at infinity (used to control long-range interactions) with the local ellipticity used to control short-range interactions, when $|x-y|$ is small.  This is the role of the function $\psi$, namely to handle the small oscillations; at the same time, since $\psi$ is bounded, \rife{maxpri} readily implies \rife{osci0} up to multiplying  the decay rate $\varpi(t)$ by a bounded factor. Finally, by a simple duality argument, an estimate like \rife{osci0} yields the decay $\| m(t)\|_{L^1(\lg x\rg^k)} \leq \varpi(t)\, \|m_0\|_{L^1({\langle x\rangle^{\bar k}})}$ for the solution of \rife{model}.

It seems remarkable that  {\it the decay of Fokker-Planck equations  will be deduced from the single estimate \rife{maxpri}} proved on the dual problem \rife{model-adjoint}. This is the main difference of our method compared to the Meyn-Tweedy, or Harris-type approach, where the local ellipticity is exploited directly in Harnack-type properties of  the Fokker-Planck equation (or strong positivity, say irreducibility conditions, on the associated semigroup).  
%Otherwise,   estimates as \rife{osci0} can be found in the probabilistic literature in the form of  contraction estimates for transition probabilities in Wasserstein's metrics,  see e.g. \cite{Eberle}, \cite{Eberle+al}, \cite{mateusz}, \cite{wang}. In that context,  the  proofs rely on probabilistic tools based on coupling methods (\cite{ChenLi}, \cite{Lin-Rog}). 
%From this perspective,  our proof can be seen as a PDE method which is alternative to the probabilistic approach.  
By contrast, \rife{maxpri} is a pointwise estimate that we derive from maximum principle, to be applied  to a suitable operator in the product space for the variables $(x,y)$. This kind of  application  of the maximum principle contains ideas that connect the doubling variables methods of viscosity solutions to probabilistic coupling methods. We will explain more about that in Section  \ref{doubcoub}, where the  key steps are given in Lemmas \ref{locmax} and \ref{new-cone}.}
%
%
%Indeed, if we denote $\mathcal U=u(t,x)-u(t,y)$ and $\mathcal \Psi:= \varpi(t) [u_0]_{\langle x\rangle^k}\, \big[ \langle x\rangle^k+ \langle y\rangle^k   + \psi(|x-y|)\big]$, estimate \rife{maxpri} follows if we build a second order positive operator $\mathcal A$ in the product space $(x,y)$ such that $\mathcal A({\mathcal U}-{\mathcal \Psi} )

Let us mention that  the idea of an analytic version of the coupling method was already suggested in  \cite{PP} for local diffusions (some related ideas  appeared later in \cite{FouPer}, in a rather informal way). On that occasion, we explained that the doubling variables approach developed in the theory of viscosity solutions was parallel, and conceptually equivalent, to the coupling method used in probability (see \cite[Appendix]{PP}). In the present article, we come back on this viewpoint and introduce similar ideas to handle nonlocal diffusions as well, getting  a crucial improvement in the doubling variable method which is used for Levy operators. 
This allows us to give entirely PDE proofs of some  contraction estimates that were only proved with probabilistic methods so far, with the advantage of relaxing the assumptions and making them ready to several possible extensions. 

%To sum up, the contribution of this article is twofold: on one hand, we put in connection the (nonlinear) PDE methods used for the long time behavior of Hamilton-Jacobi equations with the decay estimates of Fokker-Planck equations, on another hand  we suggest a new analytical approach towards oscillation decay estimates which encodes the  contraction estimates obtained with probabilistic coupling methods.  

In the present paper, in order not to increase the level of technicalities which could be induced by generalizations, we confine ourselves to the study of Levy processes, which  contain both local diffusion and jump processes, giving evidence of the generality of the method. 

We also show how the suggested strategy can be adapted to different
Lyapunov functions in order to obtain exponential or sub-exponential rates of convergence.  {\it We refer the reader to  Theorem \ref{decayFP2} for a general version of our results}, which is not limited to power-type Lyapunov functions. In particular, we show how different decay rates  are proved for the weighted norm $\| m(t)\|_{\cM_\vfi}$ (the total variation of $m(t)\,\vfi$) according to whether the Lyapunov function $\vfi$ satisfies stronger, or weaker, conditions of  super-solution.  This way,    Theorem \ref{decayFP2}  will provide with a general recipe which may please the readers inclined towards statements in a more abstract form. As an example of application of similar generalizations, we will obtain exponential convergence of solutions to \rife{lappa} even for slowly confining drifts, up to using stronger Lyapunov weights rather than powers:  we refer to Corollary \ref{kavian}, which extends recent results obtained in \cite{Kavian+al}.
\vskip2em
To conclude, let us briefly mention the organization of the paper. In Section 2, we introduce the notation and the assumptions that will be used, as well as the notions of solutions for the evolution problems to be considered.  Section 3 is devoted to the oscillation decay estimates for problems like \rife{model-adjoint}; this is where we use doubling variables methods (with coupling ideas embedded) to derive estimates as \rife{osci0}, with power-type weight functions. In Section 4, we deduce by duality the decay of Fokker-Planck equations, which generalize Theorem \ref{localthm} and \ref{nonlocalthm}.  Section 5 is devoted to a more general version of the results, where we consider Lyapunov functions of general type, and we give exponential or sub-exponential decay rates accordingly. We also mention the application of our results to the convergence of unit mass solutions towards stationary measures. Finally, in the Appendix, Section 6,  we give the proof  of the existence and uniqueness of solutions  for the Fokker-Planck equation, under  the (non standard) setting of assumptions which is used in the paper.

\section{Notation, assumptions and preliminary tools}

We denote by $\R^d$, $d\geq 1$, the $d$-dimensional Euclidean space, with Lebesgue measure $dx$. $I_d$ denotes the identity matrix in $\R^d$. Given two real numbers $a,b$, we use the standard notation $ a\vee b= \max(a,b)$ and $a\wedge b= \min(a,b)$. For a positive continuous function $\vfi(x)$, we denote by $L^1(\vfi(x))$ the usual Lebesgue space $L^1(\R^d, \mu)$ defined in terms of  the measure $\mu = \vfi(x) dx$, and by $L^\infty(\vfi(x))$ the space of measurable functions $u:\R^d\to \R$ such that $u\, \vfi $ is essentially bounded (i.e. $u\vfi \in L^\infty(\R^d)$).  

We will use systematically the weight function 
$$
\lg x\rg:= \sqrt{1+ |x|^2}\,
$$
replacing the role of $|x|$. Notice that $D\lg x\rg= \frac x{\lg x\rg}$.

The space $\cM(\R^d)$ denotes the space of  finite   (signed) Borel measures $m$ endowed with the total variation norm $\| m\|_{TV}= |m|(\R^d)$; here $|m|= m^+ + m^-$, where $m^\pm$ are positive measures giving the Hahn decomposition of $m$. The set $\cP(\R^d)$ denotes the space of probability measures, where we consider the Kantorovich-Rubinstein distance
$$
d_1(m,\tilde m):= \sup\{\intd \vfi \, d(m-\tilde m)\,,\,\, \vfi\,\hbox{Lipschitz,}\,\, \|\vfi\|_\infty, \|D\vfi\|_\infty\leq 1 \}\,.
$$
Finally, we denote by $\cM_k(\R^d)$ the subset of measures $m$ with finite $k$ moments, i.e. $\lgx^k\in L^1(\R^d, d|m|)$. For $m\in \cM_k$,   we use the notation
$$
\|m\|_{\cM_k}:= \intd \lgx^k\, d|m|
$$ 
to denote the total variation of $\lgx^k \,m$. We use a similar notation with possibly different weights, setting $\cM_\vfi(\R^d)$ the subset of measures such that $\vfi\in L^1(\R^d, d|m|)$ for some positive function $\vfi$, with $\|m\|_{\cM_\vfi}= \intd \vfi\, d|m|$.

For given $T>0$, we set $Q_T:= (0,T)\times \R^d$. We denote by $C^{1,2}(Q_T)$ the functions defined on $Q_T$ which are $C^1$ in $t$ and $C^2$ in $x$. The notation $C_b(U)$ refers to continuous bounded functions on a set $U$, and $C_c(U)$ to continuous compactly supported functions. We denote by $C([0,T];\cM(\R^d)^*)$ the functions $m$  from $[0,T]$ into $\cM(\R^d)$ which are continuous with respect to the weak$^*$ topology of measures. When working with probability measures, we denote $C([0,T];\cP(\R^d) )$ the functions $m$ which are continuous from $[0,T]$ into $\cP(\R^d)$, endowed with the $d_1$ distance.

Associated to  a continuous function $\vfi$ such that $\vfi(x) \geq 1$  for all $x\in \R^d$, 
%. For any   function $u$ such that $u= O(\vfi)$,  we have
%$$
%\sup_{x,y\in \R^{2d}}\,\, \frac{|u(x)-u(y)|}{ \vfi(x)+ \vfi(y)} = \inf_{c\in \R} \,\, \|u+c\|_{L^\infty(\vfi^{-1}dx)}
%$$
we  introduce suitable weighted seminorms, defined as
$$
[u]_\vfi= \sup_{x,y\in \R^{d}}\,\, \frac{|u(x)-u(y)|}{ \vfi(x)+ \vfi(y)}\,.
$$
In particular, we will set
$$
[u]_{\langle x\rangle^k}:=  \sup_{x,y\in \R^{d}}\,\, \frac{|u(x)-u(y)|}{\langle x\rangle^k+ \langle y\rangle^k}\,.
$$
With a slight abuse of notation, when considering  the case $k=0$, we denote
$$
[u]_0:=  \sup_{x,y\in \R^{d}}\,\, \left ( |u(x)-u(y)|\right) \,
$$
the seminorm measuring the oscillation of a bounded function $u$.   
\vskip1em
Given a continuous matrix $\Sigma\,:\R^d\to \R^{d\times \ell}$, $\Sigma^*$ denotes its transposed matrix; then we define the second order operator
\be\label{L0}
\begin{split}
\cL_0(u): & = -\lambda_0 \Delta u - {\rm tr}(\Sigma(x)\Sigma^*(x)D^2 u) 
\\
& = -\lambda_0 \sum_{i=1}^d u_{x_ix_i} - \sum_{i,j=1}^d \sum_{k=1}^\ell\Sigma_{ik}(x)\Sigma_{jk}(x) u_{x_ix_j}\
\end{split}
\ee
where $\lambda_0 \geq 0$. Without going in search of the most generality, we will assume $\Sigma$ to be bounded and Lipschitz continuous:
\be\label{sigma}
\exists \,\, \sigma_0,\sigma_1>0\,:\quad \|\Sigma (x) \| \leq \sigma_0 \,,\qquad \|\Sigma(x)-\Sigma(y)\| \leq \sigma_1 |x-y| \quad \forall x, y\in \R^d\,.
\ee
Now we define a pure jump  Levy operator as 
\be\label{Levy}
\mathcal{I} ( x, [ u ] ) := \int_{\R^d} \{u ( x+z ) -u( x ) -   (Du( x ) \cdot  z ) \mathds{1}_{| z | \leq 1} \}\nu ( dz )
\ee
%where the measure $\nu$ satisfies
%\be\label{nu1}
%\int_{\R^d} ( |z|^2 \wedge 1) \nu(dz) <\infty \,.
%\ee
%We also suppose some kind of ellipticity of the nonlocal operator. To keep things simple, we assume that 
in which the measure $\nu$ is absolutely continuous and satisfies
\be\label{nu2}
\exists \lambda, \Lambda \geq 0\, , \sigma\in (0,2)\,:\, \quad   \, \frac\lambda {|z|^{d+\sigma}} \leq \frac{d\nu}{dz}\leq   \, \frac\Lambda{|z|^{d+\sigma}}
\qquad \forall z \in \R^d\,.
\ee
%As usual, in the case $\sigma>1$ the above definition \rife{Levy} has to be meant in the sense of principal value.
%where $B_1=\{x\in \R^d: |z|\leq 1\}$.  
Typical examples of operators which satisfy  \rife{Levy}-\rife{nu2} are given by the fractional Laplacian $(-\Delta)^{\frac\sigma2}$ ($\sigma\in (0,2)$) or more generally by the generators of $\alpha$-stable jump processes. 
We define  the operator $\cL$ as
\be\label{L}
\cL[u]:= \cL_0(u) - \mathcal{I} ( x, [ u ] )\,.
\ee
It is well known (see e.g. \cite{Ap}) that $\cL[u]+ b\cdot Du$ is the infinitesimal generator of a general Levy process where the drift-diffusion part is generated by $\cL_0+ b\cdot D$ and the jump part is represented by $\cI$. This is associated with the stochastic SDE
\be\label{SDE}
\begin{cases} dX_t=-b(X_t) dt +  \sqrt{2\lambda_0} dB_t + \sqrt 2\Sigma(X_t) dB_t + dP_t & \\
X_0= x\in \R^d & 
\end{cases}
\ee
where $B_t$ is a standard Brownian motion and $P_t$ is a compound Poisson process with measure $\nu$. 

We point out that  more general conditions upon the diffusion matrix $\Sigma(x)$ and the Levy measure $\nu(z)$ would be possible, but we keep things simpler here for the reader's convenience. 
%See also Remark \ref{} in Section xxx.
 \vskip1em 
Under the assumptions \rife{sigma}, \rife{nu2}, and for a continuous vector field $b(t,x)$,  we consider the evolution problem
\be\label{pbgen}
\begin{cases}
\partial_t u + \cL[u]   + b(t,x) \cdot Du     =  0 &  \text{ in } ( 0,T ) \times \R^d, \\ 
u ( 0,x ) = u_0 (x),\qquad  \,&  \text{ in }   \R^d
\end{cases}
\ee
{ It is well known that, if $b$ is sufficiently regular, then the unique solution of \rife{pbgen} is characterized as
$u(t,x)= {\mathbb E} (X_t)$, where $X_t$ is the solution of \rife{SDE}.} 
\vskip0.4em
The associated Fokker-Planck equation is defined through  the adjoint problem of \rife{pbgen}. To this purpose, we set
\be\label{Levy*}
\mathcal{I}^* ( x, [ m ] ) := \int_{\R^d} \{m ( x+z ) -m( x ) -   (Dm( x ) \cdot  z ) \mathds{1}_{| z | \leq 1} \}\nu^* (  dz )
\ee
 where we define $\nu^*(B):=\nu(-B)$, for every Borel set $B\subseteq \R^d$. Then the operator $\cL^*$ is defined as
 \be\label{L*}
 \begin{split}
 \cL^*[m]: & = \cL_0^*(m) - \mathcal{I}^* ( x, [ m ] )
 \\
 & = - \lambda_0 \Delta m - \sum_{i,j=1}^d\sum_{k=1}^\ell  \left( \Sigma_{ik}\Sigma_{jk} m\right)_{x_ix_j} - \int_{\R^d} \{m ( x+z ) -m( x ) -   (Dm( x ) \cdot  z ) \mathds{1}_{| z | \leq 1} \}\nu^* (  dz )\,.
 \end{split}
 \ee
This allows us to introduce the Fokker-Planck equation 
\be\label{FP-gen} 
\begin{cases}
\partial_t m + \cL^*[m]   -\dive(  b(t,x) \,m)    =  0 &  \text{ in } ( 0,T ) \times \R^d, \\ 
m ( 0) = m_0 ,\qquad  \,&  \text{ in }   \R^d
\end{cases}
 \ee
where $m_0\in \cM(\R^d)$. 
\vskip1em
Suitable formulations of problems \rife{pbgen} and \rife{FP-gen} are given below.
In order to ease notations, in the following we set, for $\vfi\in C^2(\R^d)$,
\be\label{notLb}
 \cL^b[\vfi]:= \cL[\vfi]   + b(t,x) \cdot D\vfi \qquad\,; \qquad  \qquad \cL_0^b[\vfi]:= \cL_0[\vfi]   + b(t,x) \cdot D\vfi \,.
\ee
%In the results we are going to prove, 
 %results, We assume that $b: (0,T) \times \R^d \to \R^d$ is a continuous vector field. 
%
%Here we consider the case of local diffusions, and specifically the model case given by 
%\be\label{brown}
%\begin{cases}
%\partial_t u -  \Delta u  + b(t,x) \cdot Du     =  0 &  \text{ in } ( 0,T ) \times \R^d, \\ 
%u ( 0,x ) = u_0 (x),\qquad  \,&  \text{ in }   \R^d
%  \end{cases}
%\ee
%where $b: (0,T) \times \R^d \to \R^d$ is a continuous vector field and $u_0$ a continuous function in $\R^d$.
%\vskip1em
The key-assumption, throughout the whole paper,  is the following  confining condition on the drift:
\be\label{bdiss}
\exists \,\, \alpha\,,\, R>0\,:\,\quad b(t,x)\cdot x \geq  \alpha \, |x |^{\gamma}   \qquad \forall x  \in \R^d:\, |x|\geq R\,,\, \forall t>0\,,
\ee
where we assume that $\gamma >0$. 
Three different regimes will be considered: $\gamma \geq 2$, which is the coercive
case, $\gamma \in [1,2)$, which is often referred to in the literature as  a slowly confining case, and $\gamma\in (0,1)$, which may be called a degenerate confining case.
Our goal is to give a unifying approach which provides with  optimal decay estimates in all such  ranges of $\gamma$. 

We also assume a one-sided control for the accretivity of $b$. Two slightly different conditions will be assumed depending on the ellipticity of the nonlocal operator, say whether $\sigma$ in \rife{nu2} is bigger or smaller than $1$. Namely, we suppose
%\footnote{valutare estensioni tipo $(b(t,x)-b(t,y))\cdot (x-y) \geq  - c_0(1+ ( |x| \wedge 1)^{\theta}+( |y| \wedge 1)^{\theta})  |x-y| $ o simili...} 
that $b$ satisfies either 
\be\label{disso}
\exists \,\, c_0>0\,:\,\quad (b(t,x)-b(t,y))\cdot (x-y) \geq  - c_0 |x-y|    \qquad \forall x,y  \in \R^d\,, \,\forall t>0\,,
\ee
or the strongest condition (when $\sigma\leq 1$):
 \be\label{disso2}
\sigma\in (0,1]\,\, \hbox{and $\exists \,\de \in (0,1)$: \quad $(b(t,x)-b(t,y))\cdot (x-y) \geq  - c_0 |x-y| (|x-y|\wedge 1)^{1-\sigma+\de} $,}
\ee
 for  some $c_0>0$. Let us mention that further generalizations of \rife{disso}, \rife{disso2} are also possible; we refer to Remark \ref{gendrift}.

\subsection{Viscosity solutions}

It is convenient to use the framework of viscosity solutions for problem \rife{pbgen}. In our setting, where the operator is defined through \rife{L},  we could use slightly different formulations of the notion of viscosity solution, which all turn out to be  equivalent; we refer the reader to \cite{BaIm} for a discussion of this issue, and we pick just one among those possible versions. Moreover, since our focus is on regularity estimates, we will simplify some requirement by only considering continuous  solutions. For the interested reader, upper semi-continuous (USC) subsolutions, and lower semi-continuous (LSC) supersolutions, are defined accordingly in \cite{BaIm}.

\begin{definition}\label{def-visc} A continuous function $u\in C^0([0,T]\times \R^d)$ is a viscosity solution of \rife{pbgen} if $u\in L^\infty((0,T); L^\infty(\langle x\rangle^{-m}))$ for some $m\in (0,\sigma)$,   if  $u(0,x)=u_0(x)$ in $\R^d$ and if $u$ satisfies the following two requirements:
\begin{itemize}

\item[(i)] for every $\vfi\in C^{1,2}(Q_T)$ such that  $\vfi\in L^\infty((0,T); L^\infty(\langle x\rangle^{-k}))$ for some  $k\in (0,\sigma)$, if $(t_0,x_0)\in Q_T$ is a global maximum point of $u-\vfi$, then
$$
\partial_t \vfi(t_0, x_0) + \cL[\vfi](t_0, x_0)   + b(t_0,x_0) \cdot D\vfi(t_0, x_0)     \leq   0\,.
$$
\item[(ii)] for every $\vfi\in C^{1,2}(Q_T)$ such that  $\vfi\in L^\infty((0,T); L^\infty(\langle x\rangle^{-k}))$ for some  $k\in (0,\sigma)$, if $(t_0,x_0)\in Q_T$ is a global minimum point of $u-\vfi$, then
$$
\partial_t \vfi(t_0, x_0) + \cL[\vfi](t_0, x_0)   + b(t_0,x_0) \cdot D\vfi(t_0, x_0)     \geq   0\,.
$$
\end{itemize} 
\end{definition}

\begin{remark} We stress that the requirement made in Definition \ref{def-visc} concerning the growth at infinity, for  both the solution and the test function, is only due to the nonlocal term $\cI$. Due to assumption \rife{nu2}, some limitation in the growth at infinity is required in order that $\cI(x,[u])$, $\cI(x,[\vfi])$ be well-defined. Alternative definitions, formulated with a more local character,  are suggested in \cite{BaIm}, which turn out to be equivalent.  

Another alternative definition of solution can be given using sub-jets and super-jets of the function $u$, as introduced in \cite{CIL},   denoted by $ { J}^{2,-}_{Q_T} u $, $ { J}^{2,+}_{Q_T} u $ respectively. The notation   ${\overline J}^{2,\pm}_{Q_T} u$ denotes the closure of elements in  the sub/super-jets.
\end{remark}

The following lemma is the extension to nonlocal operators of the classical Ishii-Jensen's lemma in viscosity solutions's theory \cite[Thm 8.3]{CIL}. The following nonlocal version follows from  \cite[Corollary 2]{BaIm}; we use here a slightly more readable version, which is  obtained asymptotically, see e.g.  \cite[Proposition 2]{BaIm}.

\begin{theorem}\label{test} Let  $u$ be an USC sub-solution of  \rife{pbgen},  let  $v$ be a  LSC super-solution of \rife{pbgen}, and let 
$(\hat t, \hat x, \hat y) \in   (0,T) \times (\R^d)^2$ be a local maximum point of the map $(t,x,y)\mapsto [u(t,x)-v(t,y)-\Phi(t,x,y)]$, for some function $\Phi\in C^{1,2}((0,T)\times (\R^d)^2)$.  

Then, for every $n$ sufficiently large, there exist real numbers $a,b\in \R$ and matrices  $X, Y \in { \mathcal S}_d $ such that 
\begin{align*}
a-b = \partial_t \Phi(\hat t, \hat x, \hat y) & \,,\,\,  (a, D_{x} \Phi(\hat t, \hat x ,\hat y), X) 
\in {\overline J}^{2,+}_{Q_T} u (\hat t, \hat x), \,\, 
(b, -D_{y} \Phi(\hat t, \hat x ,\hat y), Y) 
\in {\overline J}^{2,-}_{Q_T} v (\hat t, \hat y)
\\ & 
   -  (n+ c_d\| A\|)I \leq \begin{pmatrix}   X    & 0 \\    0 &  -Y  \end{pmatrix} 
 \leq A+
  \frac1n A^2 
\end{align*}
where $A=D^2_{x,y}  \Phi(\hat t, \hat x ,\hat y)$. Moreover, we have
\be\label{abab}
\begin{split}
& a  - {\rm tr}( Q(\hat x) X) + b(\hat t, \hat x)\cdot D_{x} \Phi(\hat t, \hat x ,\hat y) - {\mathcal I}(\hat x, u(\hat t), D_{x} \Phi(\hat t, \hat x ,\hat y)) \leq 0  \\
& b- {\rm tr}( Q(\hat y) Y) - b(\hat t, \hat y)\cdot D_{y} \Phi(\hat t, \hat x ,\hat y) - {\mathcal I}(\hat y,v(\hat t), -D_{y} \Phi(\hat t, \hat x ,\hat y))  \geq 0 
\end{split}
\ee
where we have denoted, for $x,p\in \R^d$:
\be\label{notp}
\begin{split}
& Q(x)= \lambda_0 I_d+ \Sigma\Sigma^*(x) \\
& {\mathcal I}( x, u, p)= \intd [u(x+z)-u(x)-(p\cdot z) \mathds{1}_{{| z | \leq 1}}]d\nu(z)\,.
\end{split}
\ee
\end{theorem}

Let us stress that the previous theorem also implies the consistency of the viscosity formulation;  if $u\in C^{1,2}(Q_T)$ is a viscosity solution, then it is also a classical solution. Otherwise,  if $u_0$ is a continuous function with polynomial growth, problem \rife{pbgen} admits a unique viscosity solution; we refer to Proposition \ref{compa} for a statement of this kind. 

{   \begin{remark}  The reader which is not familiar with viscosity solutions should not be puzzled by the  statement of Theorem \ref{test}.  In fact, if  $u, v$ were smooth, the above statement would be straightforward: if $u(t,x)-v(t,y)-\Phi(t,x,y)$ has a local maximum at $(\hat t, \hat x, \hat y)$, then by elementary calculus one has
$$
\partial_t u- \partial_t v = \partial_t \Phi \,, \quad \begin{cases} D_xu= D_x \Phi & \\ D_y v= -D_y \Phi & \end{cases} \,,\quad   \begin{pmatrix}   D^2_xu    &  0 \\    0 &  -D^2_y v  \end{pmatrix} \leq \ D^2_{(x,y)} \Phi
$$
where everything is computed on $(\hat t, \hat x, \hat y)$. In this case  the statement obviously holds (and can be read  for $n\to \infty$) with $a= \partial_t u$, $b = \partial_t v$, $X= D^2_xu $ and $Y=D^2_y v$, where the two inequalities \rife{abab} are just the conditions $[\partial_t u+ \cL[u] + b \cdot Du]_{|(\hat t, \hat x)} \leq 0$   (respectively, $[\partial_t v+ \cL[v] + b \cdot Dv]_{|(\hat t, \hat y)}\geq 0$) having replaced  $D_xu(\hat t, \hat x), D_yv(\hat t, \hat y)$ by $D_x\Phi$, $-D_y \Phi $ respectively. 

The relevance  of Theorem \ref{test} is to guarantee that a similar conclusion holds for functions $u,v$ that are not differentiable (even possibly just semi-continuous); in that case, $X,Y$ should be considered as approximations (as $n\to \infty$) of the second derivatives $D^2_xu , D^2_yv$, and similarly $a,b$ as replacements of $  \partial_t u\,, \partial_t v$. Let us mention here that $X,Y$ in principle depend on $n$, but a standard compactness argument (see e.g. \cite[Rmk 2.3]{PP}) allows eventually to let $n\to \infty$ and to obtain new matrices $X,Y$ that satisfy the upper bound $\begin{pmatrix}   X    & 0 \\    0 &  -Y  \end{pmatrix} 
 \leq A$ in the limit. Thus, there is no loss in comprehension for a reader that  would consider $n\to \infty$ in a rough reading of the above statement.  
 \end{remark}
}

\subsection{Duality solutions of Fokker-Planck equations}

Solutions of the Fokker-Planck equation \rife{FP-gen} will be  defined  here in duality with the notion of viscosity solution for $\cL^b$. 
 
\begin{definition}\label{def-dual}
Let $m_0\in \cM(\R^d)$. A function $m\in C^0([0,T);\cM(\R^d)^*)$ is a solution to \rife{FP-gen} if
\begin{align*}
&
 \intd \xi\, dm(t)+ \int_0^t\intd f \,dm(\tau) \,d\tau = \intd \vfi(0,x)\, dm_0 
 \\
 & \qquad \qquad\forall \, t\in (0,T)\,,\, \xi\in C_b(\R^d), \, f\in C_b(Q_t)\,, \vfi\in C_b([0,t] \times \R^d)\quad \mbox{such that} \\
 &  \qquad\qquad\qquad \hbox{$\vfi$   is  a viscosity solution of }\quad  \mbox{\bf (B)}\,\,\begin{cases}
-\partial_t \vfi+ \cL^b[\vfi]     =  f &  \text{ in } ( 0,t ) \times \R^d, \\ 
\vfi(t,x) = \xi (x),\qquad  \,&  \text{ in }   \R^d
\end{cases} 
\end{align*}
 \end{definition}
 
 Let us point out that, since classical solutions are viscosity solutions, then  a solution $m$ according to Definition \ref{def-dual} is also a weak solution
 in the usual distributional sense, which means
 $$
 \int_0^T\intd [-\partial_t \vfi + \cL^b[\vfi] ] \,dm(t)\,dt = \intd \vfi(0,x)\, dm_0  \qquad\forall \vfi\in C^{1,2}_c([0,T)\times \R^d)\,.
 $$
The converse statement, namely that weak solutions enjoy the duality property of Definition \ref{def-dual}, is more delicate under the general setting in which 
we consider the drift term $b$ and the nonlocal diffusion $\cI$. However, this is certainly true for second order operators (i.e. when $\lambda_0>0$ in \rife{L0}) or for locally Lipschitz vector fields $b$.
In fact, weak solutions satisfy Definition \ref{def-dual} whenever they are unique. This is a consequence of the following existence and uniqueness result, which fully 
justifies the previous notion of solution in the present setting.
We recall that $\cM_k(\R^d)$ denotes the subset of measures in $\R^d$ having finite $k$-th moment.

\begin{theorem}\label{exiuniq} 
Let $m_0\in  \cM_k(\R^d)$, for some $0<k<\sigma$.  Let $\cL$ be defined by \rife{L} where $\Sigma$ satisfies  \rife{sigma},  the Levy measure in $\cI$ satisfies \rife{nu2}, and $\lambda_0+\lambda>0$.  Assume that $b\in C^0(Q_T)$ is such that $b(t,x)\cdot x$ is bounded below
and either of the following conditions hold:
\vskip0.4em
(i)  $\la_0>0$ or $\sigma\in (1,2)$,  and  \rife{disso} holds true whenever $|x-y|\leq 1$.
\vskip0.4em
(ii)  $\la_0=0$, $\sigma\in (0,1]$  and   \rife{disso2} holds true whenever $|x-y|\leq 1$.

Then there exists a unique solution $m$ of \rife{FP-gen} in the sense of Definition \ref{def-dual}, and $m(t)\in \cM_k$ for all $t>0$. Moreover, $m$ is obtained as the limit of solutions $m_\vep$ corresponding to the operator $(\cL^b)^*-\vep \Delta$.

Finally,   the application $m_0\mapsto m(t)$ preserves positivity and mass.
\end{theorem}

The proof of Theorem \ref{exiuniq} will be given in the Appendix, for the sake of completeness. It stems from the simple idea that, in duality by what happens for problem (B), the  notion of solution given in Definition \ref{def-dual} is stable by  vanishing viscosity limits. Thus, even if at first sight the formulation of Definition \ref{def-dual} sounds a bit unusual, the reader can see that this is actually a quite natural setting in which problem \rife{FP-gen} can be formulated (and proved to be well-posed) outside a standard set of assumptions.
In fact, the uniqueness result, and the 
characterization of the solution as the limit of Fokker-Planck equations with vanishing viscosity, show the consistency of this formulation with all known results  in terms of weak or classical solutions. 
In particular,  one can readily deduce that $m\in C([0,T];L^1(\R^d))$ when $\la_0>0$ (second order problems), or for the case of fractional Laplacian (or more generally, symmetric kernels) under some extra condition on $\dive (b)$ (see e.g. \cite[Thm 1.2]{Lafleche}).

Finally, we point out that the Definition \ref{def-dual} has been given for signed measures, because we will deal with zero average measures. Due to the Hahn decomposition, and the linearity of the formulation,  every solution $m$ can be split as $m= m_1-m_2$, where $m_1, m_2$ are the solutions corresponding to $m_0^+$, $m_0^-$ respectively. 
%Therefore, the analysis of problem \rife{FP-gen} will be actually reduced to the case of probability measures.

\bigskip

\section{Oscillation decay estimates}

In this Section we give the oscillation estimates for problem \rife{pbgen}, where we assume that $u_0$ is a continuous function in $\R^d$, and 
$\cL, b$ satisfy the assumptions \rife{sigma}, \rife{nu2}, \rife{bdiss} and \rife{disso} (or \rife{disso2}).
\vskip1em

%\subsection{Preliminary tools}

A crucial role here is played by the explicit Lyapunov functions $\lg x\rg^\beta:= (\sqrt{1+ |x|^2})^\beta$.
The properties of those super-solutions are listed below.

\begin{lemma}\label{lyap}  Assume \rife{bdiss}. Let us denote $\vfi(x)= \lg x\rg^\beta$, with $\beta>0$. Then we have:

\begin{itemize}

\item[(i)] for every $\vep>0$ there exists $K_\vep>0$ such that  
\be\label{lyapL0}
\mathcal {L}_0[\vfi] + b(t,x)\cdot D\vfi \geq (\alpha-\vep) \beta \, \frac{\vfi(x)}{\lg x\rg^{2-\gamma}} - K_\vep \qquad \forall x\in \R^d\,, t>0\,.
\ee

\item[(ii)] Let  the Levy measure $\nu$ satisfy  \rife{nu2}. If  $0<\beta<\sigma$,  then for every $\vep>0$ there exists $K_\vep>0$ such that  
\be\label{lyapI}
\mathcal {L}[\vfi] + b(t,x)\cdot D\vfi \geq (\alpha-\vep) \beta \, \frac{\vfi(x)}{\lg x\rg^{2-\gamma}} - K_\vep \qquad \forall x\in \R^d\,,t>0,
\ee
provided $\beta\leq 1$ or $\gamma>1$ in \rife{bdiss}.
\end{itemize}
\end{lemma}

\proof We compute $D\vfi= \beta \lgx^{\beta-1} \frac x{\lgx}$, $D^2 \vfi= \beta\lg x\rg^{\beta-2}\left( I_d+ (\beta-1) \frac x{\lg x\rg}\otimes  \frac x{\lg x\rg}\right)$ and so, if $Q(x)= \lambda_0 I_d+ \Sigma\Sigma^*(x)$, we have from \rife{sigma}
$$
 \left |{\rm tr}(Q(x) D^2 \vfi)  \right| \leq \beta\, c_\beta\, (\sigma_0^2 + \lambda_0) \lg x\rg^{\beta-2}\,
$$
for some $c_\beta$ only depending on $\beta$ and $d$. Hence
$$
\cL_0^b[\vfi] \geq  \beta \lgx^{\beta-2} \left(  b(t,x)\cdot x -c_\beta(\sigma_0^2 + \lambda_0)\right)\,.  
$$
Using \rife{bdiss}, we get \rife{lyapL0}.

As for the Levy  operator $\cI$, the difference is that we need to handle the nonlocal term for large values of $|z|$. Indeed, we have
\begin{align*}
\cI(x,[\vfi]) & =  \int_{|z|>1} \{\vfi ( x+z ) -\vfi( x )  \}\nu ( dz ) +  \int_{|z|\leq 1} \{\vfi ( x+z ) -\vfi( x ) -   (D\vfi( x ) \cdot  z ) \}\nu ( dz )
\\
& \leq    c_\beta\, \int_{|z|>1} \max(  \lg x\rg^{\beta-1} \, |z| \,,\, |z|^\beta)\nu(dz)
+ c_\beta \lg x\rg^{\beta-2} \int_{|z|\leq 1}|z|^2\, \nu ( dz )
\end{align*}
and, due to \rife{nu2}, the first term   is integrable  provided $\beta<\sigma$.  
If $\sigma\leq 1$, this also implies that $\beta\leq 1$, so we deduce that $\cI(x,[\vfi])$ is bounded above. Hence we get
$$
\cL^b[\vfi] \geq  \beta \lgx^{\beta-2}  \, b(t,x)\cdot x - C_\beta   
$$
for some different constant $C_\beta$, possibly depending on $\lambda_0$ and $\sigma_0$ as well. Using \rife{bdiss}, we conclude with \rife{lyapI}. 
The same remains true if $\sigma> 1$ and $\beta\leq 1$. 

%Notice that $\gamma+\beta-2>0$ implies $\gamma>1$ whenever $\beta<\sigma\leq 1$.

If $\sigma>1$ and $\beta>1$, we estimate, for some $c_\beta>0$,
$$
\cI(x,[\vfi]) \leq    c_\beta (\lg x\rg^{\beta-1}+1)
$$
and  we get
\be\label{longbeta}
\cL^b[\vfi] \geq  \beta \lgx^{\beta-2}  \, b(t,x)\cdot x - C_\beta   (1+ \lg x\rg^{\beta-1})\,.
\ee
Using \rife{bdiss}, the conclusion follows provided $\gamma>1$.  
\qed

\begin{remark}\label{rangegamma}
Let us observe, from the previous Lemma, the coercive case given by $\gamma\geq 2$. In this range, we have that $\liminf_{|x|\to \infty} \frac{\cL^b[\vfi] }\vfi >0$. This holds  for every  $\beta>0$ in the local case $\cL=\cL_0$, otherwise it holds for every $\beta<\sigma$ if $\cI$ satisfies \rife{nu2}. 

By contrast, if $\gamma<2$, then  $\cL^b[\vfi] \to +\infty$ as $|x|\to \infty$ only for $\beta>2-\gamma$. This introduces a natural threshold for nonlocal operators satisfying \rife{nu2}:  since $\beta<\sigma$ is required for the behavior of $\nu$ at infinity, this will lead to require $\gamma>2-\sigma$ for the confinement force to play a significant role. 

Notice that $\gamma+\beta-2>0$ already implies $\gamma>1$ whenever $\beta\leq 1$.  This always occurs if $\sigma\leq 1$, since $\beta<\sigma$.  Otherwise, when $\sigma>1$, it is possible to consider $\beta>1$;  nevertheless the restriction  $\gamma>1$ is still needed to obtain \rife{lyapI},  because one needs the drift term to dominate  the diffusion coming from the jumps involving large values of $x$ (see \rife{longbeta}).  An interpretation is that the jump part contains terms which are comparable to drift components; so they would be controlled by a confining force only if this blows-up at infinity (say, $\gamma>1$).

To sum up this remark, we observe that requiring $\lg x\rg^\beta$ to be a Lyapunov function (i.e. $\cL^b [\lg x\rg^\beta] \to +\infty$ as $|x|\to \infty$) leads to the restriction that $\gamma>1$ and $\gamma>2-\sigma$. 

However, we stress that the restriction $\beta<\sigma$ in the nonlocal case only comes from the behavior at infinity of the Levy measure $\nu(z)$, and the reader may easily think at generalizing this result  to different settings for $\nu$.
\end{remark}

\vskip1em
\subsection{Doubling variables methods}\label{doubcoub}

As we mentioned in the Introduction, two ingredients yield the oscillation estimates of this Section. While Lemma \ref{lyap}  exploits the existence of Lyapunov functions using the confinement property of the drift, the second ingredient comes out by handling  the local oscillation  of solutions through the maximum principle. Here we exploit in full strength the Ishii-Lions doubling variables method, revisited through the ideas of coupling along the lines suggested   in \cite{PP}.  

{   
Since this is the main ingredient in our approach, let us give first a rough idea of how it works, letting aside any technicality. As we explained in the Introduction, our purpose is to show an oscillation estimate of the form \rife{maxpri}, for some quantitative rate $\varpi(t)$. Rephrasing \rife{maxpri}, this amounts to show  that 
\be\label{uz}
{\mathcal U}(t,x,y) \leq {\mathcal \zeta}(t,x,y)
\ee
where 
$$
{\mathcal U}(t,x,y):= u(t,x)-u(t,y) \,,\quad \hbox{and}\quad {\mathcal \zeta}(t,x,y)= K\, \varpi(t)\left[[\vfi(x)+ \vfi(y)]+  \psi(|x-y|)\right],
$$ 
for some constant $K>0$, some Lyapunov function $\vfi(x)$ and some (modulus  of continuity) $\psi(|x-y|)$.

Now, \rife{uz} is nothing but a form of comparison, that can be established via  maximum principle applied in the product space for $(x,y)\in \R^d\times \R^d$. In fact, if $u$ solves \rife{pbgen}, then ${\mathcal U}$ satisfies
\be\label{bas}
\partial_t {\mathcal U}- {\rm tr}(Q(x)D^2u(x)- Q(y)D^2(u(y) ) -\left[  \mathcal{I} ( x, [ u(x) ] )- \mathcal{I} ( y, [ u(y) ] ) \right] +  {\mathcal B}(x,y) \cdot D_{(x,y)}{\mathcal U}=0 
\ee
where ${\mathcal B}(x,y)= (b(x), b(y))$.
Now we observe that if
${\mathcal U}-{\mathcal \zeta}$ has a maximum, then $D^2_{(x,y)}{\mathcal U}\leq D^2_{(x,y)}{\mathcal \zeta}$, which implies that 
\be\label{matt}
{\rm Tr}\left( {\mathcal A}D^2_{(x,y)}{\mathcal U}\right)(x,y)\leq {\rm Tr}\left( {\mathcal A} D^2_{(x,y)}{\mathcal \zeta}\right)(x,y)
\ee
for any matrix $\mathcal A\geq 0$. Such matrix is a degree of freedom to use in order to optimize the estimate on the maximum point.  Since $Q$ is  the  diffusion matrix  of the operator $\cL$,  it is natural to choose $\mathcal A$ of the form
$$
\mathcal A= \begin{pmatrix}   Q(x)   & C \\    C^* &  Q(y)  \end{pmatrix} 
$$
where the crossed terms are given by a generic matrix $C$ with the only constraint  that $\mathcal A\geq 0$. With this choice \rife{matt} becomes
\be\label{argloc}
{\rm tr}(Q(x)D^2u(x)- Q(y)D^2(u(y) ) \leq  {\rm Tr}\left( {\mathcal A} D^2_{(x,y)}{\mathcal \zeta}\right)(x,y)
\ee
and a clever choice of the  matrix $C$ can provide with a useful estimate of the right-hand side. This matrix $C$ plays the role of the coupling in the probabilistic approach. On account of the special form of $\zeta$,  a suitable choice of $C$ will give the estimate \rife{diffelli} of Lemma \ref{locmax} below, where $X_n, Y_n$ should be meant as approximations of $D^2u(x), D^2u(y)$ (in the viscosity solutions' approach).  This is entirely borrowed from \cite{PP}, where the reader may also find (in the Appendix) an extended discussion of the analytic and probabilistic coupling approach.
%and using the special form of $\zeta $ gives
%$$
%{\rm tr}(Q(x)D^2u(x)- Q(y)D^2(u(y) ) \leq K\, \varpi(t) \left( {\rm tr}(Q(x)D^2\vfi(x)) +  {\rm tr}(Q(y)D^2\vfi(y) ) + 
%$$

Somehow we will use a similar idea for nonlocal operators, that we roughly explain here in a model case. Suppose that the Levy measure $\nu$ is symmetric and is normalized with unit mass; then 
%given by the fractional Laplacian. Then 
$$
 \mathcal{I} ( x, [ u(x) ] )- \mathcal{I} ( y, [ u(y) ] ) \simeq    \int_{\R^d}\!\!\int_{\R^d}\![{\mathcal U}(t,x+z,y+z')   -{\mathcal U}(t,x,y)] d\pi(z,z')
 % \int_{\R^d} \{u ( x+z ) -u( x )  - u ( y+z ) +u( y )  \}\nu ( dz ) = \int\int {\mathcal U}(t,(x,y)+ (z,Az))) -{\mathcal U}(t,(x,y))
 $$
 for any coupling $\pi$ of the measures $\nu(dz), \nu(dz')$; on a maximum point of ${\mathcal U}-{\mathcal \zeta}$ this yields
$$
  \mathcal{I} ( x, [ u(x) ] )- \mathcal{I} ( y, [ u(y) ] ) \leq \int_{\R^d}\!\!\int_{\R^d}\! [{\mathcal \zeta}(t,x+z,y+z') -{\mathcal \zeta}(t,x,y)] d\pi(z,z')\,.
$$
 The choice of $\pi$ can now provide with suitable estimates. For instance, if   $\nu$ is the fractional Laplacian,  a simple choice would be  $\pi= (I,A)_{\sharp}(\nu)  $\footnote{$(I,A)_{\sharp}(\nu)$ is the push-forward measure of $\nu$ through the mapping $(I,A)$ from $\R^d$ into $\R^d\times \R^d$, meaning $\int\int \phi(z,z') d(I,A)_{\sharp}(\nu) = \int\int \phi(z,A(z))\, \nu(dz)$}, for any orthogonal matrix $A$. This will  yield the simplified estimate
\be\label{argnonloc}
  \mathcal{I} ( x, [ u(x) ] )- \mathcal{I} ( y, [ u(y) ] ) \leq \int_{\R^d}[{\mathcal \zeta}(t,x+z,y+Az) -{\mathcal \zeta}(t,x,y)] \nu(dz)
\ee
  where $A$ is any orthogonal matrix (which will play the role of the reflection coupling in the probabilistic approach).  
  
Finally, using \rife{argloc} and \rife{argnonloc} in \rife{bas}, where we can replace $\partial_t {\mathcal U}, D_{(x,y)}{\mathcal U}$ with $\partial_t {\mathcal \zeta}, D_{(x,y)}{\mathcal \zeta}$, we will obtain estimates on the maximum point of  ${\mathcal U}-{\mathcal \zeta}$ and, eventually, conclude that  ${\mathcal U}\leq {\mathcal \zeta}$.

To sum up, the strategy towards the oscillation estimate \rife{maxpri} entirely relies on the maximum principle for an (augmented) auxiliary problem in the product space. The justification of this comparison argument finds a natural place in the framework of viscosity solutions, which provides with a rigorous proof  in a context of possibly non smooth solutions. It is worth pointing out that the contributions  of  local and nonlocal diffusions (namely, \rife{argloc} and \rife{argnonloc}) is independently dealt with; this is why we give two corresponding separate statements, the first one being  already developed in our previous work \cite{PP}. 
}
%For the sake of clarity, we state separately the contribution of  local and nonlocal diffusion.

\begin{lemma}\label{locmax}  Let $u$ be a viscosity solution of \rife{pbgen} and suppose that $(x,y)$ is  a local maximum point of the function
$$
u(t,x)-u(t,y) - \zeta(t,x,y)\,,\qquad \zeta(t,\cdot, \cdot) =   [\vfi(x)+ \vfi(y)]+  \psi(|x-y|)
$$
where  $\vfi$ is a  $C^2$ function  and  $\psi:[0,\infty)\to \R_+$ is a positive, increasing   function such that  $\psi(|\xi|) $ is $C^2$  in a neighborhood of $\xi=x-y$. 

For $n>0$, let $X_n,Y_n\in {\mathcal S}^N$ be the matrices given by Theorem \ref{test} where $A= D^2_{(x,y)}\zeta$.  
%which satisfy the inequality
%\be\label{ineqbase}
% \begin{pmatrix}
%X & 0    \\
%\noalign{\medskip} 0 & -Y
% \end{pmatrix}
% \leq    D^2_{(x,y)}\zeta:= \begin{pmatrix}
% D_{xx}\zeta & D_{xy} \zeta    \\
%\noalign{\medskip} D_{xy} \zeta& D_{yy} \zeta
%\end{pmatrix}.
%\ee
Then, setting $Q(x):= \lambda_0 I_d+ \Sigma(x)\Sigma^*(x)$,  we have
\be\label{diffelli}
\begin{split}
\limsup_{n\to \infty} \,\,{\rm tr}(Q(x)X_n- Q(y)Y_n )  & \leq 4\lambda_0 \psi''(|x-y|)  + \frac{\psi'(|x-y|)}{|x-y|} \|\Sigma(x)-\Sigma(y)\|^2
\\
& \qquad \qquad  + {\rm tr}(Q(x)D^2\vfi(x)) +  {\rm tr}(Q(y)D^2\vfi(y) )\,.
\end{split}
\ee
\end{lemma}

\proof This is  a consequence of the method introduced in \cite{IL}. In the above precised form, inequality \rife{diffelli} stems from  \cite[Prop. 2.4]{PP} up to 
a specific computation involving the very definition of $D^2\zeta$ (see  \cite[formula (3.29)]{PP}).
\qed

The next tool is an extension of the previous lemma to nonlocal  operators, as defined in \rife{Levy}. Generalizations of the Ishii-Lions method to nonocal operators have been  developed in \cite{BCCI}, \cite{BCI}, but we need here an improvement which is crucial for the purpose of our results. This improvement relies on  the idea of coupling (by reflection) for jump operators, as we roughly explained  above. However, compared to the probabilistic approach,  here the  coupling   will appear as naturally inserted in the analytical doubling variables argument.  

In the following, we omit the time dependence of the function $u$, since it is irrelevant here. Moreover, we assume that all functions involved are integrable at infinity with respect to the Levy measure $\nu$ (in our context, this is satisfied by all functions which have power growth of order less than $\sigma$).

\begin{lemma} \label{new-cone}
Let the operator $\mathcal I$ be defined by \rife{Levy}, where $\nu$ satisfies \rife{nu2}.
Suppose that $(x,y)$ is  a maximum point of the function
$$
u(x)-u(y) - \zeta(x,y)\,,\qquad \zeta:=   [\vfi(x)+ \vfi(y)]+  \psi(|x-y|)
$$
where  $\vfi$ is a  $C^2$ function  and  $\psi:[0,\infty)\to \R_+$ is a positive, increasing  function such that  $\psi(|\xi|) $ is $C^2$  in a neighborhood of $\xi=x-y$. 

Then, for every $\de\leq \left(\frac{|x-y|}2 \wedge 1\right)$, we have (recall the notations in \rife{notp}):
 \be\label{nonlocal_precise}
 \begin{split}
&  \mathcal{I} ( x,u(x),D_{x} \zeta(x,y) ) -  \mathcal{I} ( y,u(y),-D_{y} \zeta(x,y) ) \leq  \mathcal{I} ( x,\vfi(x) ,D \vfi(x) ) + \mathcal{I} ( y, \vfi(y), D\vfi(y)) 
\\
& \quad \qquad \qquad 
%   \mathcal{I} ( x,\vfi(x) ,D \vfi(x) ) + \mathcal{I} ( y, \vfi(y), D\vfi(y))  
%\\
%& \quad \qquad 
+   4\lambda  \int_0^1 (1-s)\int_{B_\de}   \psi'' ( |x-y|+ 2s(\widehat{x-y}\cdot z))\,  |\widehat{x-y}\cdot z|^2 \frac{dz}{|z|^{d+\sigma}} ds ,  
\end{split}
 \ee
 where $B_\de:= \{ z\in \R^d\,:\, |z|<\de\}$ and $\widehat{x-y}= \frac{x-y}{|x-y|}$.
 \end{lemma}
 
\proof  Let us set
$$
H(x,y):= u(x)-u(y) - \zeta(x,y) = u(x)-u(y) -  [\vfi(x)+ \vfi(y)]-   \psi(|x-y|)\,.
$$
Since $H(x,y)\geq H(x+z,y+z)$ for every $z$,  we have
\be\label{maj10}
\begin{split}
& [u ( x+z ) -u( x )]- [u (y+z ) -u( y )] \\
& \qquad\qquad \leq  \left\{ \vfi( x+z ) -\vfi( x )+  \vfi (y+z ) -\vfi( y )]  \right\}\qquad \forall z\in \R^d\,.
\end{split}
\ee
In addition, we have
$$
D_x\zeta+ D_y\zeta=   D\vfi(x)+ D\vfi(y) \,.
$$
Therefore, for any  $A\subseteq \R^d$ and any Levy measure $\tilde \nu$,  we have
\be\label{locA}
\begin{split}
& \int_{A}   \{u (x+z ) -u( x ) -  \langle D_x\zeta, z  \rangle \mathds{1}_{| z | \leq 1} \} \tilde \nu ( dz ) - 
\int_{A}   \{u (y+z ) -u( y ) +  \langle D_y\zeta, z  \rangle \mathds{1}_{| z | \leq 1} \}\tilde \nu ( dz ) \\
& \quad \leq 
\int_{A}   \{\vfi (x+z ) -\vfi( x ) -  \langle D\vfi(x), z  \rangle \mathds{1}_{| z | \leq 1} \}\tilde \nu ( dz ) +
\int_{A}   \{\vfi (y+z ) -\vfi( y ) -  \langle D\vfi(y), z  \rangle \mathds{1}_{| z | \leq 1} \}\tilde \nu ( dz )\,.
\end{split}
\ee
Notice in particular that \rife{locA} with $A=\R^d$ implies \rife{nonlocal_precise} when $\lambda=0$. But now we wish to exploit the ellipticity, and to this purpose 
we split  the nonlocal terms for $|z|\geq r/2$ and $|z|< r/2$, where $r= |x-y|$. Moreover, in order to handle the {\it small } jumps we decompose the measure $\nu$ as
$$
\nu= \lambda \frac{dz}{|z|^{d+\sigma}}+ \mu\,.
$$
By assumption \rife{nu2}, $\mu$ is nonnegative, and is itself a Levy measure.  We define the ball 
$$
B:= \left\{z\in \R^d \,:\, |z| <  \left(1 \wedge \frac r2\right)\right\}
$$
and we decompose the nonlocal integrals   in three parts, using that 
$$
\nu= \nu \mathds{1}_{B^c} + \mu\mathds{1}_{B} + \lambda \frac{dz}{|z|^{d+\sigma}}\mathds{1}_{B\,.}
$$
For the first two measures, we only use \rife{locA},  applied for both $\tilde \nu= \nu \mathds{1}_{B^c}$ and  $\tilde \nu= \mu\mathds{1}_{B}$. Hence  
\be\label{coerc10}
\begin{split}
& \int_{\R^d} \{ u ( x+z ) -u( x )  -   (D_x\zeta \cdot  z  ) \mathds{1}_{| z | \leq 1} \}   \nu( dz ) 
%\\
%& \qquad \qquad 
- \int_{\R^d} \{u (y+z ) -u( y ) +   (D_y\zeta\cdot  z)  \mathds{1}_{| z | \leq 1} \}  \nu ( dz )
\\ & \leq    \int_{\R^d} \left\{ \vfi( x+z ) -\vfi( x )  - (D\vfi(x)\cdot z ) \mathds{1}_{| z | \leq 1}\right\}  \nu(dz)  
  +     \int_{\R^d} \left\{ \vfi( y+z ) -\vfi( y )  - (D\vfi(y)\cdot z) \mathds{1}_{| z | \leq 1} \right\}  \nu(dz) 
\\
& \quad 
+ \lambda \int_B \{ u ( x+z ) -u( x )  -   (D_x\zeta \cdot  z  ) \}  \frac{dz}{|z|^{d+\sigma}}
-   \lambda \int_B \{ u ( y+z ) -u( y )  +   (D_y\zeta \cdot  z  ) \}  \frac{dz}{|z|^{d+\sigma}}
\\
& \quad - \lambda \int_B \{ \vfi ( x+z ) -\vfi( x )  -   (D\vfi(x)\cdot z ) \}  \frac{dz}{|z|^{d+\sigma}}
- \lambda \int_B \{ \vfi ( y+z ) -\vfi( y )  -   (D\vfi(y)\cdot z ) \}  \frac{dz}{|z|^{d+\sigma}}\,.
\end{split}
\ee
We will use a different strategy for the measure $\lambda \frac{dz}{|z|^{d+\sigma}}\mathds{1}_{B}$, which provides with  the desired ellipticity. Since $H$ achieves  its maximum at $(x,y)$, then  $H(x,y)\geq H(x+z,y+ Az)$ where $A$ is the $d\times d$ matrix   defined as
$$
A:= I_d- 2 (\widehat{x-y}\otimes \widehat{x-y})\,.
$$
{We notice that, for a given $z\in \R^d$ (corresponding to a jump in the $x$-process), $Az$ yields the reflection of $z$ through the hyperplane $\{z\cdot \widehat{x-y}=0\}$. This encodes the idea to use a reflected process in the $y$-variable to optimize the estimate, as in the probabilistic coupling methods.}

Thus we have, by the maximality condition,
\be\label{maj20}
\begin{split}
& [u ( x+z ) -u( x )]- [u (y+Az ) -u( y )] \\
& \qquad\qquad \leq 
 [ \psi(|x-y + z-Az|)-\psi(|x-y|) ]
\\ & \qquad +
  \left\{ \vfi(  x+z ) -\vfi(  x )+  \vfi ( y+Az ) -\vfi(  y )]  \right\} \qquad \forall z\in \R^d\,.
\end{split}
\ee
Integrating \rife{maj20} on $B$ against the measure $\frac{dz}{|z|^{d+\sigma}} $ we get
\begin{align*}
& \int_{B} \left\{ u (  x+z ) -u( x )  -   D_x\zeta \cdot  z\right\}   \frac{dz}{|z|^{d+\sigma}} 
 - \int_{B} \{u ( y+Az ) -u( y ) +   D_y\zeta \cdot A z   \} \frac{dz}{|z|^{d+\sigma}} 
\\ & \leq   \int_{B} \left\{ \vfi( x+z ) -\vfi(x )  - D\vfi(x)\cdot z  \right\} \frac{dz}{|z|^{d+\sigma}} 
 +   \int_{B} \left\{ \vfi( y+Az ) -\vfi( y )  - D\vfi(y)\cdot Az  \right\} \frac{dz}{|z|^{d+\sigma}} 
\\ & \qquad +     \int_{B} \left\{  \psi(|x-y + z-Az|)- \psi(|x-y|) -  \psi'(|x-y|) \widehat{x-y}\cdot(z-Az) \right\} \frac{dz}{|z|^{d+\sigma}} \,.
\end{align*}
Now we use that   $A$ is an orthogonal matrix and the measure  $\frac{dz}{|z|^{d+\sigma}}$ is invariant by rotation. 
Therefore, the previous inequality reads as follows:
\be\label{coerc20}
\begin{split}
& \int_{B} \left\{ u (  x+z ) -u( x )  -   D_x\zeta \cdot  z\right\}   \frac{dz}{|z|^{d+\sigma}} 
 - \int_{B} \{u ( y+ z ) -u( y ) +   D_y\zeta \cdot   z   \} \frac{dz}{|z|^{d+\sigma}} 
\\ & \leq   \int_{B} \left\{ \vfi( x+z ) -\vfi(x )  - D\vfi(x)\cdot z  \right\} \frac{dz}{|z|^{d+\sigma}} 
 +   \int_{B} \left\{ \vfi( y+ z ) -\vfi( y )  - D\vfi(y)\cdot  z  \right\} \frac{dz}{|z|^{d+\sigma}} 
\\ & 
\qquad +   \int_{B} \left\{  \psi(|x-y + z-Az|)- \psi(|x-y|) -  \psi'(|x-y|) \widehat{x-y}\cdot(z-Az) \right\} \frac{dz}{|z|^{d+\sigma}} \,.
\end{split}
\ee
Putting together \rife{coerc10} and \rife{coerc20}, we deduce the following inequality
\be\label{precou}
\begin{split}
&  \mathcal{I} ( x,u(x),D_{x} \zeta(x,y) ) -  \mathcal{I} ( y,u(y),-D_{y} \zeta(x,y) ) \leq  \mathcal{I} ( x,\vfi(x) ,D \vfi(x) ) + \mathcal{I} ( y, \vfi(y), D\vfi(y)) 
\\
& \quad \qquad \qquad 
+  \lambda  \int_{B} \left\{  \psi(|x-y + z-Az|)- \psi(|x-y|) -  \psi'(|x-y|) \widehat{x-y}\cdot(z-Az) \right\} \frac{dz}{|z|^{d+\sigma}}\,.  
\end{split}
\ee
Recalling the definition of $A$, we compute
\begin{align*}
&  \int_{B} \left\{  \psi(|x-y + z-Az|)- \psi(|x-y|) -  \psi'(|x-y|) \widehat{x-y}\cdot(z-Az) \right\} \frac{dz}{|z|^{d+\sigma}}
\\ & \qquad  =   \int_{B} \left\{  \psi( |r  + 2(\widehat{x-y}\cdot z)| ) - \psi(r)  -  2 \psi'(r) (\widehat{x-y}\cdot z) \right\} \frac{dz}{|z|^{d+\sigma}}\,
\end{align*}
where, we recall, $r=|x-y|$. 
Now we estimate the right-hand side.  To this purpose,  notice that for $|z|<\frac r2$ we have $ r  + 2\widehat{x-y}\cdot z\geq 0$. Hence,  Taylor's expansion implies
\begin{align*} & 
\int_{B} \left\{  \psi( |r  + 2(\widehat{x-y}\cdot z)| ) - \psi(r)  -  2 \psi'(r) (\widehat{x-y}\cdot z) \right\} \frac{dz}{|z|^{d+\sigma}}
\\
& \leq 
4 \int_0^1 (1-s) \int_{B} \psi''( r+ 2s(\widehat{x-y}\cdot z)) \,|\widehat{x-y}\cdot z|^2 \frac{dz}{|z|^{d+\sigma}} ds \,.
%\\ & 
%\leq 4\int_0^1 (1-s)\int_{B^-} \psi''( r+ 2s(\widehat{x-y}\cdot z)) |z|^2 \frac{dz}{|z|^{d+\sigma}} ds 
\end{align*}
Finally, coming back to \rife{precou}, we deduce \rife{nonlocal_precise}.
\qed

\subsection{Exponential rate}\label{exprate}

We consider here the case that the drift is strongly confining, in the sense that $\gamma\geq 2$ in \rife{bdiss}. This provides with an exponential decay estimate for a weighted oscillation of the solution of \rife{pbgen}. For the sake of clarity, we distinguish two different statements for the case of local operators and the case of general Levy processes with jumps. 

\begin{theorem}\label{decaysup}
Let $\cL_0$ be defined by \rife{L0},  where $\lambda_0>0$ and $\Sigma$ satisfies  \rife{sigma}.  Assume that $b$ satisfies \rife{bdiss} with $\gamma\geq 2$, and \rife{disso}. Let $u$ be a  (viscosity) solution of 
\be\label{brown}
\begin{cases}
\partial_t u + \cL_0[u]   + b(t,x) \cdot Du     =  0 &  \text{ in } ( 0,T ) \times \R^d, \\ 
u ( 0,x ) = u_0 (x),\qquad  \,&  \text{ in }   \R^d
\end{cases}
\ee
such that $|u(t,x)|\leq C (1+ |x|)^{m}$ for some $m>0$ and $C>0$ (uniformly in time).   

For any $k>0$ such that $u_0 \in L^\infty(\langle x\rangle^{-k})$, there exist positive constants $K, \omega >0$ (depending on $\alpha, \gamma, \lambda_0, \sigma_0, \sigma_1, c_0, d,  k$)  such that $u$ satisfies
\be\label{decayk}
[u(t)]_{\langle x\rangle^k} \leq K \, e^{-\omega t} \, [u_0]_{\langle x\rangle^k}.
\ee
\end{theorem}

In case of nonlocal operators with $\sigma$ fractional growth, the statement is similar except for a restriction in the weight which is used (and the growth of solutions).

\begin{theorem}\label{decaysup2}
Let $\cL$ be defined by \rife{L} where $\Sigma$ satisfies  \rife{sigma},  the Levy measure in $\cI$ satisfies \rife{nu2}, and $\lambda_0+\lambda>0$.  Assume that $b$ satisfies \rife{bdiss} with $\gamma\geq 2$ and   either of the following conditions hold:

(i) $\lambda_0>0$ or $\sigma\in (1,2)$, and \rife{disso} holds true.

(ii) $\lambda_0=0$ and $\sigma\in (0,1]$, and \rife{disso2} holds true.

Let $u$ be a  (viscosity) solution of \rife{pbgen} such that $|u(t,x)|\leq C (1+ |x|)^{m}$ for some $m\in (0,\sigma)$ and $C>0$ (uniformly in time).   

For any $k\in (0,\sigma)$ such that $u_0 \in L^\infty(\langle x\rangle^{-k})$, there exist positive constants $K, \omega >0$ (depending on $\alpha, \gamma, \sigma, (\lambda_0+\lambda)^{-1}, \sigma_0, \sigma_1, \Lambda, c_0, d,  k$)  such that $u$ satisfies
\be\label{decayk2}
[u(t)]_{\langle x\rangle^k} \leq K \, e^{-\omega t} \, [u_0]_{\langle x\rangle^k}.
\ee
%Eventually, if $\nu=0$, the conclusion holds (assuming \rife{disso}) for any $k>0$, and for all solutions $u$ such that $|u(t,x)|\leq C (1+ |x|)^{m}$ for some $m>0$ (without any restriction on $m$).
\end{theorem}

{\bf Proof of Theorem \ref{decaysup} and Theorem \ref{decaysup2}.}  Let us set $\vfi(x)= \lgx^k$, $\Phi(x)= \lgx^p$, where $0<k<\sigma $ and $m<p<\sigma$ (no restriction on $k,p$ is needed where there is no jump part, i.e. for Theorem \ref{decaysup}).  In particular, we have that $u(t,x)= o(\Phi(x))$ as $|x|\to \infty$, uniformly in time.
As in \cite{CLN}, we also use the explicit auxiliary function
\be\label{proppsi}
\psi ( r ) = C_1(1 - e^{- C_{2} r^{\theta}} )
%+ L r^\theta 
%\begin{split} & \hbox{$\psi \in C^2(0,\infty)$, $\psi$ is concave, }  \\
%& 
%\psi ( r ) = 1 - e^{- C_{2} r^{\tau}}  \quad \hbox{for $r \in [ 0, r_{0} ]$}\,,\quad \hbox{$\psi$ is linear on $ [ r_{0} + 1, \infty )$} \\
%&  \min_{(0,\infty)}   \psi' ( r ) = \psi'( r_0+1) := C_{2} \tau ( r_{0} + 1 )^{\tau-1} e^{-C_{2} ( r_{0} + 1 )^{\tau}}
%\end {split}
\ee
where $\theta\in (0,1)$ and $C_1,C_2$ will be chosen later.  Notice that $\psi$ is a bounded, concave function, which is smooth in $(0,\infty)$ and satisfies $\psi \simeq  C_1C_2\, r^\theta $ for $r\to 0$. Moreover, $\psi\leq C_1$ for every $r$.

We claim that, for suitably chosen large constants  $K, C_1$ and $C_2$, and for a convenient choice of $\omega$, we have
\be\label{claim}
u(t,x)- u(t,y) \leq  e^{-\omega t}  \left\{   K\, (\vfi(x)+ \vfi(y)) + \psi(|x-y|) + \vep [\Phi(x)+ \Phi(y)] \right\}+ \frac\vep{T-t}
\ee
for every $\vep>0$ sufficiently small.

As usual, we observe that if  \rife{claim} were not true, then the function 
$$
u(t,x)- u(t,y) -  e^{-\omega t}    \underbrace{\left\{  K  [\vfi(x)+ \vfi(y) ]    +\psi(|x-y|)  + \vep [\Phi(x)+ \Phi(y)] \right\}}_{\zeta(x,y)} -  \frac\vep{T-t} 
$$
has a positive maximum in $[0,T]\times \R^{2d}$ (because $u(t,x)=o(\Phi(x))$ as $|x|\to \infty$), and this positive maximum  is attained at some point $(t, x, y)$ such that, obviously,  $t<T$ and $x\neq y$. In addition, if  $K \geq [u_0]_{\langle x\rangle^k}$, then the positive maximum cannot either be attained at $t=0$.  As a consequence, the point $(t,x,y)$ is a local maximum and we can use  the viscosity inequalities satisfied by $u(t,x)$ and $u(t,y)$, thanks to Theorem \ref{test}. In particular, we deduce the existence of matrices $X_n,Y_n\in {\mathcal S}^N$ such that
$$
 \begin{pmatrix}
X_n & 0    \\
\noalign{\medskip} 0 & -Y_n
 \end{pmatrix}
 \leq    e^{-\omega t} D^2_{(x,y)}\zeta + \frac1n \left( e^{-\omega t} D^2_{(x,y)}\zeta\right)^2
% := \begin{pmatrix}
% D_{xx}\zeta & D_{xy} \zeta    \\
%\noalign{\medskip} D_{xy} \zeta& D_{yy} \zeta
%\end{pmatrix}
$$
and we have 
\begin{align*}
& \frac\vep{(T-t)^2} - \omega e^{-\omega t}   \left(   K  [\vfi(x)+ \vfi(y) ]    +\psi(|x-y|)  + \vep [\Phi(x)+ \Phi( y)] \right) 
\\
& \quad + e^{-\omega t}    \left((b(t,x)-b(t,y)) \cdot  \widehat{x-y}\right)  \psi'(|x-y| ) 
\\
&  \quad + e^{-\omega t} \left\{ b(t,x)\cdot (K D\vfi(x)+ \vep D\Phi(x)) + b(t,y)\cdot (KD\vfi(y) + \vep D\Phi(y))\right\} 
\\
& \qquad \leq   {\rm tr}( Q( x) X_n - Q(y) Y_n ) + {\mathcal I}(  x, u(  t), e^{-\omega t} D_{x} \zeta) 
- {\mathcal I}(  y, u(t) , -e^{-\omega t}D_{y} \zeta)   
\end{align*}
where $\widehat{x-y}= \frac{x-y}{|x-y|}$ and $Q(x)= \lambda_0 I_d+ \Sigma\Sigma^*(x)$.

The right-hand side is estimated through Lemma \ref{locmax} for the local diffusion and Lemma \ref{new-cone} for the nonlocal part. Therefore, letting $n\to \infty$, we obtain
\begin{align*}
& \frac\vep{(T-t)^2} - \omega e^{-\omega t}   \left(   K  [\vfi(x)+ \vfi(y) ]    +\psi(|x-y|)  + \vep [\Phi(x)+ \Phi( y)] \right) 
\\
& \quad + e^{-\omega t}    \left((b(t,x)-b(t,y)) \cdot  \widehat{x-y}\right)  \psi'(|x-y| ) 
\\
&  \quad 
+ e^{-\omega t} \left\{ b(t,x)\cdot (K D\vfi(x)+ \vep D\Phi(x)) + b(t,y)\cdot (KD\vfi(y) + \vep D\Phi(y))\right\} 
\\
& \qquad \leq  e^{-\omega t} \left\{4\lambda_0 \psi''(|x-y|)  + \frac{\psi'(|x-y|)}{|x-y|} \|\Sigma(x)-\Sigma(y)\|^2\right\}\\
& \quad +e^{-\omega t}  \left\{ {\rm tr}(Q(x)[KD^2\vfi+ \vep D^2\Phi](x)) +  {\rm tr}(Q(y)[KD^2\vfi+ \vep D^2\Phi](y) )\right\}
\\ & \quad   +  e^{-\omega t}  K\left\{  \, \mathcal{I} ( x,\vfi(x), D \vfi(x) ) + \mathcal{I} ( y, \vfi(y), D\vfi(y))  \right\}
%\\ & \quad 
+ e^{-\omega t}  \vep \left\{  \, \mathcal{I} ( x,\Phi(x), D \Phi(x) ) + \mathcal{I} ( y, \Phi(y), D\Phi(y))  \right\}
\\ & \quad \qquad \qquad 
%   \mathcal{I} ( x,\vfi(x) ,D \vfi(x) ) + \mathcal{I} ( y, \vfi(y), D\vfi(y))  
%\\
%& \quad \qquad 
+   e^{-\omega t}\, 4\lambda  \int_0^1 (1-s)\int_{B}   \psi'' ( |x-y|+ 2s(\widehat{x-y}\cdot z))\, |\widehat{x-y}\cdot z|^2 \frac{dz}{|z|^{d+\sigma}} ds \,,  
\end{align*}
where 
$$
B:= \left\{z\in \R^d \,:\, |z| <  \left(1 \wedge \frac {|x-y|}4\right)\right\}\,.
$$
%$B:= \{ z\in \R^d\,:\, |z|<\de\}$. 
Henceforth we set  $r=|x-y|$, we use \rife{sigma} and we rearrange terms,  reconstructing the operator $\cL^b$ on the auxiliary functions $\vfi, \Phi$. Then we get
\be\label{step0}
\begin{split}
&
% \frac\vep{(T-t)^2} + e^{-\omega t} \, 
K  \left\{  \mathcal{L}^b[\vfi](x)    -   \omega \vfi(x)  +  \mathcal{L}^b[\vfi](y)   -   \omega\vfi(y)\right\} 
%\\
%& \qquad\quad 
+   \vep   \left( \mathcal{L}^b[\Phi](x)  -\omega \Phi(x) + \mathcal{L}^b[\Phi](y) -\omega \Phi(y) \right)
\\
& \quad +   \left((b(t,x)-b(t,y)) \cdot  \widehat{x-y}\right)  \psi'(r ) 
 <    4\lambda_0 \psi''(r)  +  (2\sigma_0\wedge \sigma_1 r) \sigma_1 \,\psi'(r) + \omega \,\psi(r) 
\\
& \qquad \qquad\qquad +    4\lambda  \int_0^1 (1-s)\int_{B}   \psi'' ( r+ 2s(\widehat{x-y}\cdot z))\, |\widehat{x-y}\cdot z|^2 \frac{dz}{|z|^{d+\sigma}} ds\,.   
\end{split}
\ee
Now we show that the previous inequality leads to a contradiction for suitable choices of $K,C_1, C_2, \omega$.  

We first analyze \rife{step0} for the case that $r=|x-y|\geq 1$. Using  \rife{disso} or \rife{disso2} (they reduce to the same for $|x-y|\geq 1$), and dropping the terms with $\psi''$ (which is negative), we obtain
\be\label{prep}\begin{split} 
& K  \left\{  \mathcal{L}^b[\vfi](x)    -   \omega \vfi(x)  +  \mathcal{L}^b[\vfi](y)   -   \omega\vfi(y)\right\} 
%\\
%& \qquad\quad 
+   \vep   \left( \mathcal{L}^b[\Phi](x)  -\omega \Phi(x) + \mathcal{L}^b[\Phi](y) -\omega \Phi(y) \right)
\\
& \qquad \qquad 
<     (c_0+ 2\sigma_0\sigma_1)  \psi'(r)+ \omega\, \psi(r)   
%\\
%& \leq  -   2e^{-\omega t}  \,   \psi''(|x-y|)   
% \theta\, C_1\,  C_2 e^{- C_{2} r^{\theta}} ( 1-\theta  + C_{2} r^{\theta} )   \, r^{\theta-2}    
 \,.
\end{split}
\ee 
Let us suppose that  $r=|x-y|\geq r_1$, for  a sufficiently large $r_1>1$. This implies that $\max(|x|, |y|)\geq \frac{r_1}2$.  
Due to Lemma \ref{lyap}, choosing $\omega$ sufficiently small (depending only on $\alpha, k$) we have 
\begin{align*}
 & \mathcal{L}^b[\vfi](x)    -   \omega \vfi(x)   \geq \frac{\alpha k}2 \vfi(x) \qquad \hbox{for large $|x|$}
\\
& \mathcal{L}^b[\vfi](x)    -   \omega \vfi(x) \geq -k_0 \qquad \hbox{for every $x$,}
\end{align*}
for some $k_0>0$, and so is for the function $\Phi$.  Notice that $k_0$ also depends on $\alpha, k$ (and possibly on $\sigma,d,   \sigma_0$, but this will not be recalled). Therefore, we can choose $r_1$ sufficiently large (so far,  only depending on $\alpha, k$) such that
$$
\mathcal{L}^b[\vfi](x)    -   \omega \vfi(x) +   \mathcal{L}^b[\vfi](y)   -   \omega\vfi(y) \geq 1 \qquad \forall x,y:|x-y|\geq r_1\,.
$$
If we also estimate $ \psi'(r) \leq  \psi'(r_1) $ for $r\geq r_1$ (by concavity) and $\psi(r)\leq C_1$,   
we deduce from \rife{prep}
\be\label{prep1}
 K- 2\vep k_0- (c_0+ 2\sigma_0\sigma_1) \psi'(r_1)- \omega C_1<0 \qquad \forall r\geq r_1.
\ee
We now fix 
\be\label{choiK}
K:= \omega C_1 + (c_0+ 2\sigma_0\sigma_1)  \psi' (r_1) + 2\vep k_0 + [u_0]_{\langle x\rangle^k}
\ee
and we see that   inequality \rife{prep1} cannot hold.  Therefore, we are left with the possibility that $r<r_1$.  In this case we go back to \rife{step0} and we need to distinguish two cases according to the diffusivity of the nonlocal term $\cI$. 
\vskip1em
{\bf (A) Case $\la_0>0$ or \rife{nu2} with $\sigma\in (1,2)$.}  Here we use condition \rife{disso}, and we deduce from \rife{step0} (using also   $\mathcal{L}^b[\vfi]-\omega\vfi\geq - k_0$, $\mathcal{L}^b[\Phi] -\omega \Phi\geq -k_0$)
\be\label{A0}
\begin{split}
  & -    (2Kk_0+ 2\vep k_0) <     4\lambda_0 \psi''(r) + (c_0+ 2\sigma_0\sigma_1) \psi'(r)  + \omega\, \psi(r)  
 \\ & \qquad +
    4\lambda  \int_0^1 (1-s)\int_{B}   \psi'' ( r+ 2s(\widehat{x-y}\cdot z)) \,| \widehat{x-y}\cdot z|^2 \frac{dz}{|z|^{d+\sigma}} ds\,. 
\end{split}
\ee
We observe that, using once more the concavity of $\psi$, we have $\psi'(r_1)\leq \psi'(r)$ for $r<r_1$, so  the choice of $K$ in \rife{choiK} implies
$$
K\leq  \omega C_1 + 2\vep k_0 + [u_0]_{\langle x\rangle^k}   + (c_0+ 2\sigma_0\sigma_1)  \psi' (r) \qquad \forall r<r_1\,.
$$
Therefore  \rife{A0} yields
\be\label{A1}
\begin{split}
   & -    (2 k_0(\omega C_1 +  [u_0]_{\langle x\rangle^k} ) +  2\vep k_0 (1+ 2k_0) ) <    4\lambda_0 \psi''(r) + (1+2k_0) (c_0+ 2\sigma_0\sigma_1) \psi'(r)  + \omega\, \psi(r)    
 \\ & \qquad +
    4\lambda  \int_0^1 (1-s)\int_{B}   \psi'' ( r+ 2s(\widehat{x-y}\cdot z)) \,|\widehat{x-y}\cdot z|^2 \frac{dz}{|z|^{d+\sigma}} ds\,. 
\end{split}
\ee
We point out that, if $\lambda_0>0$, we could conclude from here   using only the local term $\psi''(r)$ (this would be  similar  as in \cite{PP}). So here the main novelty is in the use of nonlocal diffusion,  neglecting  the term with $\lambda_0$. Hence, using also $\psi\leq C_1$, we  reduce \rife{A1} to the following:
\be\label{A1bis}
\begin{split}
  & -    ( (1+2 k_0)\omega C_1 +  2k_0 [u_0]_{\langle x\rangle^k}  +  2\vep k_0 (1+ 2k_0) ) <    (1+2k_0) (c_0+ 2\sigma_0\sigma_1) \psi'(r)   
 \\ & \qquad +
  4\lambda  \int_0^1 (1-s)\int_{B}   \psi'' ( r+ 2s(\widehat{x-y}\cdot z)) \, |\widehat{x-y}\cdot z|^2 \frac{dz}{|z|^{d+\sigma}} ds\,. 
\end{split}
\ee
Now  we estimate the nonlocal term as follows:
$$
  \int_0^1 (1-s) \int_{B} \psi''( r+ 2s(\widehat{x-y}\cdot z)) \,|\widehat{x-y}\cdot z|^2 \frac{dz}{|z|^{d+\sigma}} ds 
%\\ & 
\leq  \int_0^1 (1-s)\int_{B^-} \psi''( r+ 2s(\widehat{x-y}\cdot z)) \,|\widehat{x-y}\cdot z|^2 \frac{dz}{|z|^{d+\sigma}} ds 
$$
where 
\be\label{bmeno}
B^-:= \{z\in B \,:\, (\widehat{x-y}\cdot z)\leq 0 \}\,.
\ee
Setting $  e:= \widehat{x-y}$, and using polar coordinates, we observe that 
$$
\int_{B^-}  | e\cdot z|^2 \frac{dz}{|z|^{d+\sigma}} \geq  C_{d,\sigma}\, (r \wedge 1)^{2-\sigma} >0\,.
$$
Hence we also estimate, on account of the concavity of $\psi$,
\be\label{psider}\begin{split}
\psi'(r) & \leq \frac2{C_{d,\sigma}} \int_0^1 (1-s)\int_{B^-} \frac{\psi'(r)}{(r \wedge 1)^{2-\sigma}}  |e\cdot z|^2 \frac{dz}{|z|^{d+\sigma}}ds   
\\ & \leq \frac2{C_{d,\sigma}} \int_0^1 (1-s)\int_{B^-} \frac{\psi'(r+ 2s(e\cdot z) )}{(r \wedge 1)^{2-\sigma}}   |e\cdot z|^2 \frac{dz}{|z|^{d+\sigma}}ds
\\
& \leq \frac2{C_{d,\sigma}}  \int_0^1 (1-s)\int_{B^-} \frac{\psi'(r+ 2s(e\cdot z) )}{(r+ 2s(e\cdot z)) \wedge 1)^{2-\sigma}}  |e\cdot z|^2 \frac{dz}{|z|^{d+\sigma}}ds
\end{split}
\ee
where we  also used, in the last step,  that for $z\in B^-$ we have $ 0<r+ 2s(e\cdot z)\leq r$.

Using \rife{psider} in \rife{A1bis}, we obtain
\be\label{A2}
\begin{split}
  & -    ( (1+2 k_0)\omega C_1 +  2k_0 [u_0]_{\langle x\rangle^k}  + 2 \vep k_0 (1+ 2k_0) )     \\ & \qquad \leq 
    \int_0^1 (1-s)\int_{B^-} \left\{ 4\lambda \psi''(\xi_{s,z})+  C \frac{\psi'(\xi_{s,z})}{(\xi_{s,z}\wedge 1)^{2-\sigma}}\right\} |e\cdot z|^2 \frac{dz}{|z|^{d+\sigma}}ds
\end{split}
\ee
where $\xi_{s,z}= r+2 s (e\cdot z)$, and $C=  \frac2{C_{d,\sigma}}(1+2k_0) (c_0+ 2\sigma_0\sigma_1)$.
% is a constant only depending on $\alpha,k,\gamma, \sigma, \sigma_0,\sigma_1$.

Recall that    for $z\in B^-$ we have $r/2\leq \xi_{s,z}\leq   r$. In particular, we have $\xi_{s,z}\leq   r_1$ whenever $r<r_1$.  Using the definition  of $\psi$  we compute
\begin{align*}
& 4\lambda \psi''(\xi)+  C\, \frac{\psi'(\xi)}{(\xi \wedge 1)^{2-\sigma}}   = 
- C_1 \,  \theta \, C_2\, e^{-C_2\xi^\theta}\xi^{\theta-2} \left\{4\lambda(1-\theta)+ 4\lambda \theta C_2 \xi^\theta-  C (\xi^{\sigma-1}\vee \xi) \right\}\,.
\end{align*}
Now we choose $\theta <\sigma-1$, and we fix $C_2$ sufficiently large so that 
$$
  4\lambda \theta C_2 \xi^\theta-  C (\xi^{\sigma-1}\vee \xi) \geq 0 \quad \forall \xi\in (0,r_1)\,.
$$
With this choice of $C_2 $ (which only depends on $\lambda, \theta,   \sigma, C, r_1$), we finally estimate
\begin{align*}
\forall \xi\leq r\leq r_1, \qquad  4\lambda \psi''(\xi)+  C\, \frac{\psi'(\xi)}{(\xi \wedge 1)^{2-\sigma}}& 
\leq - C_1 \, \theta \, C_2\, e^{-C_2\xi^\theta}\xi^{\theta-2} 4\la(1-\theta)\\
&  \leq   - C_1  \, L\, r^{\theta-2}  \qquad \forall r<r_1
\end{align*}
for some constant $L>0$ independent of $C_1$. 
%Then from \rife{nonloc4} we deduce
%\begin{align*}
%& c_0 (1+ 2k_0) \psi'(r)  +
% \lambda    \int_{B} \left\{  \psi( |r  + 2(\widehat{x-y}\cdot z)| ) - \psi(r)  -  2 \psi'(r) (\widehat{x-y}\cdot z) \right\} \frac{dz}{|z|^{d+\sigma}}\,.
%\\
%& \quad \leq  - C_1  \, c_{d,\sigma}\, L\, r^{\theta-\sigma}  \leq  - C_1  \, c_{d,\sigma}\, L\, r_1^{\theta-\sigma} \qquad \forall r<r_1  \,, 
%\end{align*}
Then from  \rife{A2} (and using again \rife{bmeno}) we deduce
$$
 C_1  \left( c_{d,\sigma}\, L\, r_1^{\theta-\sigma}   - \omega  (1+ 2k_0 ) \right)  - 2k_0 [u_0]_{\langle x\rangle^k}-  2\vep k_0(1+ 2k_0)   
 \leq   0\,.
$$
 Choosing $\omega$ sufficiently small and $C_1$ suitably large ($ C_1 \gtrsim n ([u_0]_{\langle x\rangle^k}+\vep k_0) $ for some large $n$) we get a contradiction.
 This shows that \rife{claim} holds true; then letting $\vep \to 0$ and recalling our choices of the parameters $K, C_1$, we proved that 
\begin{align*}
 u(t,x)- u(t,y) &\leq  e^{-\omega t}  \left\{   K\, (\vfi(x)+ \vfi(y)) + \psi(|x-y|) \right\} \\
 & \quad \leq e^{-\omega t} \, (K+ C_1)(\vfi(x)+ \vfi(y)) 
\leq  e^{-\omega t} \,\tilde K \, [u_0]_{\langle x\rangle^k}(\vfi(x)+ \vfi(y)) 
\end{align*}
 for some constant $\tilde K$,  
 %depending on $\alpha, c_0$ (given by \rife{bdiss}--\rife{disso}) and on $d, k, k_0$ (given by the weight $\vfi(x)$).
 which yields  
 \rife{decayk2}.
% The contradiction show that \rife{claim} holds true, and in particular we proved that
%$$
%|u(t,x)- u(t,y)| \leq  e^{-\omega t}  \left\{   K\, (\vfi(x)+ \vfi(y)) + C_1 + \vep [\Phi(t,x)+ \Phi(t,y)] \right\}+ \frac\vep{T-t} \qquad \forall x,y\in \R^d, t<T
%$$
%with the choices of $K, C_1$ given above, and for every $\vep>0$. Letting $\vep \to 0$, and recalling the choices of parameters in \rife{choiK}, \rife{choiC}, we get
%\begin{align*}
%|u(t,x)- u(t,y)| &  \leq  e^{-\omega t}  \tilde K \, [u_0]_{\langle x\rangle^k} \left\{     (\vfi(x)+ \vfi(y)) +1\right\} 
%\\
%& \leq 2 e^{-\omega t}  \tilde K \, [u_0]_{\langle x\rangle^k}  (\vfi(x)+ \vfi(y))  \qquad \forall x,y\in \R^d, t>0\,,
%\end{align*}
%for some constant $\tilde K$ depending on $\alpha, c_0$ (given by \rife{bdiss}--\rife{disso}) and on $d, k, k_0$ (given by the weight $\vfi(x)$).
 \vskip1em
 
 {\bf (B) Case $\la_0=0$ and \rife{nu2} with $\sigma\in (0,1]$. } In this case we use \rife{disso2} in \rife{step0}  and, rather than \rife{A0}, we deduce:
 \be\label{nonloc2}
 \begin{split}
  & -    (2Kk_0+ 2\vep k_0) \leq   ( c_0 (r\wedge 1)^{1-\sigma+\de} + \sigma_1 (\sigma_1 r\wedge2\sigma_0) ) \psi'(r)  + \omega\, \psi(r)   
 \\ & \qquad +
  4\lambda  \int_0^1 (1-s)\int_{B}   \psi'' ( r+ 2s(\widehat{x-y}\cdot z)) \,|\widehat{x-y}\cdot z|^2 \frac{dz}{|z|^{d+\sigma}} ds\,. 
\end{split}
\ee
%which implies, dropping the term with $\lambda_0$ and using $\psi\leq C_1$
% \be\label{nonloc2}
%\begin{split}
% \frac\vep{(T-t)^2}  & -  e^{-\omega t}  (2Kk_0+ 2\vep k_0+ \omega C_1)  \leq   e^{-\omega t} (c_0 (r\wedge 1)^{1-\sigma+\de} + \sigma_1 (\sigma_1 r\wedge2\sigma_0)  ) \psi'(r)
% \\ & \qquad +
%    e^{-\omega t}\, 4\lambda  \int_0^1 (1-s)\int_{B}   \psi'' ( r+ 2s(\widehat{x-y}\cdot z)) |z|^2 \frac{dz}{|z|^{d+\sigma}} ds\,.
%\end{split}
%\ee
Similarly as in \rife{psider},  setting $e:= \widehat{x-y}$ we estimate
\begin{align*}
& (c_0 (r\wedge 1)^{1-\sigma+\de} + \sigma_1 (\sigma_1 r\wedge2\sigma_0)  )  \psi'(r) 
   \leq 
%\\
%& \qquad \leq 
C \int_0^1 (1-s)\int_{B^-} \frac{\psi'(r+ 2s(e\cdot z) )}{(r+ 2s(e\cdot z)) \wedge 1)^{1-\de}}  |e\cdot z|^2 \frac{dz}{|z|^{d+\sigma}}ds
\end{align*}
for some constant $C$ only depending on $d, \sigma, \sigma_0,\sigma_1, c_0$.  From \rife{nonloc2} we get, using also $\psi\leq C_1$,
\be\label{nonloc4bis}
 -   (2Kk_0+ 2\vep k_0+ \omega C_1)  \leq    
      \int_0^1 (1-s)\int_{B^-} \left\{ 4\lambda \psi''(\xi_{s,z})+  C\,   \frac{\psi'(\xi_{s,z})}{(\xi_{s,z}\wedge 1)^{1-\de}}\right\} |e\cdot z|^2 \frac{dz}{|z|^{d+\sigma}}ds
\ee
where $\xi_{s,z}= r+2 s (\widehat{x-y}\cdot z)$. 
%Recall that, for $z\in B^-$, we have $r/2\leq \xi_{s,z}\leq   r<r_1$. 
Now we compute
$$
 4\lambda \psi''(\xi)+  C\,   \frac{\psi'(\xi)}{(\xi \wedge 1)^{1-\de}}  = 
- C_1 \,  \theta \, C_2\, e^{-C_2\xi^\theta}\xi^{\theta-2} \left\{4\lambda(1-\theta)+ 4\lambda \theta C_2 \xi^\theta-  C   (\xi^{\de}\vee \xi) \right\}\,.
$$
and  we fix $C_2$ sufficiently large so that $4\lambda \theta C_2 \xi^\theta-  C    (\xi^{\de}\vee \xi)\geq 0$ for every $\xi\in (0,r_1)$. To this purpose it is enough to choose $\theta<\de$ and $C_2$ sufficiently large depending on $r_1, C, \lambda, \theta$. 

With this choice of  $C_2 $,  we finally estimate
\begin{align*}
4\lambda \psi''(\xi_{s,z})+  C   \frac{\psi'(\xi_{s,z})}{(\xi_{s,z}\wedge 1)^{1-\de}} &\leq   - 4\lambda(1-\theta) C_1 \,  \theta \, C_2\, e^{-C_2\xi_{s,z}^\theta}\xi_{s,z}^{\theta-2}  \\
& \leq   - 4\lambda(1-\theta) C_1 \,  \theta \, C_2\, e^{-C_2 r_1^\theta}r^{\theta-2} \qquad \forall r<r_1\, 
\end{align*}
since we have $r/2\leq \xi_{s,z}\leq   r<r_1$.  Then from \rife{nonloc4bis} we deduce
$$
    -   (2Kk_0+ 2\vep k_0+\omega C_1)   \leq  -  c_{d,\sigma}\, 4\lambda(1-\theta) C_1 \,  \theta \, C_2\, e^{-C_2 r_1^\theta}r_1^{\theta-\sigma}  \,.
 $$
 Recall from  \rife{choiK}, and the definition of $\psi$, that we have
 $$
 K = \omega C_1 + (c_0+ 2\sigma_0\sigma_1) C_1C_2\, e^{-C_2 r_1^\theta}\theta\, r_1^{\theta-1}  + 2\vep k_0 + [u_0]_{\langle x\rangle^k}
 $$
 so that we  can rephrase what we obtained as
\be\label{eta}
\begin{split}
    &   C_1  \left( [4\lambda(1-\theta)c_{d,\sigma}\,   r_1^{1-\sigma}-2k_0 (c_0+ 2\sigma_0\sigma_1)] C_2\, e^{-C_2 r_1^\theta}\theta\, r_1^{\theta-1}   - \omega  (1+ 2k_0 ) \right) 
\\
& \qquad\qquad \qquad  \leq  2k_0 [u_0]_{\langle x\rangle^k}+  2\vep k_0(1+ 2k_0) \,.  
\end{split}
\ee
Suppose now that $\sigma<1$; then we first choose $r_1$ sufficiently large so that $4\lambda(1-\theta)c_{d,\sigma}\,   r_1^{1-\sigma}-2k_0  (c_0+ 2\sigma_0\sigma_1)>0$. Next we choose $C_2$ as we said, depending on $r_1$, and finally, choosing $\omega$ sufficiently small and $C_1$ suitably large ($ C_1 \gtrsim n(1+2k_0) ([u_0]_{\langle x\rangle^k}+2\vep k_0) $) we get a contradiction.
 This shows that \rife{claim} holds true; then letting $\vep \to 0$ and recalling our choices of the parameters $K, C_1$, we proved  \rife{decayk2}. 
 
 Finally, we are left with the case $\sigma=1$. To handle this case  we   observe  that, with the strategy used above, we have $4\lambda \psi''(\xi)+  C   \frac{\psi'(\xi)}{(\xi \wedge 1)^{1-\de}} \leq 0$, hence
 \begin{align*}
 & \int_0^1 (1-s)\int_{B^-} \left\{ 4\lambda \psi''(\xi_{s,z})+  C    \frac{\psi'(\xi_{s,z})}{(\xi_{s,z}\wedge 1)^{1-\de}}\right\} |e\cdot z|^2 \frac{dz}{|z|^{d+\sigma}}ds \\
 & \leq  
 \int_0^1 (1-s)\int_{B^-} \left\{ 4\lambda \psi''(\xi_{s,z})+  C   \frac{\psi'(\xi_{s,z})}{(\xi_{s,z}\wedge 1)^{1-\de}}\right\} |e\cdot z|^2 \frac{dz}{|z|^{d+\sigma-\eta}}ds
% 
% - 4\lambda(1-\theta) C_1 \,  \theta \, C_2\, e^{-C_2 r^\theta}r^{\theta-2} \int_0^1 (1-s)\int_{B^-}|z|^2 \frac{dz}{|z|^{d+1}}
% \\
% & \leq  - 4\lambda(1-\theta) C_1 \,  \theta \, C_2\, e^{-C_2 r^\theta}r^{\theta-2} \int_0^1 (1-s)\int_{B^-}|z|^2 \frac{dz}{|z|^{d+1-\eta}}
% \\
% & \leq  - 4\lambda(1-\theta) C_1 \,  \theta \, C_2\, e^{-C_2 r^\theta}r^{\theta-1}c_{d, \eta} r^{\eta}
 \end{align*}
 for any $\eta>0$. This allows us to follow the same steps as before, replacing $r_1^{1-\sigma}$ with $r_1^\eta$ in \rife{eta}.
 \qed

\begin{remark}\label{gendrift} We point out that the conditions assumed on the drift $b(t,x)$ can be slightly relaxed, at the expense of a mild restriction on the power $\langle x\rangle^k $ which is used in the weighted seminorm. To be precise, we can replace assumptions \rife{disso} with 
\be\label{disso-gen}
\exists \,\, c_0>0\,:\,\quad (b(t,x)-b(t,y))\cdot (x-y) \geq  - c_0 |x-y| (1\vee |x-y|)^{q}   \qquad \forall x,y  \in \R^d\,, \,\forall t>0\,,
\ee
for some $c_0>0, q>0$,  and respectively \rife{disso2} with  
 \be\label{disso2-gen}
 \sigma \in (0,1]\, \hbox{and $\exists \,\de \in (0,1)$: \, $(b(t,x)-b(t,y))\cdot (x-y) \geq  - c_0 |x-y| (|x-y|\wedge 1)^{1-\sigma+\de} (1\vee |x-y|)^{q}$.}
\ee
Then the conclusion of Theorems \ref{decaysup} and Theorem \ref{decaysup2} remain valid for any $k$ such that $\gamma+ k-2\geq q$
(notice that the above Theorems correspond to $q=0$, where any $k>0, \gamma\geq 2$ satisfy this condition). 

To provide with this generalization, the only difference in the proof occurs for large values of $|x-y|$.
This is a straightforward variation, since using
$$
 \mathcal{L}[\vfi](x) + b(t,x)\cdot D\vfi(x)  -   \omega \vfi(x)   \geq \frac{\alpha k}2 |x|^{\ga-2} \vfi(x) \qquad \hbox{for large $|x|$}
$$
one can change  \rife{prep1}   into the following estimate
$$
 K\,\frac{\al\,k}2 \left(\frac{r}2\right)^{\ga-2+k}- 2\vep k_0- (c_0\, r^q+ 2\sigma_0\sigma_1) \psi'(r_1)- \omega C_1<0 \qquad \forall r\geq r_1.
$$
which, whenever $\gamma+ k-2\geq q$,  still implies a choice of $K$ of the same order as \rife{choiK}.

This remark shows that the classical condition which is mostly used in the probabilistic approach, say  
$$
(b(t,x)-b(t,y))\cdot (x-y) \geq  - c_0 |x-y|^2   \qquad \forall x,y  \in \R^d\,, \,\forall t>0
$$
is not only allowed in our setting, but can be  improved e.g. into 
$$
(b(t,x)-b(t,y))\cdot (x-y) \geq  - c_0 ( |x-y| \vee |x-y|^2)   \qquad \forall x,y  \in \R^d\,, \,\forall t>0\,.
$$
\end{remark}

\vskip1em

\begin{remark} 
{   An inspection of the proof of Theorem \ref{decaysup} and  \ref{decaysup2} may give some information on how   $K,\omega$ in  \rife{decayk2} depend on the constants appearing in assumptions \rife{sigma}, \rife{nu2},  \rife{bdiss}-\rife{disso2}, as well as on other relevant parameters such as the dimension $d$. In fact, an explicit choice of $K$ is provided by  \rife{choiK} and a quantitative estimate of $K$, as well as of $\omega$, will mostly depend on the parameter $r_1$ and on the choice of the function $\psi(\cdot)$. 

Notice that the parameter $r_1$ only relies on the Lyapunov function $\vfi$, and possibly depends on the dimension $d$ (in the model case with Ornstein-Ulhenbeck  drift $b(x)=x$, we have $r_1= O(\sqrt d)$). However,  the explicit choice of $\psi$ made in \rife{proppsi} may not lead to optimal estimates  of $\omega$ (neither of $K$, in turn). For  example, in case that $\cL$ is the Laplace operator and $b=x$ is the Ornstein-Ulhenbeck drift, a better choice for $\psi$ would be to take 
the solution of the ODE $4\psi''(r)+ \omega \psi(r)+ \omega L=0$ for $r\in (0,r_1)$, where $L$ is some suitable constant and $\psi(0)=0, \psi'(r_1) \simeq 0$. This would lead to  %$\omega$ smaller than the first eigenvalue in $(0,r_1)$, meaning e.g. 
$\omega= O( r_1^{-2})= O( d^{-1})$, as it is expected for this model case, but the same estimate would not be achieved  if $\psi$ is given by \rife{proppsi}. 

%This is to mention that better bounds can be obtained if the auxiliary function $\psi$ is built through more sophisticated  constructions (which involve possibly auxiliary eigenvalue problems).

This is to mention that, while  we opted for the  explicit  choice \rife{proppsi}  in order to give a unifying proof, sufficiently simple, which works in  all different settings (local and nonlocal operators, diffusive or sub-diffusive, slowly and strongly confining drifts, etc...), the same approach could yield better  bounds if the auxiliary function $\psi$ is built through more sophisticated  constructions (which involve possibly auxiliary eigenvalue problems).
%allow for refined quantitative versions,  up to modification of the function $\psi$.  
 %We are aware that, in order to optimize the constants $\omega, K$, one should build more refined versions of $\psi$.  
 A similar discussion is  done in \cite{Eberle+al}, where quantitative explicit bounds are obtained using  the  probabilistic  coupling method.  This could also be possible here, but at the expense of  either restricting to more specific examples or adding an extra source of technical details. 
%In other words, the method that we use can certainly contain refined quantitative versions, but this is outside our present scope here, and should be left to   further enquiry devoted to some   model case. }
}
\end{remark}

\subsection{Sub-exponential rate}

We now come at the slowly confining case, where, for  Lyapunov functions $\vfi$ of power type, we have $\mathcal L^b (\vfi) \to \infty$ but $\frac{\mathcal L^b (\vfi)}\vfi \to 0$. This occurs when $\gamma<2$ in \rife{bdiss}. In this case one can prove a sub-exponential decay rate (so-called sub-geometrical decay, from discrete models). We start with  a general statement of  polynomial decay of the weighted oscillation. A similar strategy as in Theorem \ref{decaysup2} will work, but provided the weight is sufficiently strong. 
As in Section \ref{exprate}, we give first a separate statement  for the purely local case: indeed, if the diffusion is local the full range of $\gamma\in (0,2)$ can be exploited and we have estimates even in  the regime of {\it degenerate confining}, say when $\gamma\in (0,1)$ in \rife{bdiss}.

\begin{theorem}\label{degconf}
Let $\cL_0$ be defined by \rife{L0} where $\Sigma$ satisfies  \rife{sigma} and  $\lambda_0>0$.  Assume that $b$ satisfies \rife{bdiss}-\rife{disso} with $\gamma \in (0, 2)$.
Let $u$ be a  (viscosity) solution of 
\be\label{pbloc}
\begin{cases}
\partial_t u + \cL_0[u]   + b(t,x) \cdot Du     =  0 &  \text{ in } ( 0,T ) \times \R^d, \\ 
u ( 0,x ) = u_0 (x),\qquad  \,&  \text{ in }   \R^d
\end{cases}
\ee
such that  $|u(t,x)|\leq C (1+ |x|)^{m}$ for some $m>0$ and $C>0$ (uniformly in time).   

 Then, for any  $k>2-\gamma$ and $u_0 \in L^\infty(\langle x\rangle^{-k})$, and for any   $\bar k >k$,  there exists a positive constant  $M>0$ (depending on $\alpha, \gamma, c_0, d,  \lambda_0, \sigma_0, \sigma_1, k, \bar k$)  such that $u$ satisfies
\be\label{decayk-ter}
[u(t)]_{\langle x\rangle^{\bar k}} \leq M \, (1+t)^{-q} \, [u_0]_{\langle x\rangle^k}\qquad \hbox{where $q=  \frac{\bar k-k}{2-\gamma}$.}
\ee
\end{theorem}
\qed
 
 \vskip1em
 
For nonlocal diffusions we have a similar result, although the presence of the   jump part (in particular,  the behavior of the kernel at infinity) restricts the result to the case that $1<\gamma <2$, as observed in Remark \ref{rangegamma}.

\begin{theorem}\label{decayslow_frac}   
Let $\cL$ be defined by \rife{L} where $\Sigma$ satisfies  \rife{sigma},  the Levy measure in $\cI$ satisfies \rife{nu2}, and $\lambda_0+\lambda>0$.  Assume that $b$ satisfies \rife{bdiss} with $\gamma \in (2-\sigma, 2)$ such that $\gamma>1$,   and that  either of the following conditions hold: 

(i) $\lambda_0>0$ or $\sigma\in (1,2)$, and \rife{disso} holds true.

(ii) $\lambda_0=0$ and  $\sigma\in (0,1]$, and \rife{disso2} holds true.

Let $u$ be a  (viscosity) solution of \rife{pbgen} such that  $|u(t,x)|\leq C (1+ |x|)^{m}$ for some $m\in (0,\sigma)$ and $C>0$ (uniformly in time).   

 Then, for any $k\in (0,\sigma)$ with $k>2-\gamma$ and $u_0 \in L^\infty(\langle x\rangle^{-k})$, and for any   $\bar k >k$ with $k<\bar k<\sigma$,  there exists a positive constant  $M>0$ (depending on $\alpha, \gamma, \sigma, (\lambda+ \lambda_0)^{-1}, \sigma_0, \sigma_1, \Lambda, c_0, d,  k, \bar k$)  such that $u$ satisfies
\be\label{decayk}
[u(t)]_{\langle x\rangle^{\bar k}} \leq M \, (1+t)^{-q} \, [u_0]_{\langle x\rangle^k}\qquad \hbox{where $q=  \frac{\bar k-k}{2-\gamma}$.}
\ee
\end{theorem}

{\bf Proof of Theorem \ref{degconf} and Theorem \ref{decayslow_frac}.} \quad We give the proof for the nonlocal case, the other being similar, except for the range of $\gamma$ which is allowed in the Lyapunov function (from \rife{lyapL0}).
 
Let us set $\vfi(x)= \lgx^k$ ($k<\sigma$), $\Phi(x)=\lgx^p$ ($m<p<\sigma$). With the same notations of Theorem \ref{decaysup2}, and $\psi$ given by \rife{proppsi}, we aim at proving that  
\be\label{claim2}
u(t,x)- u(t,y) \leq  (1+t)^{-q} \left\{   K\, (\vfi(x)+ \vfi(y)) + \psi(|x-y|) + \vep [\Phi(t,x)+ \Phi(t,y)] \right\}+ \frac\vep{T-t}
\ee
for every $x,y\in \R^d, t>0$, and for every $\vep>0$ sufficiently small.

As usual, we observe that if  \rife{claim2} were not true, then the function 
$$
u(t,x)- u(t,y) -  (1+t)^{-q}    \underbrace{\left\{  K  [\vfi(x)+ \vfi(y) ]    +\psi(|x-y|)  + \vep [\Phi(t,x)+ \Phi(t,y)] \right\}}_{\zeta(x,y)} -  \frac\vep{T-t} 
$$
admits a point $(t,x,y)$ of global maximum with $x\neq y, t<T$ and $t>0$,   provided $K\geq  [u_0]_{\langle x\rangle^k}$. 
 
Then we use Theorem \ref{test}, writing the viscosity inequalities for $u$, and we proceed as in Theorem \ref{decaysup2} using Lemma \ref{locmax} and Lemma \ref{new-cone}. We end up with a similar kind of inequality as \rife{step0}:
\be\label{newprep-frac}
\begin{split} 
& \frac\vep{(T-t)^2} +  (1+t)^{-q}  \, K  \left\{  \mathcal{L}^b[\vfi](x)  -   \frac q{(1+t)} \vfi(x)  +  \mathcal{L}^b[\vfi](y)  -    \frac q{(1+t)} \vfi(y)\right\}  \\
&   +  (1+t)^{-q}  \, \vep   \left( \mathcal{L}^b[\Phi](x)- \frac q{(1+t)} \Phi(x) + \mathcal{L}^b[\Phi](y)- \frac q{(1+t)} \Phi(y) \right)
\\
& + (1+t)^{-q}    \left((b(t,x)-b(t,y)) \cdot  \widehat{x-y}\right)  \psi'(r )  -  q(1+t)^{-q-1} \psi(r) \\ & \qquad \leq (1+t)^{-q}  \left\{4\lambda_0 \psi''(r)  + 
 (2\sigma_0\wedge \sigma_1 r) \sigma_1 \,\psi'(r)  \right\} 
\\
& \qquad + (1+t)^{-q} 4\lambda  \int_0^1 (1-s)\int_{B}   \psi'' ( r+ 2s(\widehat{x-y}\cdot z)) |\widehat{x-y}\cdot  z|^2 \frac{dz}{|z|^{d+\sigma}} ds \,.
%\\
%& \leq  -   2e^{-\omega t}    \psi''(|x-y|)  
%+ 4e^{-\omega t} \, \psi'(|x-y|) (D\vfi(x)-D\vfi(y))\cdot \widehat{x-y}
% \theta\, C_1\,  C_2 e^{- C_{2} r^{\theta}} ( 1-\theta  + C_{2} r^{\theta} )   \, r^{\theta-2}    
\end{split}
\ee
where $r=|x-y|$. 
We start by considering the case $q=0$. As we did in Theorem \ref{decaysup2}, we first suppose that $r=|x-y|\geq r_1$ for a sufficiently large $r_1>1$; in this range assumptions \rife{disso} and \rife{disso2} reduce to the same.  Moreover, since one between $|x|, |y|$ is larger than $r_1/2$, and since $k>2-\gamma$ implies $\cL^b[\vfi]\to +\infty$ as $|x|\to \infty$ (from Lemma \ref{lyap}),  we can choose $r_1>1$ such 
that
$$
 \mathcal{L}^b[\vfi](x)  +  \mathcal{L}^b[\vfi](y)\geq \frac{\alpha k}2 \left(\langle x\rangle \vee \langle y\rangle\right)^{k+\gamma-2}\geq 1 \,.
 %\qquad \forall r\geq r_1\,.
 $$
Moreover, there exists $k_0$ such that   $\mathcal{L}^b[\vfi]\geq - k_0$, and so is for $\Phi$.   Therefore,   from \rife{newprep-frac} written with $q=0$ we deduce (using \rife{disso} or  \rife{disso2} and dropping the terms with $\psi''$)
$$
  K -2\vep k_0 \leq [c_0  +  2\sigma_0 \sigma_1 ] \,\psi'(r) \qquad \hbox{if $r\geq r_1>1$}\,.
$$
Since, by concavity, $ \psi'(r)\leq \psi'(r_1)$ for $r\geq r_1$, the above inequality cannot hold if we choose 
$$
K=  2\vep k_0    + (c_0  + 2\sigma_0  \sigma_1)  \psi'(r_1) + [u_0]_{\langle x\rangle^k}\,.
$$ 
On the other hand, for $r<r_1$,   just using
$ \mathcal{L}^b[\vfi]\geq -k_0,  \mathcal{L}^b[\phi]\geq -k_0$,  we get
\begin{align*} 
&    \left((b(t,x)-b(t,y)) \cdot  \widehat{x-y}\right)  \psi'(r ) <  4\lambda_0 \psi''(r) + (2\sigma_0\wedge \sigma_1 r) \sigma_1 \,\psi'(r) \\
& \qquad + 4\lambda  \int_0^1 (1-s)\int_{B}   \psi'' ( r+ 2s(\widehat{x-y}\cdot z)) |\widehat{x-y}\cdot  z|^2 \frac{dz}{|z|^{d+\sigma}} ds + 2 K k_0+ 2\vep k_0\,.
%\\
%& \leq  -   2e^{-\omega t}    \psi''(|x-y|)  
%+ 4e^{-\omega t} \, \psi'(|x-y|) (D\vfi(x)-D\vfi(y))\cdot \widehat{x-y}
% \theta\, C_1\,  C_2 e^{- C_{2} r^{\theta}} ( 1-\theta  + C_{2} r^{\theta} )   \, r^{\theta-2}    
\end{align*}
Henceforth, using $\lambda_0+ \lambda>0$,   we argue as in Theorem \ref{decaysup2} in order to get at a contradiction with suitable choices of $C_2$ (depending on $r_1$) and $C_1\gtrsim [u_0]_{\langle x\rangle^k}$.  This proves  \rife{claim2}   for the case $q=0$. 
%\begin{theorem}\label{decayslow}
%Assume that $b$ satisfies \rife{bdiss} with $\gamma \in (0, 2)$, and \rife{disso}. Let $u$ be a  (viscosity) solution of \rife{brown} such that $|u(t,x)|\leq C (1+ |x|)^{m}$ for some $m>0$ and $C>0$ (uniformly in time).   
%
%For any $k>2-\gamma$  such that $u_0 \in L^\infty(\langle x\rangle^{-k}dx)$, and for any   $\bar k >k$,  there exists a positive constant  $M>0$ (depending on $\alpha, \gamma, c_0, d,  k, \bar k$)  such that $u$ satisfies
%\be\label{decayk}
%[u(t)]_{\langle x\rangle^{\bar k}} \leq M \, (1+t)^{-q} \, [u_0]_{\langle x\rangle^k}\qquad \hbox{where $q=  \frac{\bar k-k}{2-\gamma}$.}
%\ee
%\end{theorem}
%As in Theorem \ref{decaysup}, a suitable choice of $C_1,C_2$ in the function $\psi$ yields a contradiction, so \rife{claim2} with $q=0$ is proved. 
After letting $\vep\to 0$, we deduce that, for any $k>2-\gamma$, we have
$$
u(t,x)- u(t,y) \leq        K\, (\lgx^k+ \lg y \rg^k) + \psi(|x-y|)  \leq K\, (\lgx^k+ \lg y \rg^k) + C_1
$$
which implies (after the choices of $K, C_1$)
\be\label{osci}
[u(t)]_{\langle x\rangle^k}\leq C\, [u_0]_{\langle x\rangle^k} \qquad \forall t>0\,,
\ee
where $C$ is a constant only depending on $\alpha, \gamma, (\la +\la_0)^{-1}, c_0, d, \sigma_0, \sigma_1, \sigma, \de$ and $k$. In particular, we have proved that, for all weights $\lgx^k$ with $ 2-\gamma<k<\sigma$, the weighted oscillation of $u$ at time $t$ is controlled by the initial one.

Now we improve the estimate, by proving \rife{claim2} with $q\neq 0$. To this purpose, we observe that, on the maximum point $(t,x,y)$, we have
$$
(1+t)^{-q} [\vfi(x)+ \vfi(y)] K \leq u(t,x)-u(t,y) \leq [\lgx^{k}+ \lg y\rg^{k}] [u(t)]_{\langle x\rangle^k}\,.
$$ 
Choosing  $\vfi(x)= \lgx^{\bar k}$, with $k< \bar k<\sigma$, due to \rife{osci} we estimate
$$
(1+t)^{-q} K \leq  C\, \frac{[\lgx^{k}+ \lg y\rg^{k}]}{[\lgx^{\bar k}+ \lg y\rg^{\bar k}]} [u_0]_{\langle x\rangle^k}
$$
and choosing $q= \frac{\bar k-k}{2-\gamma}$ we get
\be\label{txy}
\frac1{1+t} \leq C\, \left(\frac{[u_0]_{\langle x\rangle^k}}K\right)^{\frac1q} \frac1 {   [\lgx^{2-\gamma}+ \lg y\rg^{2-\gamma}]} 
\ee
for a possibly different $C$ depending on the same constants as before, as well as on $q$.  Thanks to \rife{txy}, there exists a sufficiently large constant $L$  such that if we choose $K> L\, [u_0]_{\langle x\rangle^k}$ then we have
$$
\mathcal{L}^b[\vfi](x)  -   \frac q{(1+t)} \vfi(x) \to +\infty \qquad \hbox{as $|x|\to \infty$,}
$$
where we used  Lemma \ref{lyap} (and the fact that $\bar k>  2-\gamma$).  In particular we have $\mathcal{L}^b[\vfi](x)  -   \frac q{(1+t)} \vfi(x)\geq -k_0$ for some $k_0>0$. Similarly we may reason for $\Phi(x)= \lgx^p$, where we can assume, without loss of generality, that $p>2-\gamma$. Therefore, we can choose $r_1$ sufficiently large so that, whenever $r=|x-y|\geq r_1$ we have
$$
\mathcal{L}^b[\vfi](x)  -   \frac q{(1+t)} \vfi(x)  +  \mathcal{L}^b[\vfi](y)  -    \frac q{(1+t)} \vfi(y)\geq 1\, \quad \hbox{and}\, \quad \mathcal{L}^b[\Phi](x)- \frac q{(1+t)} \Phi(x)\geq -k_0\,.
$$
Hence we deduce from \rife{newprep-frac}
% (using as before the concavity of $\psi$)
\be\label{newprep2}
\begin{split}  
  K & <     2\vep k_0    + [c_0  + 2\sigma_0 \sigma_1 ] \psi'(r)+ \psi (r) \, \frac q{1+t}  
\\
& <  2\vep k_0    +  [c_0  + 2\sigma_0  \sigma_1 ] \psi'(r_1)+ C_1 \, C \left(\frac{[u_0]_{\langle x\rangle^k}}K\right)^{\frac1q}   \qquad \forall r\geq r_1 \,,
\end{split}
\ee
where we used estimate \rife{txy}, and that $\psi(r)\leq C_1, \psi'(r)\leq \psi'(r_1)$ for $r>r_1$.  Here $C$ is a  generic constant independent of $r_1, K,C_1, C_2, u_0$. 
Henceforth, we fix
\be\label{choiC_1}
C_1:= n  [u_0]_{\langle x\rangle^k}+  n^{\frac1{q+1}}\, \vep   k_0  (1+4k_0) 
\ee
where $n$ is a sufficiently large number to be chosen later. Moreover, we also fix
\be\label{choiK2}
K:= L \, (1+ n^{\frac q{q+1}})  [u_0]_{\langle x\rangle^k}  + 4\vep k_0    +   [c_0  + 2\sigma_0  \sigma_1 ] \psi'(r_1) 
\ee
where $L$ is a  sufficiently large constant. With this choice of $C_1, K$, we have
\be\label{psiq}
C_1 \, C \left(\frac{[u_0]_{\langle x\rangle^k}}K\right)^{\frac1q} \leq   C_1 \, C \frac 1{\left(L \, (1+ n^{\frac q{q+1}})\right)^{\frac1q} }\leq \frac C{L^\frac1q} \left\{ n^{\frac q{q+1}} [u_0]_{\langle x\rangle^k} +     \vep   k_0  (1+4k_0) \right\} 
\ee
and \rife{newprep2} implies
$$
K \leq 2\vep k_0    +  [c_0  + 2\sigma_0  \sigma_1 ] \psi'(r_1) + \frac C{L^\frac1q} \left\{ n^{\frac q{q+1}} [u_0]_{\langle x\rangle^k} +     \vep   k_0  (1+4k_0) \right\} \,.
$$
On account of \rife{choiK2}, this inequality cannot hold if $L$ is chosen suitably large (only depending on the constant $C$ and $k_0$). 
%
%Now we fix
%\be\label{choiK2}
%K:= L \, \left( [u_0]_{\langle x\rangle^k} + \left(C_1\, [u_0]_{\langle x\rangle^k}^{\frac1q}\right)^{\frac q{q+1}}\right)+ 4\vep k_0    + 2[c_0  + 2\sigma_0  \sigma_1 ] \psi'(r_1)   
%\ee
%where $L$ is a  sufficiently large constant; in particular, we can choose $L$ suitably large so that
%$$
%K\geq L\,  \left(C_1\, [u_0]_{\langle x\rangle^k}^{\frac1q}\right)^{\frac q{q+1}} \quad \hbox{implies} \quad C_1 \, C \left(\frac{[u_0]_{\langle x\rangle^k}}K\right)^{\frac1q}\leq \frac K2\,.
%$$
%With this latter property, and the choice of $K$, we see that \rife{newprep2} cannot hold. Observe that $L$ only depends on $C$. 

We are left to analyse the case that $r<r_1$, which we  split   in the two ranges of $\sigma$.
 \vskip0.4em
 {\bf (A) Case $\la_0>0$ or \rife{nu2} with $\sigma\in (1,2)$. }
 
In this case  we use \rife{disso} and we deduce from \rife{newprep-frac}
\be\label{psiq-1}
\begin{split} 
&   -    (2Kk_0+ 2\vep k_0) \leq   4\lambda_0 \psi''(r) + [c_0+ (2\sigma_0\wedge \sigma_1 r) \sigma_1]  \psi'(r )  +  \frac q{(1+t)}  \psi(r)  \\
& \qquad  +     4\lambda  \int_0^1 (1-s)\int_{B}   \psi'' ( r+ 2s(\widehat{x-y}\cdot z)) |\widehat{x-y}\cdot  z|^2 \frac{dz}{|z|^{d+\sigma}} ds \,.
%\\
%& \leq  -   2e^{-\omega t}    \psi''(|x-y|)  
%+ 4e^{-\omega t} \, \psi'(|x-y|) (D\vfi(x)-D\vfi(y))\cdot \widehat{x-y}
% \theta\, C_1\,  C_2 e^{- C_{2} r^{\theta}} ( 1-\theta  + C_{2} r^{\theta} )   \, r^{\theta-2}    
\end{split}
\ee
We consider the case that $\lambda_0=0$ (otherwise, the proof would follow easily by dropping the term with $\la$). Then, using $\psi'(r_1) \leq \psi'(r) $ (because $r<r_1$) we estimate $K$ from \rife{choiK2} as
$$
K\leq  L \, (1+ n^{\frac q{q+1}})  [u_0]_{\langle x\rangle^k}  + 4\vep k_0    +  [c_0  + 2\sigma_0  \sigma_1 ] \psi'(r)  
$$
%L \, \left( [u_0]_{\langle x\rangle^k} + \left(C_1\, [u_0]_{\langle x\rangle^k}^{\frac1q}\right)^{\frac q{q+1}}\right) + 4\vep k_0    + 2[c_0  + 2\sigma_0  \sigma_1 ] \psi'(r)   
and we insert this estimate in \rife{psiq-1} obtaining
\be\label{hope0}
\begin{split} 
&    -    \left\{ 2\vep k_0  (1+4k_0) + 2k_0\, L \, (1+ n^{\frac q{q+1}})  [u_0]_{\langle x\rangle^k} \right\}
\leq    (1+ 2 k_0) [c_0+ 2\sigma_0  \sigma_1]  \psi'(r )  +  \frac q{(1+t)}  \psi(r)  \\
& \qquad  +     4\lambda  \int_0^1 (1-s)\int_{B}   \psi'' ( r+ 2s(\widehat{x-y}\cdot z)) |\widehat{x-y}\cdot  z|^2 \frac{dz}{|z|^{d+\sigma}} ds \,,\qquad \forall r<r_1\,.
\end{split}
\ee
Now we observe that  \rife{txy} and \rife{psiq},   together with   $\psi(r)\leq C_1$, imply
\begin{align*}
\frac q{1+t} \psi \leq q\,C\, C_1 \left(\frac{[u_0]_{\langle x\rangle^k}}K\right)^{\frac1q}  
& \leq q\,\,   \frac C{L^\frac1q} \left\{ n^{\frac q{q+1}} [u_0]_{\langle x\rangle^k} +     \vep   k_0  (1+4k_0) \right\}   
\\
& \leq  2k_0\, L \, (1+ n^{\frac q{q+1}})  [u_0]_{\langle x\rangle^k} + 2\vep k_0  (1+4k_0) 
\end{align*}
where last inequality follows by  choosing  $L$ large enough.
%, 
%$$
%\frac q{1+t} \psi \leq  2k_0\, L \, (1+ n^{\frac q{q+1}})  [u_0]_{\langle x\rangle^k} + 2\vep k_0  (1+4k_0) \,.
%$$
Using this estimate,  we deduce from \rife{hope0}
%\begin{align*}
%& -  \left\{2\vep   k_0  (1+4k_0) + 2k_0\, L \, (1+ n^{\frac q{q+1}})  [u_0]_{\langle x\rangle^k}\right\} -  q\,\,   \frac C{L^\frac1q} \left\{ n^{\frac q{q+1}} [u_0]_{\langle x\rangle^k} +     \vep   k_0  (1+4k_0) \right\}  
%\\
%& \qquad 
%\leq   (1+ 2k_0) [c_0+ 2\sigma_0  \sigma_1]  \psi'(r )  +   4\lambda  \int_0^1 (1-s)\int_{B}   \psi'' ( r+ 2s(\widehat{x-y}\cdot z)) |z|^2 \frac{dz}{|z|^{d+\sigma}} ds  
%\end{align*}
%which implies, up to choosing a large $L$,  
\be\label{hopefully}\begin{split}
& -  4\vep   k_0  (1+4k_0) - 4k_0\, L \, (1+ n^{\frac q{q+1}})  [u_0]_{\langle x\rangle^k}  \leq   (1+ 2k_0) [c_0+ 2\sigma_0  \sigma_1]  \psi'(r )\\
& \qquad \qquad	\quad 
   +   4\lambda  \int_0^1 (1-s)\int_{B}   \psi'' ( r+ 2s(\widehat{x-y}\cdot z)) |\widehat{x-y}\cdot z|^2 \frac{dz}{|z|^{d+\sigma}} ds \qquad \forall r<r_1\,.
%\\
%& \leq 2k_0\, c_0 \psi'(r) + 2\vep   K_0  (1+2k_0) + 2k_0  [u_0]_{\langle x\rangle^k}
%\qquad \forall r<r_1\, 
\end{split}
\ee
We handle the right-hand side exactly as in Theorem \ref{decaysup2}. In particular, with a convenient choice of $C_2$ we estimate
\begin{align*}
(1+ 2 k_0) [c_0+ 2\sigma_0  \sigma_1]  \psi'(r )&   +   4\lambda  \int_0^1 (1-s)\int_{B}   \psi'' ( r+ 2s(\widehat{x-y}\cdot z)) |\widehat{x-y}\cdot z|^2 \frac{dz}{|z|^{d+\sigma}} ds
\\
& \leq - \tilde C\, C_1 \, r_1^{\theta-\sigma} \qquad \forall r<r_1
\end{align*}
for some constant $\tilde C$ independent of $C_1, K$. Inserting this inequality in \rife{hopefully}, we obtain
\be\label{newprep3} 
  -  4\vep   k_0  (1+4k_0) - 4k_0\, L \, (1+ n^{\frac q{q+1}})  [u_0]_{\langle x\rangle^k}
\leq -   C_1\, \tilde C\, r_1^{\theta-\sigma}
 \,,\qquad \forall r<r_1\,.
%\\
%& \leq 2k_0\, c_0 \psi'(r) + 2\vep   K_0  (1+2k_0) + 2k_0  [u_0]_{\langle x\rangle^k}
%\qquad \forall r<r_1\, 
\ee
Recalling the value of $C_1$ from \rife{choiC_1}, we see that the above inequality cannot hold, if we finally choose $n$ sufficiently large (depending on all previous parameters).

%Notice that last two terms in the left are of the same order as far as $C_1$ is concerned. 
%%We compute $\psi',\psi''$ and we rearrange terms; we get
%%\be\label{newprep3}\begin{split} 
%%\frac\vep{(T-t)^2} & +   (1+t)^{-q}  C_1 \left\{ \theta \, C_2 e^{- C_{2} r^{\theta}} \, r^{\theta-2} ( 2(1-\theta)  +2 C_{2} r^{\theta} - c_0(1+ 4k_0) r)-  \tilde C\left(\frac{[u_0]_{\langle x\rangle^k}}{ C_1}\right)^{\frac 1{q+1}} \right\} 
%%\\ & \qquad \leq (1+t)^{-q}\left\{ 2\vep   K_0  (1+4k_0) + 2k_0\, L  [u_0]_{\langle x\rangle^k}\right\}\,.
%%\end{split}
%%\ee
%%Here we first fix $C_2$ so that
%%$$
%%m:= \min_{r\in (0, r_1]}\left\{ \theta \, C_2 e^{- C_{2} r^{\theta}} \, r^{\theta-2} ( 2(1-\theta)  + 2C_{2} r^{\theta} - c_0(1+ 4k_0)  r)  \right\}  >0\,.
%%$$
%Then we choose $C_1$ so that
%$$
%C_1:= n ( 2\vep   k_0  (1+4k_0)+    [u_0]_{\langle x\rangle^k}) 
%$$
%where $n$ is a sufficiently large number (which in particular yields $ \left(\frac{[u_0]_{\langle x\rangle^k}}{ C_1}\right)^{\frac 1{q+1}}$ to be sufficiently small). For $n$ suitably large, we see that \rife{newprep3} cannot hold, with the above choice of $C_1$. Now the choice of $C_1$ also fixes the choice of $K$ in \rife{choiK2}. 
%Putting those values together, 
Finally, we proved that \rife{claim2} holds true, where  letting $\vep\to 0$ and recalling \rife{choiC_1}--\rife{choiK2}, we deduce that
\be\label{caso1}\begin{split}
u(t,x)- u(t,y) & \leq  (1+t)^{-q} \left\{   K\, (\vfi(x)+ \vfi(y)) + C_1 \right\}
\\
&  \leq  (1+t)^{-q} M\, [u_0]_{\langle x\rangle^k} [\lgx^{\bar k}+ \lg y\rg ^{\bar k}] \qquad \forall x,y\in \R^d\,,\, t>0\,.
\end{split}
\ee
Hence  \rife{decayk} is proved.

\vskip0.5em
 {\bf (B) Case $\la_0=0$ and \rife{nu2} with $\sigma\in (0,1]$. }
If $\sigma\in (0,1)$, we use \rife{disso2} and we get from \rife{newprep-frac}
\begin{align*} 
&   -  (2Kk_0+ 2\vep k_0) \leq     (c_0(r\wedge1)^{1-\sigma+\de}+ (2\sigma_0 \wedge \sigma_1 r)\sigma_1 ) \psi'(r )  +  \frac q{(1+t)}  \psi(r)  \\
& \qquad +   4\lambda  \int_0^1 (1-s)\int_{B}   \psi'' ( r+ 2s(\widehat{x-y}\cdot z)) |\widehat{x-y}\cdot z|^2 \frac{dz}{|z|^{d+\sigma}} ds \,.
%\\
%& \leq  -   2e^{-\omega t}    \psi''(|x-y|)  
%+ 4e^{-\omega t} \, \psi'(|x-y|) (D\vfi(x)-D\vfi(y))\cdot \widehat{x-y}
% \theta\, C_1\,  C_2 e^{- C_{2} r^{\theta}} ( 1-\theta  + C_{2} r^{\theta} )   \, r^{\theta-2}    
\end{align*}
With the same arguments as in Theorem \ref{decaysup2}, with a convenient choice of $C_2$ we estimate
\begin{align*}
(c_0(r\wedge1)^{1-\sigma+\de}+ (2\sigma_0 \wedge \sigma_1 r)\sigma_1 ) \psi'(r )& + 4\lambda  \int_0^1 (1-s)\int_{B}   \psi'' ( r+ 2s(\widehat{x-y}\cdot z)) |\widehat{x-y}\cdot z|^2 \frac{dz}{|z|^{d+\sigma}} ds 
\\ & \leq 
-4\lambda(1-\theta) c_{d,\sigma}\, C_1 \,  \theta \, C_2\, e^{-C_2r_1^\theta}r_1^{\theta-\sigma}\qquad \forall r<r_1\,.
\end{align*}
Hence we get
\be\label{mmm}
  -    (2Kk_0+ 2\vep k_0) \leq      \frac q{(1+t)}  \psi(r)  -4\lambda(1-\theta) c_{d,\sigma}\, C_1 \,  \theta \, C_2\, e^{-C_2r_1^\theta}r_1^{\theta-\sigma} \,.
\ee
Here we recall the explicit choice of $K$ \rife{choiK2}, which (using the explicit value of $\psi'$) reads as
$$
K=   L \, (1+ n^{\frac q{q+1}})  [u_0]_{\langle x\rangle^k}  + 4\vep k_0    +  [c_0  + 2\sigma_0  \sigma_1 ] C_1\,C_2\,\theta \,  e^{-C_2r_1^\theta}r_1^{\theta-1}\,.
$$
Moreover, we estimate as before $\frac q{(1+t)}  \psi$ using   \rife{txy} and \rife{psiq},  and we obtain from \rife{mmm} 
\be\label{slow-frac}
\begin{split}
 &   C_1 \left\{ \theta C_2 e^{- C_{2} r_1^{\theta}} \, r_1^{\theta-1} [ 4\lambda(1-\theta) c_{d,\sigma} r_1^{1-\sigma} -  2k_0(c_0  + 2\sigma_0  \sigma_1) ]    -  \frac{\tilde C}{L^{\frac1q}(1+ n^{\frac1{q+1}})}  \right\} 
\\ & \qquad \leq  2\vep   k_0  (1+4k_0) + 2k_0\, L   \, (1+ n^{\frac q{q+1}})  [u_0]_{\langle x\rangle^k} 
\end{split}
\ee
for some $\tilde C$ independent of $C_1, C_2, r_1$. 
Now we observe that $r_1$ can be chosen sufficiently large so that 
$$
4\lambda(1-\theta) c_{d,\sigma} r_1^{1-\sigma} -  2k_0(c_0  + 2\sigma_0  \sigma_1) >0\,.
$$
Finally, having fixed $r_1$, and then $C_2$ accordingly, and recalling the choice of $C_1$ from \rife{choiC_1}, we see that \rife{slow-frac} cannot hold for $n$  sufficiently large. This concludes the argument and proves that \rife{claim2} holds true. From that, letting $\vep\to0$, we deduce \rife{caso1}, and then \rife{decayk}.
%\begin{align*}
%u(t,x)- u(t,y) & \leq  (1+t)^{-q} \left\{   K\, (\vfi(x)+ \vfi(y)) + C_1 \right\}
%\\
%&  \leq  (1+t)^{-q} M\, [u_0]_{\langle x\rangle^k} [\lgx^{\bar k}+ \lg y\rg ^{\bar k}] \qquad \forall x,y\in \R^d\,,\, t>0\,,
%\end{align*}
%where we exploited the choice of $C_1$ and $K$ that we used so far. This gives the desired estimate. 

Finally, we only mention that the case $\sigma=1$ is dealt with as in Theorem \ref{decaysup2} obtaining a slight variation of \rife{slow-frac} where $1-\sigma$ is replaced by $\eta$ for some small $\eta>0$. This is enough to conclude exactly as before.
\qed

\begin{remark} We point out that the same generalization as mentioned in Remark \ref{gendrift} on the growth condition of the drift  also applies for the content of Theorem \ref{decaysup2}.
\end{remark}

\section{Decay of Fokker-Planck equations}\label{FoPl} 

We deduce now the  decay estimates for  the general  Fokker-Planck equation \rife{FP-gen}.  Let us recall  that,  in the setting of our assumptions, problem \rife{FP-gen} admits a unique solution, according to Theorem \ref{exiuniq}.

\vskip1em
As a first step, we recall the following   characterization  of weighted seminorms.  The following lemma is  taken from \cite{HM}; we give a slightly different statement  for our purposes (and a proof for the reader's convenience).

\begin{lemma}\label{hai} Let $\vfi$ be a continuous function such that $\vfi(x) \geq 1$  for all $x\in \R^d$. For any   function $u$ such that $u\in L^\infty(\vfi^{-1})$,  we have
$$
\sup_{x,y\in \R^{d}}\,\, \frac{|u(x)-u(y)|}{ \vfi(x)+ \vfi(y)} = \inf_{c\in \R} \,\, \|u+c\|_{L^\infty(\vfi^{-1})}
$$
\end{lemma}

\proof Let us call $M:= \sup\limits_{x,y\in \R^{d}}\,\, \frac{|u(x)-u(y)|}{  \vfi(x)+ \vfi(y)}$ and $m:= \inf\limits_{c\in \R} \,\, \|u+c\|_{L^\infty(\vfi^{-1})}$. 
Since  
\begin{align*}
|u(x)-u(y)| &  \leq |u(x)+c|+ |u(y)+c| \\
& \leq   \vfi(x) \|u+c\|_{L^\infty(\vfi^{-1})}  +   \vfi(y)  \|u+c\|_{L^\infty(\vfi^{-1})}
\end{align*}
we immediately deduce that $M\leq m$. 

For the reverse inequality, define $c:= \inf_{\R^d}\,  [M  \vfi -u]$.  
Notice that $c$ is finite because $M\vfi(x) - u(x) \geq - M\, \vfi(y)- u(y)$ for any $y\in \R^d$. 
By definition, we have 
$$
\frac{u(x)+c}{  \vfi(x)} \leq M\,.
$$
We also estimate, using the definition of $M$ and $c$,
\begin{align*}
\frac{u(x)+c}{  \vfi(x)}   & = \inf_{y\in \R^d} \,\left[ \frac{u(x)+M  \vfi(y)- u(y)}{  \vfi(x)}\right]
\\
%& 
%\geq  \inf_{y\in \R^d}  \frac{ M(1+ \vfi(y))}{1+ \vfi(x)} - \frac {|u(x|-u(y)|}{2+ \vfi(x)+ \vfi(y)}\frac{2+ \vfi(x)+ \vfi(y)}{1+ \vfi(x)} 
%\\
& \geq   \inf_{y\in \R^d}  \left[ \frac{ M  \vfi(y)}{  \vfi(x)} - M \left(\frac{  \vfi(x)+ \vfi(y)}{  \vfi(x)} \right) \right]= -M\,.
\end{align*}
Hence we deduce that $\|u+c\|_{L^\infty(\vfi^{-1})}\leq M$, which yields $m\leq M$.
\qed

\vskip1em

Thanks to Lemma \ref{hai}, the decay estimates follow by duality using the results of the previous Section.

%\be\label{linfp}
%\begin{cases}
%\partial_t m - \mathcal{I} ( x, [ m ] ) -\dive( b(t,x) m)     =  0 &  \text{ in } ( 0,T ) \times \R^d, \\ 
%m( 0,x ) = m_0 (x),\qquad  \,&  \text{ in }   \R^d
%  \end{cases}
%\ee

\begin{theorem}\label{decaym} 
Let $\cL$ be defined by \rife{L} where $\Sigma$ satisfies  \rife{sigma},  the Levy measure in $\cI$ satisfies \rife{nu2}, and $\lambda_0+\lambda>0$.  Assume that $b\in C^0(Q_T)$ satisfies \rife{bdiss} 
and either of the following conditions:
\vskip0.4em
(i)  $\la_0>0$ or $\sigma\in (1,2)$, and  \rife{disso} holds true.
\vskip0.4em
(ii)  $\la_0=0$, $\sigma\in (0,1]$   and   \rife{disso2} holds true.
\vskip0.4em
Assume that  $m_0\in  \cM_k(\R^d)$, for some $0<k<\sigma$, and that  $\intd dm_0=0$, and let $m$ be the unique solution of \rife{FP-gen} (in the sense of Definition \ref{def-dual}).  Then we have:
  \begin{itemize}
  
\item[(a)] If $\gamma\geq 2$ in \rife{bdiss}, then 
\be\label{expm}
\| m(t)\|_{\cM_k} \leq K \, e^{-\omega t} \, \|m_0\|_{\cM_k}
\ee
for some constant $K,\omega$ only depending on $\alpha, \sigma, (\lambda_0+\lambda)^{-1}, \sigma_0,\sigma_1,  \Lambda, \de, d, c_0, k$.  

\item[(b)] If $\gamma \in (\max((2-\sigma),1), 2)$ in \rife{bdiss}, then for any $k<\bar k<\sigma$  such that $m_0\in \cM_{\bar k} $ we have 
\be\label{polym}
\| m(t)\|_{\cM_k} \leq K \,  (1+t)^{-q}  \, \|m_0\|_{\cM_{\bar k}}\qquad \hbox{where $q=  \frac{\bar k-k}{2-\gamma}$}
\ee
for some $K$ depending on $\alpha, \gamma,  \sigma, (\lambda_0+\lambda)^{-1},  \sigma_0,\sigma_1, \Lambda, \de, d, c_0, k, \bar k$.  
\end{itemize}
\end{theorem}

\proof  For $t>0$ we consider the unique viscosity solution of the adjoint equation 
\be\label{adj}
\begin{cases}
-\partial_t u + \cL[ u ]  + b(t,x) \cdot Du    =  0 &  \text{ in } ( 0,t ) \times \R^d, \\ 
u( t,x ) = \xi,\qquad  \,&  \text{ in }   \R^d
  \end{cases}
\ee
where $\xi \in C(\R^d)\cap L^{\infty}(\langle x\rangle^{-k})$.  By Definition \ref{def-dual} $m$ satisfies  
\begin{align*}
\int_{\R^d} \xi \, dm(t)   = \int_{\R^d}  u(0) \, dm_0& = \int_{\R^d}   (u(0)+c)\,  dm_0 
\\
& \leq \|m_0\|_{\cM_k}\, \|u(0)+ c\|_{L^{\infty}(\langle x\rangle^{-k})} 
\end{align*}
where we used that $m_0$ has zero average. Taking the infimum with respect to $c\in \R$ and using Lemma \ref{hai} we get  
\be\label{dauam}
\int_{\R^d} \xi \, dm(t)  \leq \|m_0\|_{\cM_k}\,  [u(0)]_{\langle x\rangle^k}\,.
\ee
If $\gamma\geq 2$, we estimate $[u(0)]_{\langle x\rangle^k}$ from  Theorem \ref{decaysup2} and we obtain
\begin{align*}
\int_{\R^d} \xi \, dm(t)     & \leq K\, e^{-\omega t} \|m_0\|_{\cM_k}\,  [\xi]_{\langle x\rangle^k}  
\\
& \leq K\, e^{-\omega t} \|m_0\|_{\cM_k}\,  \|\xi\|_{L^{\infty}(\langle x\rangle^{-k})} \,.
\end{align*}
Since $\xi$ is arbitrary, we deduce  \rife{expm}. 
\vskip0.4em
If $\gamma\in (1,2)$ with $\gamma>2-\sigma$, we argue similarly but rather than \rife{dauam} we estimate as follows:
$$
\int_{\R^d} \xi \, dm(t)   \leq \|m_0\|_{\cM_{\bar k}}\,  [u(0)]_{\langle x\rangle^{\bar k}}\,
$$
for $\bar k>k>2-\gamma$ (notice that $u\in L^{\infty}(\langle x\rangle^{-k})$ implies $u\in L^{\infty}(\langle x\rangle^{-\bar k})$).  Applying now Theorem \ref{decayslow_frac} we obtain
$$
\int_{\R^d} \xi \, dm(t)   \leq \|m_0\|_{\cM_{\bar k} }\,   K \, (1+t)^{-q} \,  \|\xi\|_{L^{\infty}(\langle x\rangle^{-k})} \qquad \hbox{where $q=  \frac{\bar k-k}{2-\gamma}$}
$$
which yields \rife{polym}.
\qed

\begin{remark} 
The previous result for Levy operators includes the particular case stated for the fractional Laplacian   in the Introduction (Theorem \ref{nonlocalthm}). We only need to observe that, if the drift $b$ is locally Lipschitz in $x$ (or alternatively, under milder conditions involving only $\dive(b)$),   the weak solution is unique (so it coincides with the unique solution given by Theorem \ref{exiuniq}) and  belongs to $L^1$ (see e.g. \cite{CiGo}, \cite{Lafleche}, \cite{Wei-Tian}).
\end{remark}

\vskip1em
In the same spirit of the above remark, the solutions are $L^1$ whenever the operator is second order. In addition, if the operator is local, one can consider the full range of $\gamma$, including   $\gamma\in (0,1)$. The next result  extends Theorem \ref{localthm} which was  stated in the Introduction for the case of Laplace operator. The proof is the same as above, using duality and Theorem \ref{decaysup} or Theorem \ref{degconf}. This extends  results as  \cite[Thm 3]{tovi}, or \cite[Thm 1.2]{Kavian+al}.

\begin{theorem}\label{decaym-local} 
Let  $m_0\in L^1( \lgx^k)$ be such that $\intd m_0\,dx=0$.  Let  $m$ be  a weak solution to the Cauchy problem 
\be\label{cauchyL0}
\begin{cases}
\partial_t m + \cL_0^*[m]   - \dive( b(t,x) m)     =  0 &  \text{ in } ( 0,T ) \times \R^d, \\ 
m ( 0,x ) = m_0 (x),\qquad  \,&  \text{ in }   \R^d
\end{cases}
\ee
where we assume that $\lambda_0>0$ and \rife{sigma} holds true, and $b(t,x)$ is a continuous vector field  satisfying \rife{bdiss} and   \rife{disso}.
\vskip0.4em
Then we have:
  \begin{itemize}
  
\item[(a)] If $\gamma\geq 2$ in \rife{bdiss}, then 
\be\label{expm0}
\| m(t)\|_{L^1(\lg x\rg^k)} \leq K \, e^{-\omega t} \, \|m_0\|_{L^1({\langle x\rangle^k})}
\ee
for some constant $K,\omega$ only depending on $\alpha,   \lambda_0^{-1}, \sigma_0,\sigma_1,    d, c_0, k$.  

\item[(b)] If $\gamma \in (0, 2)$ in \rife{bdiss}, then for any $k<\bar k<\sigma$  such that $m_0\in L^1( \lgx^{\bar k})$ we have 
\be\label{polym0}
\| m(t)\|_{L^1(\lg x\rg^k)} \leq K \,  (1+t)^{-q}  \, \|m_0\|_{L^1({\langle x\rangle^{\bar k}})}\qquad \hbox{where $q=  \frac{\bar k-k}{2-\gamma}$.}
\ee
for some $K$ depending on $\alpha, \gamma,   \lambda_0^{-1}, \sigma_0,\sigma_1,  d, c_0, k, \bar k$.  
\end{itemize}
\end{theorem}

\section{Generalizations and further results}

We presented so far  the results that can be obtained using  Lyapunov functions of power type, resulting in moments decay of measure-valued solutions. Here we review our strategy in a more abstract form, showing that it can  be applied to different choices of Lyapunov functions. The Lyapunov function will determine the weighted space for the oscillation estimate and, what is more important, the rate of convergence is consequential to the weighted norm which is used. 

There are two main settings that we consider, corresponding to exponential or sub-exponential decay.  Here and below, $\vfi$ denotes a generic positive smooth function such that $\vfi(x) \to \infty$ as $|x|\to \infty$. There is no loss of generality in assuming that $\vfi(x)\geq 1$ in $\R^d$.
We first consider the simplest case of exponential decay in {\it some} weighted  norm, namely when some $\vfi$ exists satisfying   
\be\label{H1}
\tag{H1}
\exists \,\, \omega_0>0\,:\quad \liminf_{|x|\to \infty} \, \frac{\cL^b[\vfi]}{\vfi}\geq \omega_0\,.
\ee
We denote 
$$
\Lambda_1:= \{\vfi\in C^2(\R^d):\,\, \vfi(x)\mathop{\longrightarrow}^{|x|\to \infty} \infty, \,\, \hbox{$\cI(x,[\vfi])$ is well defined and $\vfi$ satisfies \rife{H1}}\}.
$$ 
Clearly, for every $\vfi\in \Lambda_1$, one can choose $\omega$ sufficiently small so that
$$
\cL^b[\vfi]-\omega\vfi \geq \frac{\omega_0}2 \vfi - k_0 \qquad \forall x\in \R^d
$$
for some $k_0>0$. Using the above property, the proof of Theorem \ref{decaysup2} follows replacing $\lgx^k$ with $\vfi(x)$, so that we can generalize the exponential decay estimate for the oscillation, as follows.
%Then we can generalize the have the following oscillation decay estimate.

\begin{theorem}\label{decay-gen1} Let $\cL$ be defined by \rife{L} where $\Sigma$ satisfies  \rife{sigma},  the Levy measure in $\cI$ satisfies \rife{nu2}, and $\lambda_0+\lambda>0$.  Assume that $b$ satisfies \rife{bdiss}  and either of the following conditions:

(i) $\la_0>0$ or  $\sigma\in (1,2)$, and \rife{disso} holds true.

(ii) $\la_0=0$ and $\sigma\in (0,1]$, and \rife{disso2} holds true.

Assume that there exists $\vfi\in \Lambda_1$.  Let $u$ be a  (viscosity) solution of \rife{pbgen} such that, for some $\Phi\in \Lambda_1$, one has  $|u(t,x)|= o(\Phi) $ as $|x|\to \infty$, uniformly for $t\in (0,T)$.   

Then there exist positive constants $K, \omega >0$ (depending on $\alpha, \gamma, \sigma, (\lambda+ \lambda_0)^{-1}, \sigma_0, \sigma_1, c_0, d,  \vfi$)  such that $u$ satisfies
\be\label{decayvfi}
[u(t)]_{\vfi} \leq K \, e^{-\omega t} \, [u_0]_{\vfi}.
\ee
\qed
\end{theorem}

%\proof The proof is just another    version of Theorem \ref{decaysup2}. The claim follows by proving that, whatever small $\vep$,  the function
%$$
%u(t,x)- u(t,y) -  e^{-\omega t}    \underbrace{\left\{  K  [\vfi(x)+ \vfi(y) ]    +\psi(|x-y|)  + \vep [\Phi(x)+ \Phi(y)] \right\}}_{\zeta(x,y)} -  \frac\vep{T-t} 
%$$
%cannot have a positive maximum. The function $\psi$ is defined in \rife{proppsi}.  Notice that the existence of a maximum on $[0,T]\times \R^d$ is ensured by the assumption that $|u(t,x)|= o(\Phi) $ as $|x|\to \infty$, uniformly for $t\in (0,T)$. Since $\vfi, \Phi\in \Lambda$, one can choose $\omega$ sufficiently small so that
%$$
%\cL^b[\vfi]-\omega\vfi \geq \frac{\omega_0}2 \vfi - k_0 \qquad ;\qquad \cL^b[\Phi]-\omega\Phi\geq   - k_0
%$$
%for some $k_0>0$. This is the only property used in Theorem \ref{decaysup2}, and the other arguments of the proof remain unchanged.
%\qed
\vskip1em

The case of sub-exponential growth is more subtle. This corresponds to the case that the  Lyapunov function $\vfi$ satisfies a  condition as follows:
\be\label{H2}
\tag{H2}
\begin{split} \cL^b[\vfi] \mathop{\to}^{|x|\to \infty} \infty \,,\quad & \hbox{and $\,\,\exists\,\,$ a decreasing function $h:\R_+\to \R_+$, with $h(r)\mathop{\to}\limits^{r\to \infty} 0$\,\,:}\quad \, 
 \\  &  \cL^b[\vfi]\geq h(\vfi) \vfi  \qquad \hbox{as $|x|\to \infty$\,.} 
 \end{split}
\ee
Correspondingly, we set
$$
\Lambda_2:= \{ \vfi\in C^2(\R^d):\,\, \vfi(x)\mathop{\longrightarrow}^{|x|\to \infty} \infty, \,\, \hbox{$\cI(x,[\vfi])$ is well defined and $\vfi$ satisfies \rife{H2}}\}.
$$
As in the estimate \rife{decayk}, the sub-exponential decay is observed only when   using different weights in the seminorms of $u(t)$ and $u(0)$. One possible choice, among others,  is to use different powers of $\vfi$, as we did in Theorem \ref{decayslow_frac}. 
%A technical condition (we will require that $\frac{|D\vfi|^2}{\vfi}=o(b\cdot D\vfi)$ as $|x|\to \infty$) is  used here to ensure that different powers of $\vfi$ are still Lyapunov functions.

\begin{theorem}\label{decaygen2}
Let $\cL$ be defined by \rife{L} where $\Sigma$ satisfies  \rife{sigma},  the Levy measure in $\cI$ satisfies \rife{nu2}, and $\lambda_0+\lambda>0$.  Assume that $b$ satisfies \rife{bdiss} with $\gamma \in (0, 2)$,   and either of the following conditions:

(i) $\la_0>0$ or  $\sigma\in (1,2)$, and \rife{disso} holds true.

(ii) $\la_0=0$ and $\sigma\in (0,1]$, and \rife{disso2} holds true.

Assume that there exists $\vfi$ satisfying \rife{H2}. Let $u$ be a  (viscosity) solution of \rife{pbgen} such that  $|u(t,x)|= o(\Phi) $ as $|x|\to \infty$ (uniformly for $t\in (0,T)$) for some $\Phi\in C^2(\R^d)$ such that $\Phi(x)\mathop{\to}\limits^{|x|\to \infty} \infty$,  $\cL^b[\Phi]\mathop{\to}\limits^{|x|\to \infty} \infty$ and $\liminf_{|x|\to \infty}  [\cL^b[\Phi]-\kappa\, h(\vfi)\Phi ]> -\infty$ for some $\kappa>0$.   

 Then, for any $\theta<1 $ such that $ h(\vfi)\vfi^\theta \to \infty$ as $|x|\to \infty$, there exist      positive constants  $L,M>0$ (depending on $\alpha, \gamma, \sigma, (\lambda_0+\lambda)^{-1},   \sigma_0, \sigma_1, c_0, d,  \vfi, \theta$)  such that $u$ satisfies
\be\label{decayvfi2}
[u(t)]_{\vfi } \leq M \, \varpi(t) \, [u_0]_{\vfi^\theta}\qquad \hbox{where $\varpi(t)$ is implicitly defined as  $\int_\varpi^1 \frac{ds}{sh(L\,s^{-\frac1{1-\theta}})}= \frac t2$,}
\ee
with the function $h$  given by \rife{H2}.
\end{theorem}

\begin{remark} The interesting part of \rife{decayvfi2} is that the decay rate $\varpi(t)$ is given in terms of the condition \rife{H2} involving the Lyapunov function $\vfi$. 
Notice that $\varpi$ is sub-exponential precisely because $h(r)\to 0$ as $r\to \infty$ ($h$ gives a quantification of the fact that $\frac{\cL^b[\vfi]}\vfi \mathop{\to}\limits^{|x|\to \infty} 0$).  The implicit condition for $\varpi(t)$ comes from the requirement that $\varpi$ satisfies the ODE
\be\label{ode}
\begin{cases}
\varpi'= -\frac 12\varpi\, h(L\varpi^{-\frac1{1-\theta}})\,,\,\, t>0 & \\
\varpi(0)= 1 
\end{cases}
\ee
Here, we did not look for sharp conditions in terms of $\theta$, and the choice of the constant $L$ is not meant to be optimal, either. Thus, there could be room for sharper conditions. However,  the  rate function which stems from \rife{ode} looks consistent with all previous results; in particular,  in the range of power functions $\vfi= \lg x\rg^k$, the above statement reduces to Theorem \ref{decayslow_frac}.
\end{remark}

\begin{remark} Despite the technical assumption appearing in the statement, the growth of solutions $u$ plays a minor role in the above result. One could readily use the statement with $\Phi= \vfi$, and apply the estimate to any bounded solution  $u$. This is how we are going to use this result later.
\end{remark}

\proof The proof is a (slightly more abstract) version of Theorem \ref{decayslow_frac}. 
%We first notice that, requiring $\frac{|D\vfi|^2}{\vfi}=o(b(t,x)\cdot D\vfi)$ as $|x|\to \infty$ (uniformly in time) lets $\vfi^p$ be a Lyapunov function itself. Indeed, we have
%$$
%D^2 \vfi^p = p \vfi^{p-1} \left(D^2 \vfi+ (p-1) \frac{D\vfi\otimes D\vfi}\vfi\right) \quad \rightsquigarrow \quad \cL_0[\vfi^p] = p \vfi^{p-1} \left(\cL_0[\vfi] + O\left(\frac{|D\vfi|^2}{\vfi}\right)\right)
%$$
%and
%$$
%\cI(x,[\vfi^p])  =  \int_{|z|>1} \{\vfi^p ( x+z ) -\vfi^p( x )  \}\nu ( dz ) +  \int_{|z|\leq 1} \{\vfi^p( x+z ) -\vfi^p( x ) -  p \vfi^{p-1} (D\vfi( x ) \cdot  z ) \}\nu ( dz )
%$$

As a first step, we observe that, for all  functions $\vfi\in C^2(\R^d)$ such that $\cL^b[\vfi]\to \infty$ as $|x|\to \infty$,  it holds  
$$
 [u(t)]_{\vfi}\leq C \,  [u_0]_{\vfi}
$$
for  a   suitable constant $C>0$. This is obtained exactly as \rife{osci} in Theorem \ref{decayslow_frac}. Now, if $\theta<1$, we have
$$
D^2 \vfi^\theta =  \theta \vfi^{\theta-1} \left(D^2 \vfi+ (\theta-1) \frac{D\vfi\otimes D\vfi}\vfi\right) \leq  \theta \vfi^{\theta-1}  D^2 \vfi\
%quad \rightsquigarrow \quad \cL_0[\vfi^p] = p \vfi^{p-1} \left(\cL_0[\vfi] + O\left(\frac{|D\vfi|^2}{\vfi}\right)\right)
$$
and, by concavity,
\begin{align*}
\cI(x,[\vfi^\theta])  & =  \int_{|z|>1} \{\vfi^\theta ( x+z ) -\vfi^\theta( x )  \}\nu ( dz ) +  \int_{|z|\leq 1} \{\vfi^\theta( x+z ) -\vfi^\theta( x ) -  \theta \vfi^{\theta-1} (D\vfi( x ) \cdot  z ) \}\nu ( dz )
\\ & \leq  \theta \vfi(x)^{\theta-1}\cI(x,[\vfi ]) \,.
\end{align*}
Hence we have
$$
\cL^b[\vfi^\theta] \geq \theta \vfi^{\theta-1} \cL^b[\vfi] \geq \theta \vfi^\theta\, h(\vfi) \to \infty\qquad \hbox{as $|x|\to \infty$}
$$
by assumption. This implies that $\vfi^\theta$ is itself a Lyapunov function and we have,  by what proved before, $ [u(t)]_{\vfi^\theta}\leq C \,  [u_0]_{\vfi^\theta}$. 
Now we go back to the maximum point of 
$$
u(t,x)- u(t,y) -  \varpi(t)   \underbrace{\left\{  K  [\vfi(x)+ \vfi(y) ]    +\psi(|x-y|)  + \vep [\Phi(x)+ \Phi(y)] \right\}}_{\zeta(x,y)} -  \frac\vep{T-t} 
$$
where $K\geq [u_0]_{\vfi^\theta} \geq [u_0]_{\vfi }$ and $\varpi$ satisfies \rife{ode}. On the maximum point we have, similarly as in \rife{newprep-frac}:
\begin{align*} 
& \frac\vep{(T-t)^2} +  \varpi(t)   \, K  \left\{  \mathcal{L}^b[\vfi](x)  +   \frac{\varpi' }\varpi \vfi(x)  +  \mathcal{L}^b[\vfi](y)  +    \frac{\varpi'}{\varpi} \vfi(y)\right\}  \\
&   +   \varpi (t) \, \vep   \left( \mathcal{L}^b[\Phi](x)+ \frac{\varpi'}{\varpi}\Phi(x) + \mathcal{L}^b[\Phi](y)+  \frac{\varpi'}{\varpi} \Phi(y) \right)
\\
& + \varpi(t)   \left((b(t,x)-b(t,y)) \cdot  \widehat{x-y}\right)  \psi'(r )  + \varpi'(t) \psi(r) \\ & \qquad \leq  \varpi(t)  \left\{4\lambda_0 \psi''(r)  + 
 (2\sigma_0\wedge \sigma_1 r) \sigma_1 \,\psi'(r)  \right\} 
\\
& \qquad + \varpi(t) 4\lambda  \int_0^1 (1-s)\int_{B}   \psi'' ( r+ 2s(\widehat{x-y}\cdot z)) |\widehat{x-y}\cdot z|^2 \frac{dz}{|z|^{d+\sigma}} ds \, 
%\\
%& \leq  -   2e^{-\varpi t}    \psi''(|x-y|)  
%+ 4e^{-\varpi t} \, \psi'(|x-y|) (D\vfi(x)-D\vfi(y))\cdot \widehat{x-y}
% \theta\, C_1\,  C_2 e^{- C_{2} r^{\theta}} ( 1-\theta  + C_{2} r^{\theta} )   \, r^{\theta-2}    
\end{align*}
where $r=|x-y|$. We consider the case that $r= |x-y| >r_1$, for a  sufficiently large $r_1$. In particular, we have $|x|\vee |y|>\frac{r_1}2$.  
We drop the terms with $\psi''$, we use \rife{disso} (or \rife{disso2}) for the drift term, and  we obtain
\begin{align*} 
&  K  \left\{  \mathcal{L}^b[\vfi](x)  +   \frac{\varpi' }\varpi \vfi(x)  +  \mathcal{L}^b[\vfi](y)  +    \frac{\varpi'}{\varpi} \vfi(y)\right\} 
   +   \vep   \left( \mathcal{L}^b[\Phi](x)+ \frac{\varpi'}{\varpi}\Phi(x) + \mathcal{L}^b[\Phi](y)+  \frac{\varpi'}{\varpi} \Phi(y) \right)
\\
&  \qquad  +  \frac{\varpi'}{\varpi} \psi(r)  < 
 ( c_0 + 2\sigma_0 \sigma_1) \,\psi'(r)   \,.
%\\
%& \leq  -   2e^{-\varpi t}    \psi''(|x-y|)  
%+ 4e^{-\varpi t} \, \psi'(|x-y|) (D\vfi(x)-D\vfi(y))\cdot \widehat{x-y}
% \theta\, C_1\,  C_2 e^{- C_{2} r^{\theta}} ( 1-\theta  + C_{2} r^{\theta} )   \, r^{\theta-2}    
\end{align*}
Using $\psi\leq C_1$ and $\psi'(r)\leq \psi'(r_1)$ we get
\begin{align*} 
&  K  \left\{  \mathcal{L}^b[\vfi](x)  +   \frac{\varpi' }\varpi \vfi(x)  +  \mathcal{L}^b[\vfi](y)  +    \frac{\varpi'}{\varpi} \vfi(y)\right\}
%\left\{  \left( h(\vfi)  +   \frac{\varpi' }\varpi \right) \vfi(x)  +  \left( h(\vfi) +    \frac{\varpi'}{\varpi}\right)  \vfi(y)\right\}   
   +   \vep   \left( \mathcal{L}^b[\Phi](x)+ \frac{\varpi'}{\varpi}\Phi(x) + \mathcal{L}^b[\Phi](y)+  \frac{\varpi'}{\varpi} \Phi(y) \right)
\\
&  \qquad    <  
 ( c_0 + 2\sigma_0 \sigma_1) \,\psi'(r_1)  -  \frac{\varpi'}{\varpi} C_1 \,.
%\\
%& \leq  -   2e^{-\varpi t}    \psi''(|x-y|)  
%+ 4e^{-\varpi t} \, \psi'(|x-y|) (D\vfi(x)-D\vfi(y))\cdot \widehat{x-y}
% \theta\, C_1\,  C_2 e^{- C_{2} r^{\theta}} ( 1-\theta  + C_{2} r^{\theta} )   \, r^{\theta-2}    
\end{align*}
We observe that, since the maximum is positive, we have
$$
K \varpi(t)[\vfi(x)+ \vfi(y) ]  \leq u(t,x)- u(t,y) \leq [\vfi^\theta(x)+ \vfi^\theta(y) ] [u(t)]_{\vfi^\theta}\leq C \, [\vfi^\theta(x)+ \vfi^\theta(y) ] \, [u_0]_{\vfi^\theta}
$$
which yields
$$
\varpi(t) \leq C \frac{[u_0]_{\vfi^\theta}}K\, \frac{\vfi^\theta(x)+ \vfi^\theta(y)}{\vfi(x)+ \vfi(y)} \leq C \frac{[u_0]_{\vfi^\theta}}K\, \frac1{(\vfi(x)\vee \vfi(y))^{1-\theta}}\,.
$$
In particular, since $h$ is decreasing, we have
$$
h(\vfi) \geq h\left( \left( C \frac{[u_0]_{\vfi^\theta}}K\right)^{\frac1{1-\theta}} \varpi^{-\frac1{1-\theta}}\right) \,.
$$
For a constant $L$ to be chosen below, such that $L\geq \left( C \frac{[u_0]_{\vfi^\theta}}K\right)^{\frac1{1-\theta}}$, we fix $\varpi$ as the solution of \rife{ode}. Then
$$
\frac{\varpi'}{\varpi} = -\frac12 h \left( L\varpi^{-\frac1{1-\theta}}\right) \geq - \frac12 h(\vfi)
$$
which implies
$$
\mathcal{L}^b[\vfi]   +   \frac{\varpi' }\varpi \vfi \geq     \frac 12 h(\vfi) \vfi \to \infty \quad \hbox{as $|x| \to \infty$.}
$$
Since  $|x|\vee |y|>\frac{r_1}2$, choosing $r_1$ large enough we have that
$$
 \mathcal{L}^b[\vfi](x)  +   \frac{\varpi' }\varpi \vfi(x)  +  \mathcal{L}^b[\vfi](y)  +    \frac{\varpi'}{\varpi} \vfi(y) \geq 1
$$
for every $x,y,$: $|x-y|> r_1$.  Thus we deduce that 
\be\label{fii}\begin{split}
& 
K   
   +   \vep   \left( \mathcal{L}^b[\Phi](x)- \frac12 h(\vfi(x) ) \Phi(x) + \mathcal{L}^b[\Phi](y)- \frac12 h(\vfi(y) ) \Phi(y) \right)
%\\  & \qquad\qquad\qquad     
<  
 ( c_0 + 2\sigma_0 \sigma_1) \,\psi'(r_1)  -  \frac{\varpi'}{\varpi} C_1\,.
\end{split}
\ee
Notice that $\varpi\leq 1$ implies $h \left( L\varpi^{-\frac1{1-\theta}}\right)\leq h(L)$, hence $-  \frac{\varpi'}{\varpi} C_1\leq \frac12 h(L)C_1$.
As for the function $\Phi$, we know by assumption that $\liminf_{|x|\to \infty}  \cL^b[\Phi]-\kappa h(\vfi)\Phi > -\infty$, for some $\kappa>0$ which we can assume, without loss of generality, such that $\kappa\geq \frac12$.  Hence
$\mathcal{L}^b[\Phi](x)- \frac12 h(\vfi(x) ) \Phi(x) \geq -k_0$ for some constant $k_0>0$.
Then, 
%since $h(\vfi)\vfi \to \infty$ as $|x|\to \infty$, choosing $r_1$ sufficiently large the   
inequality \rife{fii} implies
$$
K<  2k_0 \vep+ ( c_0 + 2\sigma_0 \sigma_1) \,\psi'(r_1)  + \frac12 h(L) C_1\,.
$$
We conclude that the maximum point cannot be reached with $r>r_1$ provided 
$$
K\geq  2k_0 \vep+ ( c_0 + 2\sigma_0 \sigma_1) \,\psi'(r_1)  + \frac12 h(L) C_1\,.
$$
This replaces \rife{choiK2} in Theorem \ref{decayslow_frac}. Once $K$ is chosen, the rest of the proof is the same, by simply using that $\mathcal{L}^b[\vfi], \mathcal{L}^b[\Phi]\geq -k_0$ for some $k_0>0$, and the conclusion is obtained by choosing $L$ and $C_1$  sufficiently large.
\qed

As in Section \ref{FoPl}, we deduce by duality the decay estimate for Fokker-Planck equations.

\begin{theorem}\label{decayFP2} Let $\cL$ be defined by \rife{L} where $\Sigma$ satisfies  \rife{sigma},  the Levy measure in $\cI$ satisfies \rife{nu2}, and $\lambda_0+\lambda>0$.  Assume that $b$ satisfies \rife{bdiss}  and either of the following conditions:

(i) $\la_0>0$ or  $\sigma\in (1,2)$, and \rife{disso} holds true.

(ii) $\la_0=0$ and $\sigma\in (0,1]$, and \rife{disso2} holds true.

For $m_0\in \cM(\R^d)$ with $\intd \,dm_0=0$, let $m$ be the  solution of \rife{FP-gen}. Then we have:
\begin{itemize}

\item[(a)] If $\vfi\in \Lambda_1$ and $m_0 \in \cM_\vfi$, then 
$$
\|m(t)\|_{\cM_\vfi} \leq K \, e^{-\omega\, t} \|m_0\|_{\cM_\vfi}
$$
for some $K,\omega$ depending on $\alpha, \gamma, (\lambda_0+\lambda)^{-1}, \sigma, \sigma_0, \sigma_1, c_0,  \delta, \vfi, d$.

\item[(b)] If $\vfi\in \Lambda_2$ and $m_0 \in \cM_\vfi$, then for all $\theta\in (0,1)$ such that   $ h(\vfi)\vfi^\theta \to \infty$ as $|x|\to \infty$  we have
% there exist      positive constants  $L,M>0$ (depending on $\alpha, \gamma, \sigma, (\lambda_0+\lambda)^{-1},   \sigma_0, \sigma_1, c_0, d,  \vfi, \theta$)
$$
\|m(t)\|_{\cM_{\vfi^\theta} }\leq K\, \varpi(t) \|m_0\|_{\cM_\vfi}
$$
where $\varpi$ is defined as in \rife{decayvfi2}, for some $K,L$ depending on $\theta, \alpha, \gamma, (\lambda_0+\lambda)^{-1}, \sigma, \sigma_0, \sigma_1, c_0,  \delta, \vfi, d$.
\end{itemize}
\qed
\end{theorem}

As an example of application of Theorem \ref{decayFP2}, when we consider the Lyapunov function
\be\label{expmu}
\vfi= e^{\mu \lg x\rg^k}
\ee
with parameters $\mu, k>0$,  we extend to the local operator $\cL_0$  similar results obtained  in the recent paper \cite{Kavian+al}  for the Laplace operator.  
%This  confirms, once more, the  efficiency of our approach. 

\begin{corollary}\label{kavian}  Let $\cL_0$ be defined by \rife{L0}, where we assume $\lambda_0>0$ and \rife{sigma}. Assume that $b(t,x)$ is   continuous and satisfies \rife{bdiss} and   \rife{disso}. 

Let $\vfi$  be  defined by \rife{expmu}, where we assume that either $k<\gamma$ (and $\mu>0$ is arbitrary) or $k=\gamma$ and $\mu<\frac\alpha\gamma$ ($\alpha$ appears in \rife{bdiss}).  
Let  $m_0\in L^1(\vfi)$ be such that $\intd m_0\,dx=0$, and let  $m$ be  the   solution to the Cauchy problem \rife{cauchyL0}.
Then we have:
  \begin{itemize}
  
\item[(a)] If $\gamma\geq 2$,  then 
\be\label{expmmu}
\| m(t)\|_{L^1(\vfi)} \leq K \, e^{-\omega t} \, \|m_0\|_{L^1(\vfi)}
\ee
for some constant $K,\omega$ only depending on $\alpha,   \lambda_0, \sigma_0,\sigma_1,  d, c_0, k, \mu$.  

\item[(b)] If $\gamma \in [1, 2)$, then the exponential decay \rife{expmmu} holds for $k\geq 2-\gamma$. 
%some $k,\mu$ (precisely, for $k  \in [2-\gamma, \gamma)$ or $k=\gamma$ and $\mu$ sufficiently small). 

\item[(c)] If $\gamma \in (0, 1)$, for any $\theta\in (0,1)$ we have  
\be\label{decaytheta}
\|m(t)\|_{L^1(\vfi^\theta)} \leq  K\, e^{- C\, t^{\frac k{2-\gamma}}}\, \|m_0\|_{L^1(\vfi)}\,
\ee 
for some $K, C$ depending on $\theta, \alpha, \gamma,   \lambda_0, \sigma_0,\sigma_1,  d, c_0, k, \mu$.  
\end{itemize}
\end{corollary}

\proof
We observe the following, through standard computations:

(i) $\lim\limits_{|x|\to \infty} \cL_0[\vfi] = \infty$ provided $k<\gamma$ or $k=\gamma$ and $\mu <\alpha/\gamma$ ($\alpha$ appears in \rife{bdiss})

(ii) if $\gamma\geq 2$, we have $\liminf\limits_{|x|\to \infty} \frac{\cL_0[\vfi] }\vfi >0 $. Hence $\|m(t)\|_{L^1(\vfi)}$ decays exponentially, according to Theorem \ref{decayFP2}.

(iii) if $\gamma \in [1, 2)$,  $\liminf\limits_{|x|\to \infty} \frac{\cL_0[\vfi] }\vfi >0 $ provided $k\in [2-\gamma, \gamma)$ (or even $k=\gamma$ if $\mu<\frac\alpha\gamma$)

(iv) If $\gamma\in (0,1)$, then we have $ \cL_0[\vfi]\to \infty$ but $\cL_0[\vfi]\simeq h(\vfi)\vfi$ for $h= O \left(\frac1{(\log \vfi)^{\frac{2-\gamma}k-1}}\right)$. Following the recipe of Theorem \ref{decayFP2}, we deduce that, for any $\theta\in (0,1)$, there exists some $C, K$ (depending on $\theta$) such that \rife{decaytheta} holds true.
\qed

In particular, the above  result shows that, in the local case,  there is exponential decay of solutions, even for $\gamma\in [1,2)$, provided one uses a suitable weighted norm.  Even for a milder confining property of the drift, say $\gamma\in (0,1)$, compared to Theorem \ref{decaym-local}  we have proved that the solution decays much faster than polynomially when considered in the above weighted norm. 

Let us observe that,  in view of assumption \rife{nu2}, the Lyapunov function \rife{expmu} grows too fast at infinity for the nonlocal part, this is why the above result was  restricted to the case of local operator. However, by simply modifying the condition at infinity on the Levy measure $\nu$, the reader can readily generalize the above statement to other classes of Levy operators, with fast decay of the kernel at infinity.

\subsection{Convergence to stationary solutions}  We stress that the above results can also be used to provide existence of stationary solutions, namely invariant measures  
of Fokker-Planck equations. We give here a sample of this  application.

 \begin{theorem} 
Under the assumptions of Theorem \ref{decayFP2}, assume that  either $\vfi\in \Lambda_1$, or $\vfi\in \Lambda_2$  and $\varpi$ defined in \rife{decayvfi2} satisfies $\varpi\in L^1(0,\infty)$.  Assume that $b=b(x)$ is a continuous, time-independent vector field. 
% Let $m_0$ be a probability measure such that $m_0\in L^1(\R^d, \lgx^k\,dx)$. Assume that $b=b(x)$ is a continuous, time-independent vector field.  Under the assumptions of Theorem \ref{decaym}, assume that problem \rife{FP-gen} admits a solution $m\in L^\infty(0,T; L^1(\R^d, \lgx^k\,dx))$, for $k<\sigma$.  
 Then there exists a unique stationary measure $\bar m \in \cP(\R^d)$, which is a solution of
\be\label{FP-stat}
%\begin{cases}
 \cL^*[\bar m]   - \dive( b(x) \bar m)     =  0  \qquad \text{ in }   \R^d,
% & \text{ in }   \R^d, \\ 
%\bar m \in L^1(\R^d, \lgx^k\,dx)\,,\quad \intd \bar m\, dx =1,\qquad  \,&  
%\end{cases}
\ee 
and such that,  
%$\bar m \in \cM_\vfi$. In addition, 
if $m_0\in \cP(\R^d)\cap \cM_\vfi$ and $m$ is the solution of \rife{FP-gen}, we have, for some $K>0$, 
%in the sense of Definition \ref{def-dual}. Moreover, we have
\begin{eqnarray}
& \hbox{if $\vfi\in \Lambda_1$, then} \quad  
\| m(t)-\bar m\|_{ \cM_\vfi} \leq K \, e^{-\omega t} \,
% \|m_0-\bar m\|_{\cM_\vfi}\,, 
& \label{stat1}
\\
& \hbox{if $\vfi\in \Lambda_2$, then} \quad  \| m(t)-\bar m\|_{\cM_{\vfi^\theta}} \leq  K\, \varpi(t) \, 
%\|m_0-\bar m\|_{\cM_\vfi} & \label{stat2}
\end{eqnarray}
for any $\theta<1$ such that $ h(\vfi)\vfi^\theta \to \infty$ as $|x|\to \infty$. 
%for any $k<\bar k<\sigma$  such that $m_0\in L^1(\R^d, \lgx^{\bar k}\,dx)$, for some $K$ depending on $\alpha, \gamma,  \sigma, \lambda_0, \sigma_0,\sigma_1, \de, d, c_0, k, \bar k$.
\end{theorem}

\proof First of all we observe that, if $m$ is the  solution of \rife{FP-gen} (in the sense of Definition \ref{def-dual}), then $m_0\in \cM_\vfi$ implies $m(t) \in \cM_\vfi$ (reasoning exactly as   in \rife{preservk} for the function $\lgx^k$).
%
%it inherits some properties of $m_0$. In particular, due to conservation of positivity and mass\footnote{$m$ satisfies 
%$$
%\intd \vfi(t)\, m(t)\,dx+  \int_0^T\intd f  \,m\,dxdt = \intd \vfi(0,x)\, m_0(x)\, dx \qquad\forall \vfi\in C^0_b([0,T]\times \R^d)\,:\, \,\,\hbox{$-\partial_t \vfi + \cL^b[\vfi] =f$ in viscosity sense, $f\in C^0_b$}
%$$
%Choosing  $\vfi$ solution of $-\partial_t \vfi + \cL^b[\vfi] =0$ in viscosity sense, with $\vfi(0)=\xi \geq 0$, we deduce the preservation of positivity. Choosing $\vfi=1$, we deduce preservation of mass. Finally, by approximation, we see that $\vfi=\lg x\rg^k$ can be used as well, and we see also preservation of finite $k$-moments.
%}, if $m_0$ is a probability measure, then such is $m(t)$, for all $t>0$. Similarly, $m_0\in \cM_\vfi$ implies $m(t) \in \cM_\vfi$. 

Suppose now that $\vfi\in \Lambda_1$, or $\vfi\in \Lambda_2$  is such that $\varpi$ defined in \rife{decayvfi2} satisfies $\varpi\in L^1(0,\infty)$.   For any $s>0$, we can apply Theorem \ref{decayFP2} to $m(t+s)- m(t)$, which yields
\be\label{ts}
\begin{split}
\hbox{if $\vfi\in \Lambda_1$, then}\,\, \| m(t+s)-m(t)\|_{\cM_\vfi}  & \leq K \, e^{-\omega t} \, \|m(s)-m_0\|_{\cM_\vfi}\,,
\\
\hbox{if $\vfi\in \Lambda_2$, then}\,\, \| m(t+s)-m(t)\|_{\cM_{\vfi^\theta}}  & \leq K \, \varpi(t) \, \|m(s)-m_0\|_{\cM_\vfi}\,.
\end{split}
\ee
In particular, setting $m_n(x):= m(n,x)$, we have
\be\label{telesco}
m_n= m_0 + \sum_{k=0}^{n-1} [m_{k+1}-m_k]
\ee
where last series is strongly convergent in $\cM_\vfi$ (or in $\cM_{\vfi^\theta}$ if $\vfi\in \Lambda_2$), due to \rife{ts}. Hence there exists a measure $\bar m\in \cP(\R^d)$ such that $m_n\to \bar m$ in $\cM_\vfi$ (or in $\cM_{\vfi^\theta}$ if $\vfi\in \Lambda_2$). In fact,  $\bar m$ is itself a probability measure. Now we observe that, for any stationary $\xi$, we have 
$$
\intd \xi\, \, d(m_n-m_{n-1})   +  \int_{n-1}^{n}\intd   \cL^b[\xi]   \,dm(t)\, dt = 0
$$
which implies
\be\label{eqn}
\intd   \cL^b[\xi]   \,dm_n  =  \int_{n-1}^{n}\intd   \cL^b[\xi]   \,d(m_n-m(t))\, dt   - \intd \xi\, d(m_n-m_{n-1})  \,.
\ee
If $\vfi\in \Lambda_1$, by \rife{ts} we have  $\|m_n-m_{n-1}\|_{\cM_\vfi}\to 0$. The first integral in the right-hand side also goes to zero; indeed, again by \rife{ts} we have
$$
\left |\int_{n-1}^{n}\intd   \cL^b[\xi] \,d(m(t)-m_n)\, dt\right| \leq \| \cL^b[\xi]\|_{L^\infty(\vfi^{-1})} \int_{n-1}^n \|m(t)-m_n\|_{\cM_\vfi}dt \leq C\, \int_{n-1}^n e^{-\omega t}  dt \to 0\,.
$$
If $\vfi\in \Lambda_2$, the same argument applies in $\cM_{\vfi^\theta}$ and replacing $e^{-\omega t}$ with $\varpi(t)$, according to \rife{ts}.
Therefore, passing to the limit in \rife{eqn}, recalling the convergence of $m_n$, we conclude that $\bar m$ is a  stationary solution. Finally, we can now apply  Theorem \ref{decayFP2}  to $m-\bar m$, and we deduce that the whole sequence $m(t)$ converges to $\bar m$. Let us point out that, in case $\vfi\in \Lambda_2$, one should apply Theorem \ref{decayFP2} between the spaces $\cM_{\vfi^\theta}$ and $\cM_{\vfi^{\theta+\vep}}$, since a priori one only knows that $\bar m\in \cM_{\vfi^\theta}$ for any $\theta<1$. 
\qed

We point out that, even in the case   when the rate function $\varpi(t)$ in \rife{decayvfi2} is not integrable at infinity, the above argument produces the existence of a stationary solution $\bar m\in \cP(\R^d)$. Indeed, in the inequality  \rife{eqn} it is enough to use that $\varpi(t)\to 0$ as $t\to \infty$. By the compactness of $\cP(\R^d)$, the sequence $m_n$ will converge to a stationary solution. What is less obvious, in this case, is to show that the stationary solution enjoys
the further regularity $\bar m\in \cM_{\vfi^\theta}$ for any $\theta<1$. In fact, we don't know that the sum \rife{telesco} is  converging in $\cM_{\vfi^\theta}$, if $\varpi$ is no longer in $L^1(0,\infty)$. In specific cases, one should use extra arguments
to show first that $m(t)$ is bounded in $\cM_{\vfi^\theta}$ (uniformly in time). We mention that, for the case of power weights $\vfi= \lgx^k$ and fractional Laplace operator, one can use e.g. \cite[Thm 1.3]{Lafleche} to deduce that $\bar m\in L^1(\lgx^{\bar k})$ for all ${\bar k}\leq 1$; this allows for instance to apply to $m(t)-\bar m$  the decay estimate \rife{polym} for some $k\in (2-\gamma,1 \wedge \sigma)$.

%\newpage 

\section{Appendix: existence and uniqueness of solutions}

In this section we give a proof of the existence and uniqueness of solutions of the Fokker-Planck equation \rife{FP-gen}. Once more, this can be essentially deduced by duality from the well-posedness of viscosity solutions of the adjoint problem. This is why we start by giving a proof of the comparison principle, and uniqueness, of viscosity solutions.

\begin{proposition}\label{compa} Let $u_0\in C(\R^d)\cap \limk$, $f\in C(Q_T)\cap L^\infty(0,T;L^\infty(\lg x\rg^{-k}))$, for some $k<\sigma$. 
Let $\cL$ be defined by \rife{L} where $\Sigma$ satisfies  \rife{sigma},  the Levy measure in $\cI$ satisfies \rife{nu2}, and $\lambda_0+\lambda>0$.  Assume that $b\in C(Q_T)$ is such that $b(t,x)\cdot x$ is bounded below
and either of the following conditions hold:
\vskip0.4em
(i)  $\la_0>0$ or $\sigma\in (1,2)$,  and  \rife{disso} holds true whenever $|x-y|\leq 1$.
\vskip0.4em
(ii)  $\la_0=0$, $\sigma\in (0,1]$  and   \rife{disso2} holds true whenever $|x-y|\leq 1$.
\vskip0.4em
Then there exists a unique viscosity solution of  
\be\label{pbgenf}
\begin{cases}
\partial_t u + \cL[u]   + b(t,x) \cdot Du     =  f &  \text{ in } ( 0,T ) \times \R^d, \\ 
u ( 0,x ) = u_0 (x),\qquad  \,&  \text{ in }   \R^d
\end{cases}
\ee
such that  $|u(t,x)|\leq C (1+ |x|)^{m}$ for some $m<\sigma $.
\end{proposition}

\proof We prove first the uniqueness of solutions.  
Let us take some $\beta<\sigma$ such that $\beta \geq  k \vee m$. We first remark that (see the computations in Lemma \ref{lyap}) $\cL^b[ \lgx^\beta]\geq \beta \lgx^{\beta-2}(b\cdot x) - c_\beta (1+ \lgx^{\beta-1})$; hence there exists $K_0>0$ such that
$$
\lgx^\beta + \cL^b[ \lgx^\beta]\geq \frac12\lgx^\beta- K_0
$$
where $K_0$ only depends on $\sigma, d, \|\Sigma(x)\|_\infty , \Lambda$ and the lower bound of $b\cdot x$. We deduce that, choosing $ K$ large enough, the function
$W:= e^{  t} \lgx^\beta + K t+ K_1+ \frac1{T-t}$ is a supersolution, provided $ \lgx^\beta + K_1\geq |u_0|$. Let now $u,v$ be, respectively, a sub and super-solution of \rife{pbgenf}. By linearity, we have that $u_\vep:= u- \vep W$ and $v_\vep:=v+\vep W$ are still sub/supersolutions, and, if $\beta>m$, there exists a compact set $K_\vep\subset \R^d$ such that $u_\vep \leq 0 \leq v_\vep$ for $t\in (0,T)$, $x\in \R^d\setminus K_\vep$.  Of course, if we prove the comparison for $u_\vep, v_\vep$, then we get the comparison for $u,v$ by letting $\vep \to 0$.  Therefore, up to replacing $u,v$ with $u_\vep, v_\vep$,   we can assume  that $u\leq 0\leq v$   outside a compact set $K$. 
Now, suppose that $M:= \sup_{Q_T}(u-v)>0$, which means that $M$ is attained at some point in $(0,T)\times K$, with $t>0$. Doubling variables, we consider
$$
M_\vep:= \sup_{(x,y)}\,\,  \{u(t,x)-v(t,y)- \frac{|x-y|^2}{2\vep} - A_\vep  t^{\eta}\}
$$
for some $  \eta, A_\vep>0$ to be chosen, with $A_\vep \to 0$.  Since $M_\vep \geq M-A_\vep T^\eta$,  then $M_\vep$ is positive for $\vep$ small. Recalling that $u\leq 0\leq v$ outside $K$, this implies that  $M_\vep$ is assumed at points $(x_\vep, y_\vep)$ which lie in a compact set, which is actually an  $\vep$-neighborhood of $K$,
since $\frac{|x_\vep-y_\vep|^2}{2\vep}\leq \hat M= \sup u - \inf v  <\infty$.  In particular we have
\be\label{hatM}
 |x_\vep-y_\vep| \leq \sqrt{2\vep\, \hat M}
\ee
and by usual arguments, since $M_\vep \to M$ as $\vep \to 0$, we also get 
$\frac{|x_\vep-y_\vep|^2}\vep\to 0$. Using  Theorem \ref{test}, the Lipschitz character of $\Sigma(x)$ and standard  arguments in the doubling variable method, we get 
\be\label{aeps}
 \eta A_\vep \, t^{\eta-1}  + (b(t,x_\vep) - b(t,y_\vep) )\cdot \frac{(x_\vep-y_\vep)}\vep \leq f(t,x_\vep)-f(t,y_\vep)+ C\, \sigma_1\frac{|x_\vep-y_\vep|^2}\vep \,.
\ee 
At this stage we use the following regularity result, which can be obtained with similar arguments as used in Section 3 (see  the Appendix  in \cite{EJP}): using conditions \rife{disso}, or \rife{disso2}, for $|x-y|\leq 1$, every viscosity solution of \rife{pbgenf} satisfies
\be\label{lipt}
\begin{split}
 & | u(t,x)-u(t,y)|  \leq  C_T \frac{[\langle x\rangle^k+ \langle y\rangle^k]}{t^{\frac1{\sigma}}} \, |x-y| \,,\quad \forall x,y\,,\qquad \hbox{if $\sigma\in (1,2)$,}
 \\
 & | u(t,x)-u(t,y)|  \leq  C_T \frac{[\langle x\rangle^k+ \langle y\rangle^k]}{t^{\frac\theta\sigma}} \, |x-y|^\theta \,,\quad \forall x,y\,,\qquad \hbox{if $\sigma\in (0,1]$, $\theta<\sigma$.}
 \end{split}
 \ee
Therefore,  assuming that at least one  between $u,v$ is a solution, estimate \rife{hatM} can be improved accordingly. Recalling that the points $(x_\vep, y_\vep)$  lie in a compact set, we deduce that they satisfy, if $\sigma>1$,
 $$
 \frac{|x-y|^2}{2\vep} \leq \frac {C_{T,K} }{t^\frac1{\sigma }} |x-y| \leq C \frac {\sqrt\vep}{t^\frac1{\sigma }} 
\,. $$
 Notice that this constant $C_{T,K}$ only depends on the local bounds of $u,v$, and in particular is determined independently of $A_\vep$. Hence, \rife{aeps} implies   
\begin{align*}
  \eta A_\vep \, t^{\eta-1}    & \leq  \frac{2 C_{T,K}} {t^\frac1{\sigma}}\, |b(t,x_\vep) - b(t,y_\vep) | +  f(t,x_\vep)-f(t,y_\vep)+ C\, \sigma_1\frac{|x_\vep-y_\vep|^2}\vep
\\ & \leq C\, \omega( \vep) \left( \frac 1{t^\frac1{\sigma }}+ 1\right) + \frac C {t^\frac1{\sigma }}\, \sqrt\vep
\end{align*}
for some $C>0$ and some modulus of continuity $\omega(\cdot)$ depending    on  $b,f$. Choosing $\eta= 1-\frac1{\sigma}$, if  $A_\vep= L\, (\omega(  \vep)+ \sqrt \vep)$  we obtain a contradiction for $L$ sufficiently large. 

If $\sigma\in (0,1]$, we estimate
$$
 \frac{|x-y|^2}{2\vep} \leq \frac C {t^\frac\theta\sigma} |x-y|^\theta
$$
then we
use \rife{disso2} and we get
\begin{align*}
  \eta A_\vep \, t^{\eta-1}    & \leq c_0 \, |x_\vep-y_\vep|^{1-\sigma+ \de} \, \frac{|x_\vep-y_\vep|}\vep + f(t,x_\vep)-f(t,y_\vep) + C\, \sigma_1\frac{|x_\vep-y_\vep|^2}\vep \\
  & \leq  \frac C {t^\frac\theta\sigma}\, |x_\vep-y_\vep|^{\theta-\sigma+\de} +  f(t,x_\vep)-f(t,y_\vep)+ \sigma_1 \frac C {t^\frac\theta\sigma}\, |x_\vep-y_\vep|^\theta\,.
%\\ & \leq C\, \omega(\sqrt \vep) \left( \frac 1{t^\frac\kappa\sigma}+ 1\right)
\end{align*}
Choosing $\theta= \sigma-\frac\de 2$ and $\eta= 1-\frac\theta\sigma$, we conclude as before with a suitable choice of $A_\vep$. Having obtained a contradiction, we infer that 
$\sup_{Q_T}(u-v)\leq 0$, i.e. $u\leq v$. Reversing the role, we get the desired uniqueness.
The existence of solutions can be proved through the vanishing viscosity method (using \rife{lipt} to infer local  compactness, see also the proof of Theorem \ref{exiuniq} below).
\qed

The following lemma is a  classical statement in viscosity solutions theory to ensure equi-continuity in time from continuity in space (see \cite[Lemma 9.1]{BBL}). We give a proof for the reader's convenience, since the setting of our assumptions does not fall into previous versions.

\begin{lemma}\label{equic} Let $K\subset \R^d$ be a compact set, and let $u$ be a  viscosity solution of \rife{pbgenf}, where $b,f,u_0$ are continuous. 
 For any $\vep>0$, there exists $\de>0$ such that 
 $$
 \|u(t)-u_0\|_{L^\infty(K)} \leq \vep \qquad \forall t<\de\,,
 $$
 where $\de$ depends only on the bounds of $b, f, u$ and on the modulus of continuity of $u_0$, for $(t,x)$ in a  neighborhood of $\{0\}\times K$ (and of course on $\sigma, d, \|\Sigma(x)\|_\infty, \Lambda$ as well).
\end{lemma}

\proof For $\eta>0$ small enough, $y\in B_R$, $x\in B_{2R}$, let us consider the function
$$
\Phi(x):= u_0(y)+ C\, \zeta(|x-y|)+ kt+\eta
$$
where  $\zeta(s)=\cT_{R}(s^2)$; here $\cT_R(\cdot)$ denotes a $C^2$, bounded and  concave function such that $\cT_R(s)=s$ for $|s|<(3R)^2$. By concavity of $\cT_R$, we have
\begin{align*}
\cI(x,[\zeta(|x-y|)]) & \leq   \int_{|z|>1}  \{ \zeta ( |x+z-y| ) -\zeta( |x-y| )\} \nu(dz) \\
& \qquad +  \cT_{R}'(|x-y|^2) \int_{|z|\leq 1}   \{|x+z-y|^2 -|x-y|^2 -  2 ((x-y) \cdot  z ) \mathds{1}_{| z | \leq 1} \}\nu ( dz )
\\ & \leq M_R\,.
\end{align*}
In the local terms, we just observe that $\zeta(|x-y|)= |x-y|^2$; so, using the local bound of $b$, we have
$$
\cL_0^b(\zeta) \geq - M_R
$$
for a possibly different constant, only depending on $R$, and local bounds of $Q(x)$, $b$. Hence we have
$$
\Phi_t + \cL^b([\Phi])- f(x) \geq  k - C\, M_R - \|f\|_{L^\infty(B_{2R})}
$$
which implies that $\Phi$ is a supersolution for $x\in B_{2R}$, up to choosing $k \gtrsim CM_R+ \|f\|_{L^\infty(B_{2R})}$. If $x\in \partial B_{2R}$, we have $|x-y|>R$, and $\Phi\geq u$ provided  $C\, R^2+ u_0(y)\geq \|u(t)\|_{L^\infty(B_{2R})}$; in particular, this holds if 
$$
C\, R^2 \geq 2\|u\|_{L^\infty((0,T)\times B_{2R})}\,.
$$
At $t=0$, we have $u_0\leq \Phi$ if $u_0(x)-u_0(y)\leq C |x-y|^2 + \eta$. This holds if $\omega(|x-y|) -\eta \leq C |x-y|^2$, where $\omega(\cdot)$ is a modulus of continuity of $u_0$ in $B_{2R}$. Without loss of generality, we can assume that
$\omega(r) \geq L\, r$ near the origin, so there exists $ \eta_C:= \sup_{r\in (0,3R)} [\omega(r)-Cr^2]>0$ and $\eta_C\to 0$ as $C\to \infty$. Choosing  $\eta=\eta_C$, we have
the comparison between $u$ and $\Phi$ and we conclude that $u(t,x)\leq \Phi(t,x)$ for $x\in B_{2R}$. Choosing $x=y$ gives
$$
u(t,y)\leq u_0(y) + C\, M_R\, t +  \|f\|_{L^\infty(B_{2R})}\, t + \eta_C\,.
$$
Using $\eta_C\to 0$ as $C\to \infty$, we deduce that
$$
\| u(t)-u_0\|_{L^\infty(B_R)} \mathop{\to}^{t \to 0} 0 \,,
$$
and this limit is uniform for all solutions $u$ which are uniformly bounded in $B_{2R}$.
\qed

%\subsection{Existence and uniqueness of solutions to Fokker-Planck equations}

%We show, for completeness, the existence and uniqueness of solutions in the sense of Definition \ref{def-dual}. This is mostly relying on the stability of viscosity solutions of the dual problem. 

As a Corollary of the previous lemma, we infer a compactness principle for solutions of the Fokker-Planck equation.

\begin{lemma}\label{compactm}
Let $m\in C^0([0,T];\cP(\R^d))$ be a solution of \rife{FP-gen} in the sense of Definition \ref{def-dual}, where $m_0\in \cP(\R^d)\cap \cM_k(\R^d)$, and $\cL, b$ satisfy the assumptions of Theorem \ref{exiuniq}. Then,   for every $\vep>0$ there exists $\de>0$ such that 
\be\label{compactmd1}
d_1(m(t), m(s))\leq \vep \qquad \forall t,s\,,\,|t-s|\leq \de
\ee
where $\de$ only depends on $\|m_0\|_{\cM_k}$, on local uniform bounds of $b$, on lower   bound  of $b(t,x)\cdot x$,  as well as on the upper bounds $\|\Sigma(x)\|_\infty, \Lambda$.

In particular, any sequence  $\{m_n\} \subset  C^0([0,T]; \cP(\R^d))$ of solutions of \rife{FP-gen}, corresponding to $b_n, m_{0n}$,  is relatively compact in $C^0([0,T]; \cP(\R^d)) $ as soon as $b_n$ is locally uniformly bounded in $Q_T$, $b_n\cdot x$ is uniformly bounded below, and $m_{0n}$ is bounded in $\cM_k(\R^d)$.
\end{lemma}

\proof As in Proposition \ref{compa}, we take the function $W\simeq e^{(t-s)}\lgx^k+ K(t-s)$ so that  $-W_s + \cL^b[W]\geq 0$ in $(0,t)$. Any (smooth) truncation of this function remains a supersolution and can be used as test function in Definition \ref{def-dual}. We deduce
\be\label{preservk}
\intd \lgx^k\, dm(t) \leq C_t \intd \lgx^k\,  dm_0 \,,
\ee
so we have estimates of $m(t)$ in $\cM_k(\R^d)$ whenever $m_0 \in \cM_k(\R^d)$.

Now, let us estimate $d_1(m(t), m(s))$. To this purpose, take $\xi\,: \|\xi\|_{W^{1,\infty}(\R^d)}\leq 1$ and assume $t>s$. 
%
%
%We have
%\begin{align*}
%\forall \xi\,: \|\xi\|_{W^{1,\infty}(\R^d)}\leq 1\,,\qquad 
%\intd \xi\, d(m(t)-m(s))&  \leq   \int_{B_R^c} d(m(t)+m(s))  + \int_{B_R} \xi \, d(m(t)-m(s))
%\\  & \leq  \frac{2C_T}{R^k} |m_0|_{\cM_k} +  \int_{B_R} \xi \, d(m(t)-m(s))
%\end{align*}
In Definition \ref{def-dual}, we consider the solution $\vfi$ of (B) with $f=0$. Then we have
\begin{align*}
\intd \xi\, d(m(t)-m(s)) & = \intd (\vfi(s)-\xi)\, dm(s)  \leq   2 \int_{B_R^c} d m(s) + \int_{B_R} (\vfi(s)-\xi)\, dm(s)
\\
& \leq 2\frac{C_T}{R^k} \|m_0\|_{\cM_k} + \int_{B_R} (\vfi(s)-\xi)\, dm(s)
\end{align*}
where we used \rife{preservk}. 
By Lemma \ref{equic}, for fixed $R$ we have
$$
\int_{B_R} (\vfi(s)-\xi)\, dm(s) \leq \omega_R(t-s)
$$
for some quantity $\omega_R(r)$ converging to zero as $r\to 0$. Moreover $\omega_R(\cdot)$ only depends on $R$ and $\|b\|_{L^\infty((0,T)\times B_{2R})}$, and in particular it is uniform for all $\xi$ with $\|\xi\|_{W^{1,\infty}(\R^d)}\leq 1$. We conclude that
$$
d_1(m(t), m(s)) \leq \frac{C }{R^k} \|m_0\|_{\cM_k} + \omega_R(t-s)\,.
$$
Hence \rife{compactmd1} is proved. By Ascoli-Arzel\`a  theorem, this means that $m$ lies in  a relatively compact set of $C^0([0,T]; \cP(\R^d)) $ as $m_0$ is in a bounded set of $\cM_k$, and $b$ varies in  a set of functions which are  locally uniformly bounded with $b\cdot x$ uniformly bounded below.
\qed

Now we are ready to give the proof of Theorem \ref{exiuniq}.
\vskip1em

{\bf Proof of Theorem \ref{exiuniq}.} \quad Since any $m_0$ can be decomposed as $m_0^+-m_0^-$, using the Hahn decomposition, and since the problem is linear, it is enough to consider the case that $m_0$ is a positive measure. Without loss of generality, we assume $m_0\in \cP(\R^d)$. Indeed, it is immediate to check from Definition \ref{def-dual} that  $\into dm(t)=\into dm_0$, and $m_0\geq 0$ implies $m(t)\geq 0$; hence $m(t) \in \cP(\R^d)$. 

We construct a solution through vanishing viscosity approximation. So we consider the problem
\be\label{eps-FP}
\begin{cases}
\partial_t m -\vep \Delta m + \cL^*[m]   - \dive( b(t,x) m)     =  0 &  \text{ in } ( 0,T ) \times \R^d, \\ 
m ( 0,x ) = m_0 (x),\qquad  \,&  \text{ in }   \R^d
\end{cases}
\tag{$(FP)_\vep^b$}
\ee
where $\vep>0$. If $b$ is a smooth bounded vector field (e.g. $b\in C^2(Q_T)$ with $b, Db$ bounded), then for any smooth $m_0\in \cP(\R^d)$ there exists a smooth solution $m_\vep\in \cP(\R^d)$.
%\footnote{if $\sigma>1$, quote Olav-Espen. If $\sigma\in (0,1)$, do a fixed point as lower order perturbation of Laplacian. Other possibilities: truncate the singularity of the kernel, do a fixed point as lower order perturbation of Laplacian, then prove Lipschitz estimates for both $m$ and $Dm$, in which case $m$ is bounded $C^{1,\alpha}$ for all $\alpha$ and $W^{2,\infty}$; this latter implies a bound in $C^{0,\alpha}$ for $cI$, then it's done}. 
Of course $m_\vep$ is a weak as well as a duality solution, and satisfies
\be\label{eps-pb}
\intd m_\vep(t)\, \xi + 
\int_0^t\intd f \, m_\vep \, dx \,dt = \intd \vfi(0,x)\, dm_0 
\ee
 for every $ t\in (0,T)\,,\, \xi\in C_b(\R^d), \, f\in C_b(Q_t)$, where $\vfi $ is the (unique)  viscosity solution of 
\be\label{eps-HJ}
  \begin{cases}
-\partial_t \vfi  -\vep \Delta \vfi + \cL[\vfi] + b\cdot D\vfi      =  f &  \text{ in } ( 0,t ) \times \R^d, \\ 
\vfi(t,x) = \xi (x),\qquad  \,&  \text{ in }   \R^d\,.
\end{cases} 
\tag{$(H)_\vep^b$}
\ee
Now, we consider $b\in C^0_b$, and a sequence $\{b_n\}$ of smooth functions, uniformly bounded, which  converges to $b$ locally uniformly.
By standard results on the stability of viscosity solutions, the solutions $\vfi_\vep^n$ of $(H)^{b_n}_\vep$ are relatively compact and converge    locally uniformly in $[0,t]\times \R^d$ towards  a function $\vfi$ which is a  viscosity solution of \rife{eps-HJ}. This implies $\intd \vfi_\vep^n(0,x)\, dm_0   \to \intd \vfi(0,x)\, dm_0  $. As for the solutions $m_\vep^n$,   Lemma \ref{compactm} implies that the sequence $m_\vep^n$ is relatively compact in $C^0([0,T];\cP(\R^d))$. Hence, one can pass to the limit in the equality \rife{eps-pb}, so we build  a solution $m$ for $b\in C^0_b$.  In a similar way, we can get rid of the global boundedness of $b$; it is enough to take a sequence $\zeta_R$ of cut-off functions converging  to $1$ and to consider $b_R:= b\, \zeta_R$. This  is a sequence of bounded functions converging to $b$ locally uniformly, and such that $b_R\cdot x$ is uniformly bounded below. By Lemma \ref{compactm}, we have again that $m_R$, the corresponding solutions, are relatively compact in $C^0([0,T];\cP(\R^d))$. In parallel, the viscosity solutions $\vfi_R$ of $(H)^{b_R}_\vep$, are uniformly bounded and  relatively compact in $C^0([0,t]\times K)$ for all compact subsets $K$; hence, up to subsequences, $\vfi_R $ converges to a  viscosity solution of \rife{eps-HJ}. 

Now we can let $\vep \to 0$. The viscosity solutions of \rife{eps-HJ} are still uniformly bounded, and they are also relatively compact in $C^0([0,t]\times K)$ for all compact subsets $K$, by using estimates \rife{lipt} and Lemma \ref{equic} for the time equi-continuity. By  Lemma \ref{compactm}, the sequence $m_\vep$ is also relatively compact in $C^0([0,T];\cP(\R^d))$ and converges to some $m\in C^0([0,T];\cP(\R^d))$. We deduce that  $m$ satisfies Definition \ref{def-dual}, where $\vfi$ is the viscosity solution  of (B) which has been obtained in the vanishing viscosity limit. But this is the unique viscosity solution of problem (B), from Proposition \ref{compa}. Thus, the existence of $m$, satisfying Definition \ref{def-dual}, is proved. The uniqueness is straightforward, because problem (B) is solvable for all $\xi\in C^0_b$ (we just proved it, through vanishing viscosity limit).

Finally, since $m_0\in \cM_k$,  we have $m(t) \in \cM_k(\R^d)$ from \rife{preservk}. Then the formulation of $m$ extends to  test functions $f\in L^\infty((0,T);\limk),\xi\in \limk$ with an easy truncation argument, using dominated convergence theorem; it is enough, once more, to observe that the viscosity solutions of (B) are stable as bounded sequences $f_n,\xi_n$ approximate $f\in L^\infty((0,T);\limk),\xi\in \limk$, and in addition the limit solution is unique. 
%The preservation of mass and positivity is an easy exercise.
\qed

\end{document}